\documentclass[11pt,letterpaper]{amsart}
\usepackage{amsmath,amssymb,amsthm}
\numberwithin{equation}{section}
\usepackage{hyperref}
\hypersetup{colorlinks=true,linkcolor=blue,citecolor=blue,urlcolor=blue}
\newtheorem{theorem}{Theorem}[section]
\newtheorem{lemma}[theorem]{Lemma}
\newtheorem{proposition}[theorem]{Proposition}
\newtheorem{corollary}[theorem]{Corollary}
\theoremstyle{definition}
\newtheorem{definition}[theorem]{Definition}
\newtheorem{remark}[theorem]{Remark}
\newtheorem{example}[theorem]{Example}
\newtheorem*{pairingdefinition}{Definition}
\newtheorem{question}{Question}[section]

\title[Twisted convolution with conic support]{Twisted convolution with conic support and the parabolic algebra}
\author{Jie Qin}
\thanks{The author was supported by the National Natural Science Foundation of China (12501154).}
\address{School of Mathematics and Statistics,
Chongqing Technology and Business University,
Chongqing 400067, China}
\email{qinjie24520@163.com}
\subjclass[2020]{Primary 46F10 ; Secondary 47L75, 47B35}
\keywords{Twisted convolution, conic support, tempered distributions,
zero divisors, parabolic algebra, Toeplitz operators}
\begin{document}
\begin{abstract}
We prove that tempered distributions on \(\mathbb R^{2n}\) whose supports satisfy a common strict cone condition form a unital associative algebra without zero divisors under twisted convolution \(\star\). The supports may be noncompact, and the associated Weyl maps need not be bounded. As applications, we show that the parabolic algebra and the norm-closed algebra generated by positive translations and modulations have no zero divisors, confirming Levene's conjecture. We also obtain a zero-product theorem for Toeplitz operators on the one-dimensional Fock space with bounded symbols whose distributional Fourier transforms satisfy a common strict cone condition. The proof of the distribution theorem uses Schwartz smoothing between distributional Weyl maps to convert a hypothetical zero product of nonzero distributions into a smooth function \(h\not\equiv0\) satisfying the same cone condition and \(h\star h=0\). The cone bound implies the \(Q\)-support condition, meaning that some linear functional is bounded above on \(\operatorname{supp}h\) and has compact upper level sets there. Truncating \(h\) near the maximum of this functional and applying the support identity for compactly supported twisted squares rule out \(h\star h=0\). For the operator applications, uniformly bounded Weyl-conjugation averages identify twisted convolution with bounded operator composition under the common cone condition.
\end{abstract}
\maketitle

\section{Introduction}\label{sec:1}

On \(L^2(\mathbb R)\), consider the semigroups of translations and modulations in the positive directions. For \(\mu,\lambda\ge0\), the corresponding operators \(D_\mu\) and \(M_\lambda\) are defined by
$$
(D_\mu u)(x)=u(x-\mu),\quad
(M_\lambda u)(x)=e^{i\lambda x}u(x)
$$
for \(x\in\mathbb R\).

The parabolic algebra \(\mathcal A_p\) is the weak-operator closed algebra generated by these two semigroups. The weak-operator topology (WOT) is determined by the scalar maps \(A\mapsto\langle Au,v\rangle_2\) for \(u,v\in L^2(\mathbb R)\). Here \(\langle u,v\rangle_2=\int_{\mathbb R}u(x)\overline{v(x)}\,dx\), with associated norm \(\|\cdot\|_2\). We use the same convention on \(L^2(\mathbb R^m)\) for each integer \(m\ge1\).

An associative complex algebra has no zero divisors if a product of two nonzero elements is always nonzero. Levene conjectured in 2004 that \(\mathcal A_p\) has this property \cite[Conjecture 3.13, p. 37]{Levene2004}. Power asked whether \(\mathcal A_p\) is an integral domain \cite{Power2006}, and Kastis and Power later posed the question in terms of zero divisors \cite[Section 7, Question 4, pp. 585--586]{KastisPower2024}.
\begin{question}\label{Q1}
Does \(\mathcal A_p\) have zero divisors?
\end{question}

Let \(\mathcal A_{\mathrm{alg}}\) be the complex linear span of \(\{D_\mu M_\lambda:\mu,\lambda\ge0\}\). The relation $D_\mu M_\lambda=e^{-i\mu\lambda}M_\lambda D_\mu$ gives $$(D_\mu M_\lambda)(D_\nu M_\omega)=e^{i\lambda\nu}D_{\mu+\nu}M_{\lambda+\omega}$$ for \(\mu,\lambda,\nu,\omega\ge0\), so \(\mathcal A_{\mathrm{alg}}\) is an algebra.

Denote the closure of \(\mathcal A_{\mathrm{alg}}\) in the operator norm \(\|A\|=\sup_{\|u\|_2=1}\|Au\|_2\) by \(\mathcal A_{\mathrm n}\), and write \(\mathcal M=\mathrm{WOT}(\overline{\mathcal A_{\mathrm{alg}}})\) for its weak-operator closure. For each \(C\in\mathcal A_{\mathrm{alg}}\), weak-operator continuity of right multiplication gives \(\mathcal M C\subseteq\mathcal M\). Hence \(A\mathcal A_{\mathrm{alg}}\subseteq\mathcal M\) for \(A\in\mathcal M\), and continuity of left multiplication by \(A\) gives \(A\mathcal M\subseteq\mathcal M\). This identifies \(\mathcal M\) with the weak-operator closed algebra generated by the two semigroups. Since norm convergence implies weak-operator convergence, we have
\begin{equation}
\mathcal A_{\mathrm{alg}}\subseteq\mathcal A_{\mathrm n}\subseteq\mathcal A_p
=\mathrm{WOT}(\overline{\mathcal A_{\mathrm{alg}}}).
\label{E1}
\end{equation}
A related question asks whether a finite linear combination of distinct time-frequency shifts can annihilate a nonzero vector in \(L^2(\mathbb R)\) when its coefficients are not all zero \cite{HeilRamanathanTopiwala1996}. Linnell ruled this out for shifts in a discrete subgroup of phase space \cite{Linnell1999}. Heil and Speegle relate the question to zero divisors in the Heisenberg group ring \cite{HeilSpeegle2015}, and Enstad and van Velthoven deduce the discrete-subgroup result from the absence of zero divisors in twisted group rings of discrete locally indicable groups \cite[Proposition 3.4 and Theorem 4.1]{EnstadVanVelthoven2025}.

Question \ref{Q1} concerns the entire range of a second operator, since \(AB=0\) precisely when \(\operatorname{ran}B\subseteq\ker A\). To see why two nonzero finite sums in \(\mathcal A_{\mathrm{alg}}\) have nonzero product, order the parameter pairs with nonzero coefficients lexicographically. The sum of the two least indices occurs only once in the product, with a nonzero coefficient by the multiplication rule above. It follows that the product is nonzero, since the spreading distributions of distinct shifts are nonzero multiples of point masses at distinct points by \eqref{E24} and Lemma \ref{L3}.

The least-index argument breaks down in the norm closure, where the indices of the nonzero generalized Fourier coefficients need not have a least element. Kastis and Power expected zero divisors already in \(\mathcal A_{\mathrm n}\) \cite[Section 7, pp. 585--586]{KastisPower2024}. Kastis ruled out one possible source by proving that \(\mathcal A_{\mathrm n}\) has no nontrivial idempotents \cite[Corollary 3.6]{Kastis2019}. We approach the closure problem through a support property that holds throughout \(\mathcal A_p\). Under the Weyl correspondence, multiplication on this algebra becomes twisted convolution of spreading distributions supported in a common quadrant.

We use \(\mathcal S(\mathbb R^m)\) and \(\mathcal S'(\mathbb R^m)\) for the Schwartz space and tempered distributions, and \(\mathcal D(\mathbb R^m)=C_c^\infty(\mathbb R^m)\) and \(\mathcal D'(\mathbb R^m)\) for compactly supported smooth test functions and distributions. For \(T\in\mathcal S'(\mathbb R^m)\) and \(\varphi\in\mathcal S(\mathbb R^m)\), or \(T\in\mathcal D'(\mathbb R^m)\) and \(\varphi\in\mathcal D(\mathbb R^m)\), the distributional pairing is given by
\begin{equation}
\langle T,\varphi\rangle=T(\varphi).
\label{E43}
\end{equation}
This pairing is complex-linear in \(\varphi\). The continuous dense inclusion \(\mathcal D(\mathbb R^m)\subset\mathcal S(\mathbb R^m)\) identifies \(\mathcal S'(\mathbb R^m)\) with a subspace of \(\mathcal D'(\mathbb R^m)\) by restriction \cite[Lemma 7.1.8]{Hormander2003}. We identify each \(f\in L^1_{\mathrm{loc}}(\mathbb R^m)\) with the distribution given by \eqref{E44}. The topologies of the test spaces and the convergence conventions are specified in Subsection \ref{sec:notation}.

\begin{definition}[Distributional support]\label{def:distribution-support}
The support \(\operatorname{supp}T\) of \(T\in\mathcal D'(\mathbb R^m)\) is the complement of the largest open set on which it vanishes on all test functions supported in that set.
\end{definition}

For \(x,y\in\mathbb R^m\), write \(x\cdot y=\sum_{j=1}^m x_jy_j\) and let \(|x|\) be the Euclidean norm. Fix an integer \(n\ge1\), put \(d=2n\), and identify \(\mathbb R^d\) with \(\mathbb R^n\times\mathbb R^n\). We use the symplectic form defined for \(a=(s,t),b=(s',t')\in\mathbb R^n\times\mathbb R^n\) by
\begin{equation}
\sigma(a,b)=s\cdot t'-t\cdot s'.
\label{E55}
\end{equation}
The twisted convolution of \(f,g\in C_c^\infty(\mathbb R^d)\) is defined for \(z\in\mathbb R^d\) by
\begin{equation}
(f\star g)(z)=\int_{\mathbb R^d}
f(x)g(z-x)e^{i\sigma(x,z-x)}\,dx.
\label{eq:twisted-product}
\end{equation}

For a compactly supported test \(\varphi\), the function \((a,b)\mapsto\varphi(a+b)\) need not have compact support. In testing a twisted product of distributions, it suffices to localize near the pairs in the two supports for which \(a+b\in\operatorname{supp}\varphi\).

\begin{definition}[Properness of the addition map]\label{def:proper-addition}
For closed sets \(K_1,K_2\subset\mathbb R^d\), the addition map \((x,y)\mapsto x+y\) from \(K_1\times K_2\) to \(\mathbb R^d\) is \emph{proper} if
\begin{equation}
\{(x,y)\in K_1\times K_2:x+y\in L\}
\quad\hbox{is compact}
\label{E34}
\end{equation}
for every compact set \(L\subset\mathbb R^d\).
\end{definition}

For \(U,V\in\mathcal D'(\mathbb R^d)\), let \(U\otimes V\) be their tensor product as in \eqref{E36}. If addition is proper on \(\operatorname{supp}U\times\operatorname{supp}V\) in the sense of \eqref{E34}, it pushes \(e^{i\sigma}(U\otimes V)\) forward to the twisted product \(U\star V\). This is the product of \cite[Definition 1, p. 43]{LoupiasMiracleSole1967} with their symplectic form taken to be \(-\sigma\). Its local formula \eqref{E23} uses the convolution construction of \cite[Section 4.2, p. 104]{Hormander2003} with the smooth multiplier \(e^{i\sigma}\) inserted by \eqref{E35}. For tempered factors, the following cone condition ensures properness and a unique tempered extension of this local product (Proposition \ref{prop:tempered-product}).

\begin{definition}[Cone support]\label{def:cone-algebra}
Fix \(e\in\mathbb R^{2n}\) and \(\eta>0\). Let \(\mathcal C_{e,\eta}\) be the class of all \(U\in\mathcal S'(\mathbb R^{2n})\) for which there is \(M_U\ge0\) such that
\begin{equation}
e\cdot a\ge\eta|a|-M_U
\label{eq:cone}
\end{equation}
for every \(a\in\operatorname{supp}U\).

\end{definition}

\begin{theorem}\label{thm:distribution}

For every \(e\in\mathbb R^{2n}\) and \(\eta>0\), twisted convolution makes \(\mathcal C_{e,\eta}\) a unital associative complex algebra with no zero divisors. In particular, we have
$$
U,V\in\mathcal C_{e,\eta},\quad U\star V=0
\quad\Longrightarrow\quad U=0\ \text{or}\ V=0.
$$

\end{theorem}

When \(e=0\), \(\mathcal C_{0,\eta}\) consists of all compactly supported distributions. The theorem also covers noncompact supports, with no boundedness assumption on the associated Weyl maps. Examples arise from tempered distributions supported in translates of a fixed closed convex pointed cone, as considered by Schwartz \cite[pp. 32--33]{Schwartz1951}. For ordinary convolution, Meril and Struppa studied convolution equations for distributions supported in a cone \cite{MerilStruppa1987}, and Schmidt and Dierolf constructed a predual for the algebra allowing translates of the cone \cite{SchmidtDierolf1994}. Bargetz gives the translated-cone formulation for cones with nonempty interior \cite[Proposition 7, p. 62]{Bargetz2011}.

To verify that one pair \(e,\eta\) works for all translates of a closed convex cone \(C\subset\mathbb R^d\) with \(C\cap(-C)=\{0\}\), write \(C^*=\{y\in\mathbb R^d:y\cdot x\ge0\text{ for all }x\in C\}\). The polar in \cite[Corollary 14.6.1]{Rockafellar1970} is \(-C^*\). Since the lineality space \(C\cap(-C)\) is zero, the corollary gives \(\dim C^*=d\), so this convex cone has nonempty interior. Choose \(e\) and \(\eta>0\) so that \(\overline B(e,\eta)\subset C^*\). Testing \(e-\eta x/|x|\in C^*\) against \(x\in C\setminus\{0\}\) yields \(e\cdot x\ge\eta|x|\) on \(C\). Each translate \(b+C\) then satisfies the same bound up to a constant, since for \(a\in b+C\) we have
$$
e\cdot a\ge e\cdot b+\eta|a-b|
\ge\eta|a|-(\eta|b|+|e\cdot b|).
$$

For functions on the line with supports bounded below, Titchmarsh's convolution theorem determines the lower endpoint of the support of the convolution from those of its factors \cite[Theorem VII, p. 286]{Titchmarsh1926}. Domar extended this result to certain weighted spaces on the line \cite{Domar1983}. In several variables, Lions proved that \(S*T\ne0\) for nonzero distributions supported in translates of a common closed convex pointed cone \cite[Corollary 1, p. 1531]{Lions1951} and studied support identities under decay assumptions in a second note \cite{Lions1951b}. For nonzero compactly supported distributions, the convex hull of the support of \(S*T\) is the sum of the convex hulls of the two supports \cite[Theorem 4.2, p. 292]{Hormander1963}.

For group convolution, Garncarek excluded zero divisors among compactly supported finite measures on connected, simply connected supersolvable Lie groups \cite{Garncarek2013} and gave a direct proof for compactly supported continuous functions on the Heisenberg group \cite[Section 3]{Garncarek2014}. Passing from Heisenberg convolution to \eqref{eq:twisted-product} fixes a central Fourier frequency, so nonvanishing on the group alone does not settle the question at that frequency. Properness of addition suffices to define the twisted product in \eqref{E23} but still allows zero divisors, as Example \ref{Eg1} shows.

Our proof of Theorem \ref{thm:distribution} starts with a counterpart of the classical convex-hull identity for a twisted square. Write \(\operatorname{conv}K\) for the convex hull of \(K\subset\mathbb R^d\) and \(2K=\{2x:x\in K\}\). Lemma \ref{lem:compact-square} establishes
$$
\operatorname{conv}\operatorname{supp}(c\star c)
=2\operatorname{conv}\operatorname{supp}c
$$
for every \(c\in C_c^\infty(\mathbb R^d)\) with \(c\not\equiv0\).
For a nonzero smooth \(h\), condition \((Q)\) of Definition \ref{def:upper-caps} provides a linear functional \(\ell\) that is bounded above on \(\operatorname{supp}h\) and has compact upper level sets there. Proposition \ref{prop:smooth-square} applies the identity to a truncation near the maximum of \(\ell\) on this support to deduce \(h\star h\ne0\).

For nonzero \(U,V\in\mathcal C_{e,\eta}\), Weyl smoothing and approximation of a rank-one map yield \(g\in C_c^\infty(\mathbb R^d)\) for which \(h=V\star g\star U\) is a nonzero smooth function in \(\mathcal C_{e,\eta}\). Its cone bound implies \((Q)\) with \(\ell(z)=-e\cdot z\) by \eqref{E80}, so its smooth twisted square is nonzero. Proposition \ref{prop:tempered-product} and \eqref{E86} identify that smooth square, on compactly supported tests, with the distributional square \(h\star h\). If \(U\star V=0\), associativity would make this very square vanish, since \(U\star V\) occurs between the two copies of \(g\).

For bounded operators \(A,B\) on \(L^2(\mathbb R^n)\), let \(U_A,U_B\) denote their spreading distributions, so \(W(U_A)=A\) and \(W(U_B)=B\) on Schwartz functions (see \cite{DoerflerTorresani2010}). Under the common cone condition, Proposition \ref{prop:bounded-composition} uses uniformly bounded Weyl-conjugation averages to identify \(AB\) with \(W(U_A\star U_B)\) on this space. Injectivity of the Weyl map transfers the zero-product conclusion to \(A,B\) in Corollary \ref{cor:bounded-zero-product}. For the parabolic algebra, the two invariant nests place the spreading distributions in a common quadrant, as shown in Subsection \ref{sec:6}.

\begin{corollary}[Answer to Question \ref{Q1}]\label{cor:parabolic}

The parabolic algebra \(\mathcal A_p\) has no zero divisors. More precisely, we have
$$
A,B\in\mathcal A_p,\quad AB=0
\quad\Longrightarrow\quad A=0\ \text{or}\ B=0.
$$
This confirms Levene's conjecture, answers Power's question \cite{Power2006} affirmatively, and gives a negative answer to Question \ref{Q1}. By \eqref{E1}, the same conclusion holds for \(\mathcal A_{\mathrm n}\) and \(\mathcal A_{\mathrm{alg}}\).

\end{corollary}

Our second application concerns the one-dimensional bounded-symbol case of Bauer and Le's zero-product problem \cite[Section 6, Question A]{BauerLe2011}. Following \cite[Section 2.1 and Corollary 2.3]{Zhu2012}, identify \(z=x+iy\in\mathbb C\) with \((x,y)\in\mathbb R^2\), let \(dz\) denote planar Lebesgue measure, and put \(d\gamma(z)=\pi^{-1}e^{-|z|^2}\,dz\). The Fock space in this normalization is
$$
\mathcal F^2(\mathbb C)=
\{h\in L^2(\mathbb C,d\gamma):h\text{ is entire}\}.
$$
This is a closed subspace of \(L^2(\mathbb C,d\gamma)\), with orthogonal projection \(P:L^2(\mathbb C,d\gamma)\to\mathcal F^2(\mathbb C)\). For an essentially bounded symbol \(f\in L^\infty(\mathbb C)\), the Toeplitz operator is defined by \(T_fh=P(fh)\) for \(h\in\mathcal F^2(\mathbb C)\), and \(\|T_f\|\le\|f\|_\infty\).

We regard each bounded symbol as a tempered distribution by the \(L^\infty\) case of \eqref{E44} and take its Fourier transform in that sense. For \(f\in\mathcal S(\mathbb R^m)\) and \(\xi\in\mathbb R^m\), our Fourier-transform convention is
\begin{equation}
\widehat f(\xi)=\int_{\mathbb R^m}f(x)e^{-ix\cdot\xi}\,dx.
\label{E7}
\end{equation}
With the pairing convention \eqref{E43}, the Fourier transform of \(T\in\mathcal S'(\mathbb R^m)\) is defined by
\begin{equation}
\langle\widehat T,\varphi\rangle=\langle T,\widehat\varphi\rangle
\label{E15}
\end{equation}
for \(\varphi\in\mathcal S(\mathbb R^m)\). Under the Bargmann transform, \(T_f\) corresponds to a bounded operator on \(L^2(\mathbb R)\) with spreading distribution \(U_f\). By \eqref{eq:fock-support}, a fixed invertible linear map carries \(\operatorname{supp}\widehat f\) onto \(\operatorname{supp}U_f\), so a common Fourier support cone transfers to the two spreading distributions. Formula \eqref{E2} then translates vanishing of either distribution back to almost-everywhere vanishing of its symbol.

\begin{corollary}\label{cor:fock-zero-product}

Let \(f,g\in L^\infty(\mathbb C)\). Suppose that their distributional Fourier transforms satisfy \(\widehat f,\widehat g\in\mathcal C_{e,\eta}\) for some \(e\in\mathbb R^2\) and \(\eta>0\). If \(T_fT_g=0\) on \(\mathcal F^2(\mathbb C)\), then \(f=0\) or \(g=0\) almost everywhere.

\end{corollary}

\section{Notation and support theorems for smooth twisted squares}\label{sec:2}

\subsection{Notation}\label{sec:notation}

\makeatletter
\begingroup
\phantomsection
\let\@currentlabel\@empty
\def\@currentlabelname{Test spaces and pairings}
\label{def:test-spaces}
\endgroup
\makeatother
For background on test spaces, distributions, Fourier transforms, and the Schwartz kernel theorem, see \cite{Hormander2003}. Write \(\mathbb N_0=\{0,1,\ldots\}\). For an integer \(m\ge1\), a multiindex is \(\alpha=(\alpha_1,\ldots,\alpha_m)\in\mathbb N_0^m\), with \(|\alpha|=\alpha_1+\cdots+\alpha_m\) and \(\partial^\alpha=\partial_{x_1}^{\alpha_1}\cdots\partial_{x_m}^{\alpha_m}\), where \(x=(x_1,\ldots,x_m)\in\mathbb R^m\). Put \(\langle x\rangle=(1+|x|^2)^{1/2}\). Recall that the Schwartz space \(\mathcal S(\mathbb R^m)\) consists of the smooth complex-valued functions \(f\) satisfying
\begin{equation}
p_{N,k}(f)=\max_{|\alpha|\le k}\sup_{x\in\mathbb R^m}
\langle x\rangle^N|\partial^\alpha f(x)|<\infty
\label{E51}
\end{equation}
for every \(N,k\in\mathbb N_0\).
These seminorms define the Schwartz topology. When an estimate involves different sets of variables, the superscript in \(p_{N,k}^{(x)}\) or \(p_{N,k}^{(x,y)}\) specifies the variables used in \eqref{E51}. The space \(\mathcal S'(\mathbb R^m)\) of tempered distributions is its continuous complex-linear dual, with the pairing \eqref{E43}. Let \(T\in\mathcal S'(\mathbb R^m)\) and \(\varphi\in\mathcal S(\mathbb R^m)\). By \cite[Definition 7.1.9 and Theorem 7.1.10]{Hormander2003} and Fourier inversion in the convention \eqref{E7}, $$\widehat{\bigl[\widehat\varphi(-\cdot)\bigr]}=(2\pi)^m\varphi.$$  Applying \eqref{E15} to the test function \(\widehat\varphi(-\cdot)\) yields
\begin{equation}
\langle T,\varphi\rangle
=(2\pi)^{-m}\langle\widehat T,\widehat\varphi(-\cdot)\rangle.
\label{E13}
\end{equation}

The space \(\mathcal D(\mathbb R^m)=C_c^\infty(\mathbb R^m)\) carries its usual test-function topology. Its continuous complex-linear dual \(\mathcal D'(\mathbb R^m)\) is the space of distributions. A sequence \(\varphi_j\) converges to \(\varphi\) in \(\mathcal D(\mathbb R^m)\) precisely when the supports of \(\varphi\) and all \(\varphi_j\) lie in a common compact set and \(\partial^\alpha\varphi_j\to\partial^\alpha\varphi\) uniformly for every multiindex \(\alpha\) \cite[Corollary 1 to Proposition 14.6]{Treves1967}. We identify a locally integrable function \(f\in L^1_{\mathrm{loc}}(\mathbb R^m)\) with the distribution given by
\begin{equation}
\langle f,\varphi\rangle=\int_{\mathbb R^m}f(x)\varphi(x)\,dx
\label{E44}
\end{equation}
for \(\varphi\in\mathcal D(\mathbb R^m)\). All unweighted spaces \(L^p(\mathbb R^m)\), \(1\le p\le\infty\), are taken with respect to Lebesgue measure, with norm \(\|\cdot\|_p\). If \(f\in L^p(\mathbb R^m)\) for some \(1\le p\le\infty\), H\"older's inequality shows that \eqref{E44} also holds for \(\varphi\in\mathcal S(\mathbb R^m)\) and defines an element of \(\mathcal S'(\mathbb R^m)\). Under this identification, the inclusions \cite[Sections 2.1 and 7.1]{Hormander2003} and the density theorem \cite[Theorem 15.4, p. 160]{Treves1967} give
\begin{equation}
\begin{aligned}
C_c^\infty(\mathbb R^m)\subset\mathcal S(\mathbb R^m)\subset\mathcal S'(\mathbb R^m) \quad\text{and}\quad
\overline{C_c^\infty(\mathbb R^m)}^{\,\mathcal S(\mathbb R^m)}=\mathcal S(\mathbb R^m).
\end{aligned}
\label{E18}
\end{equation}
By Definition \ref{def:distribution-support}, the distributional support of a continuous function \(h\in C(\mathbb R^m)\) is given by
\begin{equation}
\operatorname{supp}h=\overline{\{x\in\mathbb R^m:h(x)\ne0\}}.
\label{E37}
\end{equation}

Multiplication by \(a\in C^\infty(\mathbb R^m)\) assigns to each \(T\in\mathcal D'(\mathbb R^m)\) the distribution \(aT\in\mathcal D'(\mathbb R^m)\) defined by
\begin{equation}
\langle aT,\varphi\rangle=\langle T,a\varphi\rangle
\label{E35}
\end{equation}
for \(\varphi\in\mathcal D(\mathbb R^m)\). If \(T\in\mathcal S'(\mathbb R^m)\) and the map \(\varphi\mapsto a\varphi\) is continuous from \(\mathcal S(\mathbb R^m)\) to itself, the same formula for \(\varphi\in\mathcal S(\mathbb R^m)\) defines \(aT\in\mathcal S'(\mathbb R^m)\). Suppose that, for each \(\alpha\in\mathbb N_0^m\), there are \(C_\alpha>0\) and \(N_\alpha\in\mathbb N_0\) such that \(|\partial^\alpha a(x)|\le C_\alpha\langle x\rangle^{N_\alpha}\) for all \(x\in\mathbb R^m\). For fixed \(k\), set \(M_k=\max_{|\alpha|\le k}N_\alpha\). Leibniz' rule bounds \(p_{N,k}(a\varphi)\) by \(C_kp_{N+M_k,k}(\varphi)\), where \(C_k\) depends only on \(a,k\) and the dimension. This proves continuity of the multiplier
\begin{equation}
\mathcal S(\mathbb R^m)\ni\varphi\longmapsto a\varphi
\in\mathcal S(\mathbb R^m).
\label{E89}
\end{equation}

In the pairing \eqref{E43}, we indicate the variable \(x\in\mathbb R^m\) by writing
$$
\langle U(x),\varphi(x)\rangle
=\langle U,x\mapsto\varphi(x)\rangle
=U(\varphi)
$$
for \(U\in\mathcal D'(\mathbb R^m)\) and \(\varphi\in\mathcal D(\mathbb R^m)\). The same notation applies to \(U\in\mathcal S'(\mathbb R^m)\) and \(\varphi\in\mathcal S(\mathbb R^m)\).
For integers \(r,s\ge1\), \(U\in\mathcal D'(\mathbb R^r)\), and \(V\in\mathcal D'(\mathbb R^s)\), the tensor product \(U\otimes V\in\mathcal D'(\mathbb R^{r+s})\) is defined by the iterated pairings \cite[Theorem 5.1.1]{Hormander2003}
\begin{equation}
\begin{aligned}
\langle U\otimes V,\varphi\rangle
&=\big\langle U(x),\langle V(y),\varphi(x,y)\rangle\big\rangle=\big\langle V(y),\langle U(x),\varphi(x,y)\rangle\big\rangle
\end{aligned}
\label{E36}
\end{equation}
for \(\varphi\in\mathcal D(\mathbb R^{r+s})\), with \(x\in\mathbb R^r\) and \(y\in\mathbb R^s\). We shall use the support identity \cite[(5.1.3), p. 128]{Hormander2003}
\begin{equation*}
\operatorname{supp}(U\otimes V)=\operatorname{supp}U\times\operatorname{supp}V.
\tag{\getrefnumber{E36}a}\label{E121}
\end{equation*}
If \(U\in\mathcal S'(\mathbb R^r)\) and \(V\in\mathcal S'(\mathbb R^s)\), the continuous bilinear form \((f,g)\mapsto\langle U,f\rangle\langle V,g\rangle\) extends to a continuous functional on \(\mathcal S(\mathbb R^{r+s})\) by the completed tensor product identification in \cite[Theorem 51.6, (51.3)]{Treves1967}. For \(\varphi\in\mathcal S(\mathbb R^{r+s})\), put \(F_\varphi(x)=\langle V,\varphi(x,\cdot)\rangle\). Choose \(C_V>0\) and \(M,j\in\mathbb N_0\) so that \(|\langle V,\psi\rangle|\le C_Vp_{M,j}^{(y)}(\psi)\) for \(\psi\in\mathcal S(\mathbb R^s)\). Taylor's formula in each \(x\)-coordinate, with remainder estimated in \(p_{M,j}^{(y)}\), gives \(\partial_x^\alpha F_\varphi(x)=\langle V,\partial_x^\alpha\varphi(x,\cdot)\rangle\) for every multi-index \(\alpha\). Applying the bound for \(V\) to this slice supplies the \(y\)-weight and the \(j\) additional derivatives in the next estimate. Since \(\langle x\rangle^N\langle y\rangle^M\le\langle(x,y)\rangle^{N+M}\), \eqref{E51} gives
\begin{equation*}
p_{N,k}^{(x)}(F_\varphi)
\le C_Vp_{N+M,k+j}^{(x,y)}(\varphi)
\qquad(N,k\in\mathbb N_0).
\tag{\getrefnumber{E36}b}\label{E127}
\end{equation*}
Thus \(U\) can act on \(F_\varphi\), and the iterated pairing is continuous in \(\varphi\). The same estimate with the variables exchanged applies in the reverse order. Both pairings agree with the tensor product on the dense span of elementary tensors, proving \eqref{E36} for Schwartz tests.
With \(U\in\mathcal S'(\mathbb R^r)\), \(V\in\mathcal S'(\mathbb R^s)\), \(\chi_1\in C_c^\infty(\mathbb R^r)\), and \(\chi_2\in C_c^\infty(\mathbb R^s)\), equations \eqref{E35} and \eqref{E36} give
\begin{equation}
\big\langle(\chi_1U)\otimes(\chi_2V),\psi\big\rangle
=\big\langle U\otimes V,\chi_1(x)\chi_2(y)\psi(x,y)\big\rangle
\label{E81}
\end{equation}
for every \(\psi\in\mathcal S(\mathbb R^{r+s})\).

The tensor product of functions \(f:\mathbb R^r\to\mathbb C\) and \(g:\mathbb R^s\to\mathbb C\) is the function \(f\otimes g:\mathbb R^r\times\mathbb R^s\to\mathbb C\) defined by
\begin{equation}
(f\otimes g)(x,y)=f(x)g(y)
\label{E17}
\end{equation}
for \(x\in\mathbb R^r\) and \(y\in\mathbb R^s\). If \(f\) and \(g\) are locally integrable, \eqref{E44} and Fubini's theorem identify \(f\otimes g\) with the tensor product of their associated distributions. Tensor products of more than two functions are defined by iteration. For \(f\in\mathcal S(\mathbb R^r)\) and \(g\in\mathcal S(\mathbb R^s)\), the derivatives in \eqref{E17} act on separate factors. Since \(\langle(x,y)\rangle\le\langle x\rangle\langle y\rangle\), the seminorms \eqref{E51} satisfy
\begin{equation*}
p_{N,k}^{(x,y)}(f\otimes g)\le p_{N,k}^{(x)}(f)\,p_{N,k}^{(y)}(g)
\qquad(N,k\in\mathbb N_0).
\tag{\getrefnumber{E17}a}\label{E118}
\end{equation*}

The point mass at \(p\in\mathbb R^d\) is the distribution \(\delta_p\in\mathcal S'(\mathbb R^d)\) defined by \(\delta_p(\varphi)=\varphi(p)\) for \(\varphi\in\mathcal S(\mathbb R^d)\). Weak convergence \(T_j\to T\) in \(\mathcal S'\), with \(T,T_j\in\mathcal S'(\mathbb R^m)\) and \(j=1,2,\ldots\), means
$$
\langle T_j,\varphi\rangle\longrightarrow\langle T,\varphi\rangle
$$
as \(j\to\infty\) for every \(\varphi\in\mathcal S(\mathbb R^m)\).
Weak convergence in \(\mathcal D'\) tests only functions in \(C_c^\infty\subset\mathcal S\), so it follows from weak \(\mathcal S'\) convergence.

\subsection{Compactly supported twisted squares}

To locate the support of a smooth twisted square, we first concentrate both factors at the same exposed point of the support. Normalization on the Fourier side permits this localization even when the function vanishes to arbitrary order at the supporting boundary.

\begin{lemma}\label{lem:localization} 
Let \(c\in C_c^\infty(\mathbb R^d)\) with \(c\not\equiv0\). Suppose that a real-linear functional \(\ell:\mathbb R^d\to\mathbb R\) attains its maximum \(m\) on \(\operatorname{supp}c\) at exactly one point \(p\). There are sequences \((R_\nu)_{\nu\ge1}\subset(0,\infty)\) and \((\tau_\nu)_{\nu\ge1}\subset\mathbb R^d\), with \(R_\nu\to\infty\) as \(\nu\to\infty\), such that the functions \(a_\nu:\mathbb R^d\to\mathbb C\) defined by
$$
a_\nu(x)=
\frac{e^{R_\nu\ell(x)-i\tau_\nu\cdot x}c(x)}
{\widehat{\bigl[e^{R_\nu\ell}c\bigr]}(\tau_\nu)}
$$            
satisfy
\postdisplaypenalty=10000\relax
\begin{equation}
a_\nu\longrightarrow\delta_p,
\qquad a_\nu\otimes a_\nu\longrightarrow\delta_p\otimes\delta_p
\label{eq:localization}
\end{equation}
weakly in \(\mathcal S'(\mathbb R^d)\) and \(\mathcal S'(\mathbb R^{2d})\), respectively, and
\begin{equation}
R_\nu^{-1}\log\left|\widehat{\bigl[e^{R_\nu\ell}c\bigr]}(\tau_\nu)\right|
\longrightarrow m. \label{eq:fourier-growth}
\end{equation}
All \(a_\nu\) are supported in the fixed compact set \(\operatorname{supp}c\).

\end{lemma}

\begin{proof}

For \(R>0\), set \(u_R=e^{R\ell}c\), \(f_R=\widehat{u_R}\), and \(M_R=\|f_R\|_\infty\). We choose the denominator in \eqref{E14} at a maximizing frequency so that the Fourier transform of the normalized function equals one at zero and has modulus at most one everywhere. On each nonempty compact subset of \(\operatorname{supp}c\setminus\{p\}\), the maximum of \(\ell\) is strictly less than \(m\). To retain this exponential gap after normalization, we first prove that \(R^{-1}\log M_R\to m\). For each \(\varepsilon>0\), the open set \(\{\ell>m-\varepsilon\}\) contains \(p\) and hence meets \(\{c\ne0\}\) by \eqref{E37}. Continuity of \(c\) ensures that the constant \(C_{1,\varepsilon}=\int_{\{\ell>m-\varepsilon\}}|c(x)|^2\,dx\) is positive and independent of \(R\). Plancherel's theorem \cite[Theorem 7.1.11]{Hormander2003} in the normalization \eqref{E7} and the bound \(e^{2R\ell(x)}\ge e^{2R(m-\varepsilon)}\) on \(\{\ell>m-\varepsilon\}\) give
\begin{equation}
\|f_R\|_2^2=(2\pi)^d\|u_R\|_2^2
\ge(2\pi)^d C_{1,\varepsilon} e^{2R(m-\varepsilon)}.
\label{E9}
\end{equation}
In particular, \(M_R>0\) for every \(R>0\). To obtain a quantitative lower bound for \(M_R\) from \eqref{E9}, we estimate the Fourier tail using the gradient of \(u_R\). By the finite-dimensional Riesz representation theorem \cite[Theorem 6.42]{Axler2024}, there is a fixed vector \(b_\ell\in\mathbb R^d\) such that \(\ell(x)=b_\ell\cdot x\). With \(\nabla\) denoting the gradient in the real variables, the product rule reads
$$
\nabla u_R=e^{R\ell}(\nabla c+R b_\ell c).
$$
Since \(\operatorname{supp}\nabla c\subset\operatorname{supp}c\) and \(\ell\le m\) on \(\operatorname{supp}c\), the \(L^2\) triangle inequality gives \(\|\nabla u_R\|_2\le e^{Rm}(\|\nabla c\|_2+R|b_\ell|\|c\|_2)\), and hence
\begin{equation}
\|\nabla u_R\|_2^2\le C_2(1+R)^2e^{2Rm},
\label{E46}
\end{equation}
where \(C_2>0\) depends only on \(c,\ell,d\).
Since \(u_R\in C_c^\infty(\mathbb R^d)\), we may integrate by parts in \eqref{E7} to obtain
\begin{equation}
\widehat{\bigl[\partial_{x_j}u_R\bigr]}(\xi)=i\xi_j f_R(\xi)
\label{E8}
\end{equation}
for \(\xi\in\mathbb R^d\) and \(j=1,\ldots,d\). Summing the squared \(L^2\) norms in \eqref{E8} and applying Plancherel's theorem to each derivative, we can express the second Fourier moment as
\begin{equation}
\begin{aligned}
\int_{\mathbb R^d}|\xi|^2|f_R(\xi)|^2\,d\xi
&=\sum_{j=1}^d\|\widehat{\bigl[\partial_{x_j}u_R\bigr]}\|_2^2=(2\pi)^d\|\nabla u_R\|_2^2.
\end{aligned}
\label{E11}
\end{equation}

We choose \(L_R=e^{2\varepsilon R}\) so that the Fourier tail is small compared with the lower bound in \eqref{E9}. On \(\{|\xi|>L_R\}\) we have \(1\le L_R^{-2}|\xi|^2\). Expressing the resulting Fourier moment by \eqref{E11} and applying \eqref{E46} with \(L_R^{-2}=e^{-4\varepsilon R}\), we obtain
\begin{equation}
\begin{aligned}
\int_{|\xi|>L_R}|f_R(\xi)|^2\,d\xi
&\le L_R^{-2}\int_{\mathbb R^d}|\xi|^2|f_R(\xi)|^2\,d\xi=(2\pi)^dL_R^{-2}\|\nabla u_R\|_2^2\\
&\le C_3(1+R)^2e^{2R(m-2\varepsilon)},
\end{aligned}
\label{E10}
\end{equation}
where \(C_3=(2\pi)^dC_2\). For fixed \(\varepsilon>0\), dividing \eqref{E10} by the lower bound for \(\|f_R\|_2^2\) in \eqref{E9} gives
$$
\frac{\displaystyle\int_{|\xi|>L_R}|f_R(\xi)|^2\,d\xi}{\|f_R\|_2^2}
\le\frac{C_3}{(2\pi)^dC_{1,\varepsilon}}
(1+R)^2e^{-2\varepsilon R}\longrightarrow0
$$
as \(R\to\infty\).

For each fixed \(\varepsilon>0\), the preceding limit shows that the integral over \(\{|\xi|\le L_R\}\) is at least \(\frac12\|f_R\|_2^2\) for all sufficiently large \(R\). Set \(C_{4,\varepsilon}=\frac12(2\pi)^dC_{1,\varepsilon}>0\), and let \(C_d\) be the volume of the unit ball in \(\mathbb R^d\). Using \eqref{E9}, the bound \(|f_R|\le M_R\), and the ball volume \(C_dL_R^d=C_de^{2d\varepsilon R}\), we obtain
\begin{equation}
C_{4,\varepsilon} e^{2R(m-\varepsilon)}
\le\frac12\|f_R\|_2^2\le\int_{|\xi|\le L_R}|f_R(\xi)|^2\,d\xi
\le C_d e^{2d\varepsilon R}M_R^2.
\label{E12}
\end{equation}
With \(\varepsilon\) still fixed, taking logarithms of the inequality between the first and last terms in \eqref{E12} and dividing by \(2R\) gives
$$
R^{-1}\log M_R
\ge m-(d+1)\varepsilon+\frac{1}{2R}\log\frac{C_{4,\varepsilon}}{C_d}.
$$
The last term tends to zero as \(R\to\infty\), since \(C_{4,\varepsilon}\) and \(C_d\) are independent of \(R\). Hence
$$
\liminf_{R\to\infty}R^{-1}\log M_R
\ge m-(d+1)\varepsilon.
$$
Using \eqref{E7} and \(\ell\le m\) on \(\operatorname{supp}c\), we have \(M_R\le\|u_R\|_1\le e^{Rm}\|c\|_1\), which implies \(\limsup_{R\to\infty}R^{-1}\log M_R\le m\). Letting \(\varepsilon\downarrow0\) in the preceding lower bound and combining the two bounds, we obtain
\begin{equation}
R^{-1}\log M_R\longrightarrow m
\label{E3}
\end{equation}
as \(R\to\infty\).

Choose \(\tau_R\in\mathbb R^d\) with \(|f_R(\tau_R)|=M_R>0\). Indeed, \(f_R\) is continuous and tends to zero at infinity, so it attains its maximum modulus. Define \(a_R:\mathbb R^d\to\mathbb C\) by
\begin{equation}
a_R(x)=\frac{e^{R\ell(x)-i\tau_R\cdot x}c(x)}{f_R(\tau_R)}.
\label{E14}
\end{equation}
The normalization \eqref{E14} and the choice \(|f_R(\tau_R)|=M_R\) give
\begin{equation}
g_R(\xi):=\widehat{a_R}(\xi)
=\frac{f_R(\xi+\tau_R)}{f_R(\tau_R)},\qquad
\|g_R\|_\infty=g_R(0)=1.
\label{E59}
\end{equation}

To retain \(g_R(0)=1\) in the Fourier transform of a weak limit of the \(a_R\), we seek a subsequence on which \(g_R\) converges locally uniformly. The common support bound in \eqref{E14} allows us to fix \(\chi\in C_c^\infty(\mathbb R^d)\) equal to one near \(\operatorname{supp}c\). Write $$(f*g)(\xi)=\int_{\mathbb R^d}f(\zeta)g(\xi-\zeta)\,d\zeta$$ for the ordinary convolution of \(f,g\in\mathcal S(\mathbb R^d)\). Since \(a_R=\chi a_R\) and \(\|\chi\|_1\|g_R\|_1<\infty\), Fourier inversion for \(a_R\) and Fubini's theorem give
$$
\begin{aligned}
g_R(\xi)
&=\int_{\mathbb R^d}\chi(x)a_R(x)e^{-ix\cdot\xi}\,dx\\
&=(2\pi)^{-d}\int_{\mathbb R^d}g_R(\eta)
\left[\int_{\mathbb R^d}\chi(x)e^{-ix\cdot(\xi-\eta)}\,dx\right]\,d\eta\\
&=(2\pi)^{-d}(\widehat\chi*g_R)(\xi).
\end{aligned}
$$
In the last equality, we identify the inner integral as \(\widehat\chi(\xi-\eta)\) by \eqref{E7} and set \(\zeta=\xi-\eta\) to write the result as a convolution.
In the integral with variable \(\eta\), differentiation in \(\xi\) acts on \(\widehat\chi(\xi-\eta)\). For \(\xi\) in a fixed compact set, each such derivative is bounded by a constant times \(\langle\eta\rangle^{-d-1}\), since \(\widehat\chi\) is Schwartz. We may therefore differentiate under the integral and use \(\|g_R\|_\infty=1\) from \eqref{E59} to obtain
$$
\|\partial^\alpha g_R\|_\infty=\|(2\pi)^{-d}(\partial^\alpha\widehat\chi)*g_R\|_\infty
\le(2\pi)^{-d}\|\partial^\alpha\widehat\chi\|_1
$$
for every multiindex \(\alpha\in\mathbb N_0^d\) and \(R>0\).
For \(\xi,\zeta\in\mathbb R^d\), the fundamental theorem of calculus along the segment joining them and the preceding derivative bounds give
$$
\begin{aligned}
|g_R(\xi)-g_R(\zeta)|
&\le\sum_{j=1}^d|\xi_j-\zeta_j|
\int_0^1|(\partial_{\xi_j}g_R)(\zeta+t(\xi-\zeta))|\,dt\\
&\le(2\pi)^{-d}|\xi-\zeta|
\sum_{j=1}^d\|\partial_{\xi_j}\widehat\chi\|_1.
\end{aligned}
$$
The uniform Lipschitz bound and \(\|g_R\|_\infty=1\) make the family equicontinuous and uniformly bounded on every closed ball. Apply the Arzel\`{a}--Ascoli theorem \cite[Theorem 14.1]{Treves1967} to the functions \(g_R\), \(R=1,2,\ldots\), successively on \(\{\xi\in\mathbb R^d:|\xi|\le j\}\), \(j=1,2,\ldots\), and take a diagonal subsequence. We obtain \(R_\nu\to\infty\) and a function \(g:\mathbb R^d\to\mathbb C\) such that \(g_{R_\nu}\to g\) locally uniformly as \(\nu\to\infty\), with \(\|g\|_\infty\le1\) and \(g(0)=\lim_{\nu\to\infty}g_{R_\nu}(0)=1\) by \eqref{E59}.

We use Fourier inversion and dominated convergence for both limits in \eqref{eq:localization}. Let \(q\ge1\) be an integer and let \(T_\nu,T\in\mathcal S'(\mathbb R^q)\) have Fourier transforms represented by functions \(G_\nu,G\) with \(|G_\nu|\le1\) and \(G_\nu\to G\) pointwise. For \(\varphi\in\mathcal S(\mathbb R^q)\), we use \eqref{E44} to write the Fourier pairing in \eqref{E13} as an integral. Reflecting the integration variable and applying dominated convergence with majorant \(|\widehat\varphi|\), we obtain
\begin{equation*}
\begin{aligned}
\langle T_\nu,\varphi\rangle
&=(2\pi)^{-q}\int_{\mathbb R^q}G_\nu(-\xi)\widehat\varphi(\xi)\,d\xi\\
&\longrightarrow(2\pi)^{-q}\int_{\mathbb R^q}G(-\xi)\widehat\varphi(\xi)\,d\xi
=\langle T,\varphi\rangle.
\end{aligned}
\tag{\getrefnumber{E59}a}\label{E97}
\end{equation*}
Here the last equality is \eqref{E13} for \(T\). Set \(\tau_\nu=\tau_{R_\nu}\) and \(a_\nu=a_{R_\nu}\). By \eqref{E14}, \(a_\nu\in C_c^\infty(\mathbb R^d)\), so \eqref{E18} identifies it with a tempered distribution. The bounded function \(g\) defines a tempered distribution by \eqref{E44}. Let \(a\) be its inverse Fourier transform. Since \(\widehat a_\nu=g_{R_\nu}\to g\) locally uniformly, with \(|g_{R_\nu}|\le1\) by \eqref{E59}, we may apply \eqref{E97} with \(q=d\) to conclude that \(a_\nu\to a\) weakly in \(\mathcal S'\).

To locate the support of \(a\), take \(\varphi\in C_c^\infty(\mathbb R^d)\) with \(p\notin\operatorname{supp}\varphi\). Then compactness of \(\operatorname{supp}\varphi\cap\operatorname{supp}c\) and uniqueness of the maximum give \(\ell\le m-\delta\) on this intersection for some \(\delta>0\). By \eqref{E3},
$$
M_{R_\nu}\ge e^{R_\nu(m-\delta/2)}
$$
for all sufficiently large \(\nu\). After inserting the normalization \eqref{E14} into \eqref{E44}, we use \(\ell\le m-\delta\) on \(\operatorname{supp}(c\varphi)\) in the numerator and the lower bound for \(M_{R_\nu}=|f_{R_\nu}(\tau_\nu)|\) in the denominator to obtain
$$
|\langle a_\nu,\varphi\rangle|
\le \frac{e^{R_\nu(m-\delta)}\|c\varphi\|_1}{M_{R_\nu}}
\le\|c\varphi\|_1e^{-\delta R_\nu/2}\longrightarrow0
$$
as \(\nu\to\infty\). By weak convergence, \(a\) annihilates every such test, so \(\operatorname{supp}a\subset\{p\}\) by Definition \ref{def:distribution-support}. The continuity of \(\widehat a=g\) and the normalization \(g(0)=1\) show that \(a\ne0\), so its support is exactly \(\{p\}\). As a compactly supported distribution, \(a\) has finite order by \cite[(2.3.1)]{Hormander2003}. The structure theorem for distributions supported at one point \cite[Theorem 2.3.4]{Hormander2003} therefore gives in \(\mathcal D'(\mathbb R^d)\) the finite sum
$$
a=\sum_{|\alpha|\le N}b_\alpha\partial^\alpha\delta_p
$$
for some \(N\in\mathbb N_0\) and coefficients \(b_\alpha\in\mathbb C\). The distributional derivatives satisfy \(\langle\partial^\alpha\delta_p,\varphi\rangle=(-1)^{|\alpha|}\partial^\alpha\varphi(p)\) for every \(\varphi\in\mathcal S(\mathbb R^d)\), so each derivative of \(\delta_p\) is tempered by \eqref{E51}. Since \(a\) is also tempered, density in \eqref{E18} extends the finite-sum identity to \(\mathcal S'(\mathbb R^d)\). Taking its Fourier transform with the convention \eqref{E15}, we obtain
$$
\widehat a(\xi)=e^{-ip\cdot\xi}\sum_{|\alpha|\le N}b_\alpha(i\xi)^\alpha.
$$
This Fourier expression and \(g\) are continuous and define the same distribution, so they agree pointwise by \eqref{E37} applied to their difference. Since the exponential has modulus one, the polynomial is bounded on \(\mathbb R^d\) and hence constant. Its value at the origin is \(g(0)=1\). Thus \(g(\xi)=\widehat a(\xi)=e^{-ip\cdot\xi}=\widehat{\delta_p}(\xi)\) for every \(\xi\), and Fourier inversion gives \(a=\delta_p\). This proves the first convergence in \eqref{eq:localization}.

For the tensor limit we use the same Fourier bound in both variables. Each \(a_\nu\) belongs to \(C_c^\infty(\mathbb R^d)\subset L^1(\mathbb R^d)\), and \(\|a_\nu\otimes a_\nu\|_1=\|a_\nu\|_1^2<\infty\) by \eqref{E17}. The Fourier convention \eqref{E7} and Fubini's theorem give
$$
\begin{aligned}
\widehat{\bigl[a_\nu\otimes a_\nu\bigr]}(\xi,\zeta)
&=\iint a_\nu(x)a_\nu(y)e^{-i(x\cdot\xi+y\cdot\zeta)}\,dx\,dy\\
&=\widehat a_\nu(\xi)\widehat a_\nu(\zeta)
=g_{R_\nu}(\xi)g_{R_\nu}(\zeta)
\end{aligned}
$$
for \(\xi,\zeta\in\mathbb R^d\). By \eqref{E59} this product has modulus at most one, and at \((-\xi,-\zeta)\) it tends to \(e^{ip\cdot(\xi+\zeta)}\) because \(g_{R_\nu}\to\widehat{\delta_p}\) locally uniformly. For \(\psi\in\mathcal S(\mathbb R^{2d})\), apply \eqref{E97} with \(q=2d\) and \(T_\nu=a_\nu\otimes a_\nu\). Fourier inversion evaluates the limiting pairing as
$$
\begin{aligned}
\langle a_\nu\otimes a_\nu,\psi\rangle
&\longrightarrow(2\pi)^{-2d}\iint
 e^{ip\cdot(\xi+\zeta)}\widehat\psi(\xi,\zeta)\,d\xi\,d\zeta\\
&\qquad=\psi(p,p)=\langle\delta_p\otimes\delta_p,\psi\rangle
\end{aligned}
$$
as \(\nu\to\infty\), completing the proof of \eqref{eq:localization}. The identity \(|f_{R_\nu}(\tau_\nu)|=M_{R_\nu}\) and \eqref{E3} prove \eqref{eq:fourier-growth} along the same sequence.
\end{proof}

If a linear functional separates the upper support levels of \(c\star c\) and \(2\operatorname{supp}c\), we can perturb it so that the gap persists and the new functional has a unique maximizer \(p\) on \(\operatorname{supp}c\). For the normalized functions of Lemma \ref{lem:localization} in that direction, common compact support and \(\sigma(p,p)=0\) turn the tensor limit \eqref{eq:localization} into weak convergence of their twisted squares to \(\delta_{2p}\). The gap and \eqref{eq:fourier-growth} would force the same squares to tend to zero. For subsets \(K,L\) of a Euclidean space, write \(K+L=\{x+y:x\in K,\ y\in L\}\).

\begin{lemma}\label{lem:compact-square}

For every \(c\in C_c^\infty(\mathbb R^d)\) with \(c\not\equiv0\), we have \(c\star c\not\equiv0\) and
\begin{equation}
\operatorname{conv}\operatorname{supp}(c\star c)
=2\operatorname{conv}\operatorname{supp}c. \label{eq:compact-square}
\end{equation}

\end{lemma}

\begin{proof}

\begingroup
\postdisplaypenalty=10000\relax
Write \(K=\operatorname{supp}c\), \(B_0=\max_{x\in K}|x|>0\), and define \(h_K(v)=\max_{x\in K}v\cdot x\) for \(v\in\mathbb R^d\). For \(f,g\in C_c^\infty(\mathbb R^d)\), we have \(\sigma(x,z-x)=\sigma(x,z)\) by \eqref{E55}. Each \(z\)-derivative in \eqref{eq:twisted-product} therefore only differentiates \(g(z-x)\) or contributes a polynomial factor in \(x\). For each multi-index \(\alpha\), the differentiated integrand is bounded by \(C_{\alpha,f,g}|f(x)|\), uniformly in \(z\), since the derivatives of \(g\) are bounded and \(x\in\operatorname{supp}f\) wherever the integrand is nonzero. This justifies differentiation under the integral at every order. The integral vanishes outside the compact sum of the two supports, so \eqref{E37} gives
\begin{equation*}
f\star g\in C_c^\infty(\mathbb R^d),\qquad
\operatorname{supp}(f\star g)\subset\operatorname{supp}f+\operatorname{supp}g.
\tag{\getrefnumber{eq:compact-square}a}\label{E98}
\end{equation*}
For \(x,y\in K\), the midpoint \((x+y)/2\) belongs to \(\operatorname{conv}K\), while \(2K\subset K+K\). Taking convex hulls in \eqref{E98} with \(f=g=c\) therefore gives \(\operatorname{conv}\operatorname{supp}(c\star c)\subset\operatorname{conv}(K+K)=2\operatorname{conv}K\), which proves one inclusion in \eqref{eq:compact-square}. To prove nonvanishing and the reverse inclusion, suppose, for a contradiction, that \(c\star c\equiv0\), or that \(c\star c\not\equiv0\) and
\begin{equation}
\operatorname{conv}\operatorname{supp}(c\star c)
\subsetneq2\operatorname{conv}K.
\label{E42}
\end{equation}
We turn the strict inclusion \eqref{E42} into a strict upper bound for \(v\cdot z\) on \(\operatorname{supp}(c\star c)\) in a suitable direction \(v\). For a nonempty compact set \(F\subset\mathbb R^d\) and a linear functional \(\ell\), every convex combination of points \(x_j\in F\), with \(\lambda_j\ge0\) and \(\sum_j\lambda_j=1\), satisfies
$$
\ell\!\bigg(\sum_j\lambda_jx_j\bigg)
=\sum_j\lambda_j\ell(x_j)\le\max_F\ell.
$$
Together with \(F\subset\operatorname{conv}F\), this proves
\begin{equation*}
\max_{\operatorname{conv}F}\ell=\max_F\ell.
\tag{\getrefnumber{E42}a}\label{E122}
\end{equation*}
In the case of strict inclusion, choose \(z_0\in2\operatorname{conv}K\setminus\operatorname{conv}\operatorname{supp}(c\star c)\) by \eqref{E42}. By \eqref{E98}, \(\operatorname{supp}(c\star c)\) is compact, and its convex hull is compact by \cite[Theorem 17.2]{Rockafellar1970}. This hull is nonempty in the case under consideration, so strict separation from \(\{z_0\}\) by \cite[Corollary 11.4.2]{Rockafellar1970} yields \(v\in\mathbb R^d\) such that
$$
\max_{z\in\operatorname{supp}(c\star c)}v\cdot z
<v\cdot z_0\le2h_K(v).
$$
The strict inequality comes from separation, and the upper bound follows from \(z_0/2\in\operatorname{conv}K\) and \eqref{E122} with \(\ell(x)=v\cdot x\). We can therefore choose \(\delta>0\) such that
\begin{equation}
v\cdot z\le2m-\delta
\label{eq:support-gap}
\end{equation}
for every \(z\in\operatorname{supp}(c\star c)\), where \(m=h_K(v)\).
If \(c\star c\equiv0\), its support is empty, so \eqref{eq:support-gap} holds for any \(v\in\mathbb R^d\) and \(\delta>0\).
\par\endgroup

The gap \eqref{eq:support-gap} holds in the direction \(v\), whereas Lemma \ref{lem:localization} requires a unique maximizer on \(K\). We perturb \(v\) so that the two conditions hold in the same direction. For \(t>0\), choose \(p_t\in K\) maximizing \(|x+tv|\), and put \(\mu_t=v+p_t/t\). Maximality gives
$$
0\ge |x+tv|^2-|p_t+tv|^2
=2t\mu_t\cdot(x-p_t)+|x-p_t|^2
$$
for every \(x\in K\). Then \(\mu_t\cdot(x-p_t)\le-|x-p_t|^2/(2t)\), so \(\mu_t\) uniquely exposes \(p_t\). To see how much of the gap survives this perturbation, write \(m_t=h_K(\mu_t)\) and note that
$$
m_t=\max_{x\in K}\mu_t\cdot x
\le m+\max_{x\in K}|(\mu_t-v)\cdot x|.
$$
Interchanging \(v\) and \(\mu_t\)
in this estimate bounds the difference in the opposite direction, so
$$
|m_t-m|\le\max_{x\in K}|(\mu_t-v)\cdot x|
\le B_0|\mu_t-v|\le\frac{B_0^2}{t}.
$$
For \(z\in\operatorname{supp}(c\star c)\), the inclusion \(z\in K+K\) from \eqref{E98} implies \(|z|\le2B_0\) and hence \(|(\mu_t-v)\cdot z|\le2B_0^2/t\). Combining this bound with \eqref{eq:support-gap} and \(m\le m_t+B_0^2/t\), we obtain
$$
\begin{aligned}
\mu_t\cdot z
&\le v\cdot z+\frac{2B_0^2}{t}
\le2m-\delta+\frac{2B_0^2}{t}\le2m_t-\delta+\frac{4B_0^2}{t}.
\end{aligned}
$$

Fix \(t>8B_0^2/\delta\), so that \(4B_0^2/t<\delta/2\). The perturbed direction then retains the gap
\begin{equation}
\mu_t\cdot z\le2m_t-\frac{\delta}{2}
\label{E16}
\end{equation}
for every \(z\in\operatorname{supp}(c\star c)\). The direction \(\mu_t\) now serves both parts of the argument. Its unique maximizer \(p_t\) gives the point of localization, while \eqref{E16} bounds the exponential weight on \(\operatorname{supp}(c\star c)\) strictly below its value at \(2p_t\), where \(\mu_t\cdot(2p_t)=2m_t\). Applying Lemma \ref{lem:localization} with \(\ell(x)=\mu_t\cdot x\) gives sequences \(R_\nu,\tau_\nu\) and functions \(a_\nu\in C_c^\infty(\mathbb R^d)\) supported in \(K\). To pass from the tensor limit in \eqref{eq:localization} to a weak limit for \(a_\nu\star a_\nu\), we use the common support \(K\) to localize the test function. Choose \(\chi\in C_c^\infty(\mathbb R^d)\), independent of \(\nu\), equal to one near \(K\). For \(\psi\in\mathcal S(\mathbb R^d)\), define \(\varphi_\psi\in C_c^\infty(\mathbb R^{2d})\) by
$$
\varphi_\psi(x,y)=\chi(x)\chi(y)e^{i\sigma(x,y)}\psi(x+y).
$$
For \(f,g\in C_c^\infty(\mathbb R^d)\), the smooth product in \eqref{E98} defines a tempered distribution by \eqref{E18}. After substitution of \eqref{eq:twisted-product} into \eqref{E44}, the integrand has \(L^1(dx\,dz)\) norm at most \(\|\psi\|_\infty\|f\|_1\|g\|_1\). Fubini's theorem and the change \(y=z-x\) yield
\begin{equation*}
\begin{aligned}
\langle f\star g,\psi\rangle
&=\int\left[\int f(x)g(z-x)e^{i\sigma(x,z-x)}\,dx\right]\psi(z)\,dz\\
&=\iint f(x)g(y)e^{i\sigma(x,y)}\psi(x+y)\,dx\,dy.
\end{aligned}
\tag{\getrefnumber{E16}a}\label{E99}
\end{equation*}
With \(f=g=a_\nu\), \eqref{E99} and \eqref{E17} identify \(\langle a_\nu\star a_\nu,\psi\rangle\) with \(\langle a_\nu\otimes a_\nu,\varphi_\psi\rangle\), since \(\chi=1\) near \(\operatorname{supp}a_\nu\subset K\). For the fixed \(t\), \(\varphi_\psi\in\mathcal S(\mathbb R^{2d})\) is independent of \(\nu\), so the tensor convergence in \eqref{eq:localization} gives
$$
\begin{aligned}
\langle a_\nu\otimes a_\nu,\varphi_\psi\rangle
&\rightarrow\langle\delta_{p_t}\otimes\delta_{p_t},\varphi_\psi\rangle=\varphi_\psi(p_t,p_t)=\psi(2p_t).
\end{aligned}
$$
The last equality uses \(\chi(p_t)=1\) and the identity \(\sigma(p_t,p_t)=0\) from \eqref{E55}. Since \(\psi\) was arbitrary, we conclude that
\begin{equation}
a_\nu\star a_\nu\longrightarrow\delta_{2p_t}
\label{E19}
\end{equation}
weakly in \(\mathcal S'(\mathbb R^d)\) as \(\nu\to\infty\).

To contradict this limit, we combine the support gap \eqref{E16} with \eqref{eq:fourier-growth} to show that \(\|a_\nu\star a_\nu\|_1\to0\). For \(x,z\in\mathbb R^d\), the normalization \eqref{E14} gives
$$
a_\nu(x)a_\nu(z-x)
=\frac{e^{R_\nu\mu_t\cdot z-i\tau_\nu\cdot z}}
{\widehat{\bigl[e^{R_\nu\mu_t\cdot(\cdot)}c\bigr]}(\tau_\nu)^2}
c(x)c(z-x).
$$
Substitution into \eqref{eq:twisted-product} expresses the normalized square as
\begin{equation}
(a_\nu\star a_\nu)(z)
=\frac{e^{R_\nu\mu_t\cdot z-i\tau_\nu\cdot z}}
{\widehat{\bigl[e^{R_\nu\mu_t\cdot(\cdot)}c\bigr]}(\tau_\nu)^2}
(c\star c)(z).
\label{E20}
\end{equation}

The multiplier in \eqref{E20} never vanishes, so the two squares have the same support. On that support \eqref{E16} bounds the exponential in the numerator. For the denominator we use \eqref{eq:fourier-growth}, which for all sufficiently large \(\nu\) gives \(R_\nu^{-1}\log|\widehat{[e^{R_\nu\mu_t\cdot(\cdot)}c]}(\tau_\nu)|\ge m_t-\delta/8\). Exponentiating after multiplication by \(2R_\nu\) yields
$$
\left|\widehat{\bigl[e^{R_\nu\mu_t\cdot(\cdot)}c\bigr]}(\tau_\nu)\right|^2
\ge e^{2R_\nu(m_t-\delta/8)}
=e^{R_\nu(2m_t-\delta/4)}.
$$
Taking absolute values in \eqref{E20} and integrating with these two bounds yields
$$
\begin{aligned}
\|a_\nu\star a_\nu\|_1
&\le\|c\star c\|_1
\frac{e^{R_\nu(2m_t-\delta/2)}}
{\left|\widehat{\bigl[e^{R_\nu\mu_t\cdot(\cdot)}c\bigr]}(\tau_\nu)\right|^2}\le\|c\star c\|_1e^{-\delta R_\nu/4}
\longrightarrow0.
\end{aligned}
$$
For every \(\psi\in\mathcal S(\mathbb R^d)\), \eqref{E44} and the last estimate give
$$
|\langle a_\nu\star a_\nu,\psi\rangle|
\le\|\psi\|_\infty\|a_\nu\star a_\nu\|_1\longrightarrow0.
$$
Choose \(\psi\in C_c^\infty(\mathbb R^d)\) with \(\psi(2p_t)=1\). The last estimate forces the pairing with this test to tend to zero, whereas \eqref{E19} gives the limit \(\psi(2p_t)=1\). This rules out both \(c\star c\equiv0\) and the strict inclusion \eqref{E42}, so the inclusion obtained by taking convex hulls in \eqref{E98} is the equality \eqref{eq:compact-square}.
\end{proof}

\subsection{Smooth twisted squares under condition \texorpdfstring{\((Q)\)}{(Q)}}

To apply \eqref{eq:compact-square} to a smooth function with noncompact support, we truncate near the maximum of a linear functional on its support. The next condition ensures that this part of the support is compact.

\begin{definition}[$Q$-support condition]\label{def:upper-caps}

In this paper, we say that a closed set \(K\subset\mathbb R^d\) satisfies the \emph{Q-support condition}, or condition \((Q)\), in the direction of a real-linear functional \(\ell:\mathbb R^d\to\mathbb R\) if \(\sup_{x\in K}\ell(x)<\infty\) and
\begin{equation}
K\cap\{x:\ell(x)\ge s\}\quad\hbox{is compact}
\label{eq:proper-caps}
\end{equation}
for every \(s\in\mathbb R\).

\end{definition}

Put \(\ell(x)=-e\cdot x\) for \(x\in\mathbb R^d\). For every \(S\in\mathcal C_{e,\eta}\) with support constant \(M_S\), \eqref{eq:cone} gives
\begin{equation}
\begin{aligned}
\ell(x)&\le M_S-\eta|x|\le M_S
\qquad(x\in\operatorname{supp}S),\\
\operatorname{supp}S\cap\{\ell\ge s\}
&\subset\{x\in\mathbb R^d:\eta|x|\le M_S-s\}
\qquad(s\in\mathbb R).
\end{aligned}
\label{E80}
\end{equation}
Then \(\operatorname{supp}S\) satisfies condition \((Q)\), since its closed upper level sets are bounded.

Truncating a smooth function leaves a remainder whose support may be noncompact. To control the terms involving that remainder in the twisted square, we first locate the support of a product of two smooth functions whose supports satisfy \((Q)\) in a common direction.

\begin{lemma}\label{L2}

Let \(f,g\in C^\infty(\mathbb R^d)\) have supports \(K_1,K_2\) satisfying condition \((Q)\) in the direction of a common real-linear functional \(\ell:\mathbb R^d\to\mathbb R\). Then the integral in \eqref{eq:twisted-product} converges absolutely for every \(z\in\mathbb R^d\) and defines \(f\star g\in C^\infty(\mathbb R^d)\), with
\begin{equation}
\operatorname{supp}(f\star g)\subset K_1+K_2.
\label{E22}
\end{equation}
Addition is proper on \(K_1\times K_2\), and \(f\star g\) represents the pushforward of \(e^{i\sigma(x,y)}(f\otimes g)\) under addition on compactly supported tests. For fixed \(\ell\), twisted convolution is bilinear on the class of smooth functions with supports satisfying \((Q)\) and takes values in the same class.

\end{lemma}

\begin{proof}

Choose \(m_1,m_2\in\mathbb R\) with \(\ell\le m_j\) on \(K_j\) for \(j=1,2\), and let \(L\subset\mathbb R^d\) be an arbitrary nonempty compact set. For \((x,y)\in K_1\times K_2\) with \(x+y\in L\), linearity gives \(\ell(x)+\ell(y)=\ell(x+y)\ge\min_L\ell\). Then the upper bounds \(\ell(x)\le m_1\) and \(\ell(y)\le m_2\) imply
$$
\ell(x)\ge\min_L\ell-m_2\quad\text{and}\quad
\ell(y)\ge\min_L\ell-m_1,
$$
so \((x,y)\in K_{1,L}\times K_{2,L}\), where \(K_{1,L}=K_1\cap\{\ell\ge\min_L\ell-m_2\}\) and \(K_{2,L}=K_2\cap\{\ell\ge\min_L\ell-m_1\}\) are compact by \eqref{eq:proper-caps}. The pairs in \(K_1\times K_2\) with sum in \(L\) form a closed subset of this compact product, so
\begin{equation}
\{(x,y)\in K_1\times K_2:x+y\in L\}
\quad\hbox{is compact}.
\label{E38}
\end{equation}
By Definition \ref{def:proper-addition}, \eqref{E38} is precisely the properness of addition on \(K_1\times K_2\). It also proves that \(K_1+K_2\) is closed. Indeed, any convergent sequence in the sum and its limit lie in one closed ball \(L\). The intersection \((K_1+K_2)\cap L\), being the image under addition of the compact set in \eqref{E38}, is compact and therefore contains the limit.

To define \(f\star g\) near \(z_0\in\mathbb R^d\), take \(L\) to be a compact neighborhood of \(z_0\). Choose \(\chi_1,\chi_2\in C_c^\infty(\mathbb R^d)\) with \(\chi_j=1\) on a neighborhood of \(K_{j,L}\) for \(j=1,2\). For \(z\in L\), \(f(x)g(z-x)\ne0\) requires \(x\in K_{1,L}\) and \(z-x\in K_{2,L}\), so
$$
\begin{aligned}
\int_{\mathbb R^d}|f(x)g(z-x)|\,dx
&=\int_{\mathbb R^d}|(\chi_1f)(x)(\chi_2g)(z-x)|\,dx\\
&\le\|\chi_1f\|_1\|\chi_2g\|_\infty<\infty.
\end{aligned}
$$
Thus the integral in \eqref{eq:twisted-product} converges absolutely for \(z\in L\) and satisfies
\begin{equation}
\begin{aligned}
(f\star g)(z)
&=\int_{\mathbb R^d}(\chi_1f)(x)(\chi_2g)(z-x)e^{i\sigma(x,z-x)}\,dx\\
&=((\chi_1f)\star(\chi_2g))(z).
\end{aligned}
\label{E21}
\end{equation}
Formula \eqref{E98} applies to \(\chi_1 f,\chi_2 g\), so the right-hand side of \eqref{E21} is smooth. Since \(z_0\) was arbitrary, \(f\star g\) is smooth on \(\mathbb R^d\). To obtain \((Q)\) for the product, it suffices to prove the support inclusion \eqref{E22} and verify \((Q)\) for the closed sum \(K_1+K_2\). Indeed, if a closed set \(K\) satisfies \((Q)\) and \(F\subset K\) is closed, the upper bound for \(\ell\) on \(K\) holds on \(F\), and \(F\cap\{\ell\ge s\}\) is a closed subset of the compact set \(K\cap\{\ell\ge s\}\) for each \(s\). Hence
\begin{equation*}
K\text{ satisfies }(Q),\quad F\subset K\text{ closed}
\quad\Longrightarrow\quad F\text{ satisfies }(Q).
\tag{\getrefnumber{E21}a}\label{E100}
\end{equation*}
To prove \eqref{E22}, let \(T\) denote the pushforward of \(e^{i\sigma(x,y)}(f\otimes g)\) under addition, as in Section \ref{sec:1}. Given \(\varphi\in C_c^\infty(\mathbb R^d)\), take a compact neighborhood \(L\) of \(\operatorname{supp}\varphi\) and choose \(\chi_j=1\) near \(K_{j,L}\) as in \eqref{E21}. Since \(\chi_1(x)\chi_2(y)=1\) near the pairs in \eqref{E38}, inserting these factors in the pushforward pairing leaves \(\langle T,\varphi\rangle\) unchanged. The resulting double integral is the right-hand side of \eqref{E99} for \(\chi_1 f,\chi_2 g\). By \eqref{E21}, their twisted product equals \(f\star g\) on \(L\), so
\begin{equation}
\langle T,\varphi\rangle
=\int_{\mathbb R^d}(f\star g)(z)\varphi(z)\,dz,
\qquad \varphi\in C_c^\infty(\mathbb R^d).
\label{E86}
\end{equation} Since addition is proper on \(K_1\times K_2\) by \eqref{E38}, \cite[Theorem 1, p. 44]{LoupiasMiracleSole1967} places \(\operatorname{supp}T\) in \(K_1+K_2\). By \eqref{E86}, \(f\star g\) represents \(T\) on compact tests, so \eqref{E37} gives
$$
\operatorname{supp}(f\star g)
=\overline{\{z\in\mathbb R^d:(f\star g)(z)\ne0\}}
=\operatorname{supp}T\subset K_1+K_2,
$$
which proves \eqref{E22}. To verify \((Q)\) for \(K_1+K_2\), note first that \(\ell\le m_1+m_2\) on this sum. If \(z=x+y\) with \(x\in K_1\), \(y\in K_2\), and \(\ell(z)\ge s\), then \(\ell(x)=\ell(z)-\ell(y)\ge s-m_2\) and \(\ell(y)\ge s-m_1\). Hence, for each \(s\in\mathbb R\),
$$
(K_1+K_2)\cap\{\ell\ge s\}
\subset
\bigl(K_1\cap\{\ell\ge s-m_2\}\bigr)+
\bigl(K_2\cap\{\ell\ge s-m_1\}\bigr).
$$
The set on the left is closed by the closedness of \(K_1+K_2\) proved after \eqref{E38}. The two sets on the right are compact by \((Q)\), so their sum is compact. This proves \((Q)\) for \(K_1+K_2\), and \eqref{E100} transfers it to the closed support in \eqref{E22}.
For the vector space assertion, finite unions preserve the upper bound on \(\ell\), and each upper level set of a finite union is a finite union of compact sets. The support of a linear combination is closed and lies in the union of the individual supports, so it satisfies \((Q)\) by \eqref{E100}. On this space the absolutely convergent integral \eqref{eq:twisted-product} is bilinear.
\end{proof}

For a nonzero smooth \(h\) whose support satisfies \((Q)\) in the direction \(\ell\), choose a compactly supported truncation \(c\) that agrees with \(h\) near the highest level of \(\ell\) on \(\operatorname{supp}h\). The identity \eqref{eq:compact-square} determines the highest level reached by \(\operatorname{supp}(c\star c)\), and \eqref{E22} keeps the terms involving \(h-c\) strictly below it.

\begin{proposition}\label{prop:smooth-square}

Suppose that \(h\in C^\infty(\mathbb R^d)\) and that \(\operatorname{supp}h\) satisfies condition \((Q)\) in a direction \(\ell\). If the locally defined twisted square \(h\star h\) vanishes, then \(h\equiv0\).

\end{proposition}

\begin{proof}

Assume that \(h\star h\equiv0\), and suppose for a contradiction that \(h\not\equiv0\). Let \(m=\sup_{\operatorname{supp}h}\ell\), which is finite by condition \((Q)\), and fix \(\delta>0\). By the definition of \(m\), the set \(K_\delta=\operatorname{supp}h\cap\{\ell\ge m-\delta\}\) is nonempty and satisfies \(\sup_{K_\delta}\ell=m\). Since \(K_\delta\) is compact by \eqref{eq:proper-caps}, continuity of \(\ell\) gives a point \(p\in K_\delta\) with \(\ell(p)=m\). To retain this maximum while lowering the support bound of the remainder by \(\delta\), fix \(\chi_0\in C_c^\infty(\mathbb R^d)\) with \(0\le\chi_0\le1\), equal to one on \(\{|x|\le1\}\) and zero on \(\{|x|\ge2\}\). Choose \(R_0>\max_{x\in K_\delta}|x|\), and define
$$
\chi(x)=\chi_0(x/R_0),\qquad c=\chi h\in C_c^\infty(\mathbb R^d).
$$
The cutoff satisfies \(\chi=1\) near \(K_\delta\) because \(K_\delta\subset\{|x|<R_0\}\). In particular, \(c=h\) near \(p\), so \(p\in\operatorname{supp}c\subset\operatorname{supp}h\), \(c\not\equiv0\), and \(\max_{\operatorname{supp}c}\ell=m\). The remainder \(k=h-c=(1-\chi)h\) satisfies
$$
\{x:k(x)\ne0\}\subset\operatorname{supp}h\setminus K_\delta
\subset\{x:\ell(x)<m-\delta\}.
$$
Taking closures in this inclusion and applying \eqref{E37} gives $$\operatorname{supp}k\subset\{\ell\le m-\delta\}.$$

Both \(\operatorname{supp}c\) and \(\operatorname{supp}k\) are closed subsets of \(\operatorname{supp}h\), so they inherit condition \((Q)\) by \eqref{E100}. Bilinearity in Lemma \ref{L2} now allows us to expand \(h\star h=0\) using \(h=c+k\), obtaining
$$
c\star c=-c\star k-k\star c-k\star k.
$$
Every term on the right contains \(k\), and \eqref{E22} places its support in the sum of the supports of its factors. Adding the bounds for \(\ell\) on these supports gives
$$
\begin{aligned}
\operatorname{supp}(c\star k)\cup\operatorname{supp}(k\star c)
&\subset\{\ell\le m+(m-\delta)\}=\{\ell\le2m-\delta\},\\
\operatorname{supp}(k\star k)
&\subset\{\ell\le2(m-\delta)\}\subset\{\ell\le2m-\delta\}.
\end{aligned}
$$
Thus \(h\star h=c\star c\) on the open half-space \(\{\ell>2m-\delta\}\), where all terms containing \(k\) vanish. Under the assumption \(h\star h=0\), the support inclusion for the sum \(-c\star k-k\star c-k\star k\) reads
\begin{equation}
\begin{aligned}
\operatorname{supp}(c\star c)
&\subset\operatorname{supp}(c\star k)\cup\operatorname{supp}(k\star c)\cup\operatorname{supp}(k\star k)\\
&\subset\{\ell\le2m-\delta\}.
\end{aligned}
\label{E47}
\end{equation}

Since \(c\not\equiv0\), Lemma \ref{lem:compact-square} gives \(c\star c\not\equiv0\) and determines the convex hull of its support. To compare its maximum with \eqref{E47}, apply \eqref{E122} first to \(\operatorname{supp}(c\star c)\) and then to \(\operatorname{supp}c\). Between these two applications, \eqref{eq:compact-square} identifies the convex hulls up to dilation, and \(\ell(2x)=2\ell(x)\). Hence
$$
\begin{aligned}
\max_{\operatorname{supp}(c\star c)}\ell
&=\max_{\operatorname{conv}\operatorname{supp}(c\star c)}\ell
=\max_{2\operatorname{conv}\operatorname{supp}c}\ell=2\max_{\operatorname{supp}c}\ell=2m.
\end{aligned}
$$
This equality and \eqref{E47} imply
$$
2m=\max_{z\in\operatorname{supp}(c\star c)}\ell(z)
\le2m-\delta,
$$
which contradicts \(\delta>0\). Therefore \(h\equiv0\), which completes the proof.
\end{proof}

\section{The cone-supported distribution algebra}\label{sec:3}

\subsection{Tempered products}

To reduce a zero product of distributions to the smooth-square criterion of Proposition \ref{prop:smooth-square}, we first construct an associative twisted product on \(\mathcal C_{e,\eta}\), starting with its action on compactly supported tests.

\begin{pairingdefinition}[Local twisted convolution]
Let \(U,V\in\mathcal D'(\mathbb R^d)\), and suppose that addition is proper on \(\operatorname{supp}U\times\operatorname{supp}V\). For \(\varphi\in C_c^\infty(\mathbb R^d)\), \eqref{E34} allows us to choose \(\chi_\varphi\in C_c^\infty(\mathbb R^{2d})\) equal to one near the compact set
$$
(\operatorname{supp}U\times\operatorname{supp}V)
\cap\{(a,b):a+b\in\operatorname{supp}\varphi\}.
$$
Using the tensor product in \eqref{E36}, we define \(U\star V\in\mathcal D'(\mathbb R^d)\) by
\begin{equation}
\langle U\star V,\varphi\rangle
=\big\langle U(a)\otimes V(b),
\chi_\varphi(a,b)e^{i\sigma(a,b)}\varphi(a+b)\big\rangle.
\label{E23}
\end{equation}
\end{pairingdefinition}
The tests obtained from two such cutoffs agree near \(\operatorname{supp}U\times\operatorname{supp}V\). By \eqref{E121} this is the support of \(U\otimes V\), so the pairing in \eqref{E23} is independent of \(\chi_\varphi\). For tests supported in a fixed compact \(L\subset\mathbb R^d\), properness supplies a common cutoff \(\chi_L\), so \eqref{E23} is linear. If \(\varphi_j\to0\) in \(\mathcal D(\mathbb R^d)\), choose \(L\) containing all their supports. By the product and chain rules, \(\chi_L(a,b)e^{i\sigma(a,b)}\varphi_j(a+b)\) tends to zero with all derivatives uniformly on the fixed support \(\operatorname{supp}\chi_L\). Continuity of \(U\otimes V\) now implies \(\langle U\star V,\varphi_j\rangle\to0\), as required by \cite[Theorem 2.1.4]{Hormander2003}.

For tempered factors satisfying the common cone condition \eqref{eq:cone}, the sum controls both input variables. We use this control to extend the local product from compactly supported tests to Schwartz functions.

\begin{proposition}\label{prop:tempered-product}

Let \(U,V\in\mathcal C_{e,\eta}\), with support constants \(M_U,M_V\) as in \eqref{eq:cone}. The local product \eqref{E23} has a unique tempered extension \(U\star V\). There exist smooth cutoffs \(\rho_U,\rho_V\) equal to one near the respective supports such that, for every \(\varphi\in\mathcal S(\mathbb R^d)\), the function $$j_\varphi:(a,b)\mapsto\rho_U(a)\rho_V(b)e^{i\sigma(a,b)}\varphi(a+b)$$ belongs to \(\mathcal S(\mathbb R^{2d})\). The extension acts by \(\varphi\mapsto\langle U\otimes V,j_\varphi\rangle\), as in \eqref{eq:distribution-product}.
Moreover, the support of \(U\star V\) satisfies
\begin{equation}
\operatorname{supp}(U\star V)
\subset\operatorname{supp}U+\operatorname{supp}V
\subset\{z:e\cdot z\ge\eta|z|-M_U-M_V\}. \label{eq:product-support}
\end{equation}
For a third factor \(T\) in the same cone class, associativity takes the form
$$
(U\star V)\star T=U\star(V\star T),\quad U,V,T\in\mathcal C_{e,\eta}.
$$

\end{proposition}

\begin{proof}

By \eqref{E80}, both supports satisfy condition \((Q)\) for \(\ell(a)=-e\cdot a\). The compactness argument proving \eqref{E38} in Lemma \ref{L2} thus shows that addition is proper on \(\operatorname{supp}U\times\operatorname{supp}V\), so \eqref{E23} defines the local twisted product. To extend the local product to Schwartz tests, we need the decay of \(\varphi\) in \(a+b\) to control the weights and phase derivatives in \((a,b)\). The pullback \(\varphi(a+b)\) itself is constant along each fiber of addition. Smooth cutoffs equal to one near the respective distributional supports suffice if their supports satisfy a cone bound controlling \(|a|+|b|\) in terms of \(|a+b|\). Choose \(\theta\in C^\infty(\mathbb R)\) with \(0\le\theta\le1\), \(\theta(r)=0\) for \(r\le0\), and \(\theta(r)=1\) for \(r\ge1\). For any \(S\in\mathcal C_{e,\eta}\) with support constant \(M_S\) in \eqref{eq:cone}, set \(L_S=M_S+\eta/2+2\). With \(\langle a\rangle=(1+|a|^2)^{1/2}\), define
\begin{equation*}
q_S(a)=e\cdot a-\frac\eta2\langle a\rangle+L_S,\qquad
\rho_S(a)=\theta(q_S(a)),\qquad a\in\mathbb R^d.
\tag{\getrefnumber{eq:product-support}a}\label{E101}
\end{equation*}
On \(\operatorname{supp}S\), \eqref{eq:cone} and \(\langle a\rangle\le|a|+1\) give
$$
q_S(a)\ge\eta|a|-M_S-\frac\eta2(|a|+1)+L_S
=\frac\eta2|a|+2\ge2.
$$
To bound the derivatives of \(\rho_S\), we use the estimate
$$
|\partial_a^\beta\langle a\rangle|
\le C_{d,\beta}\langle a\rangle^{1-|\beta|},
$$
which follows by repeated differentiation for \(a\in\mathbb R^d\) and \(\beta\in\mathbb N_0^d\) with \(|\beta|\ge1\). The constant \(C_{d,\beta}>0\) depends only on \(d,\beta\).
Since \(\langle a\rangle\ge1\) and \(|\beta|\ge1\), the definition of \(q_S\) and the preceding estimate imply
$$
|\partial_a^\beta q_S(a)|
\le |e|+\frac\eta2 C_{d,\beta}\langle a\rangle^{1-|\beta|}
\le |e|+\frac\eta2 C_{d,\beta}.
$$
The preceding derivative bounds for \(q_S\) and the bounds \(\|\theta^{(k)}\|_\infty<\infty\), \(k=1,2,\ldots\), control the positive-order derivatives of \(\rho_S\) by the chain rule. Together with \(0\le\rho_S\le1\), this yields for every \(\alpha\in\mathbb N_0^d\)
\begin{equation}
\sup_{a\in\mathbb R^d}|\partial_a^\alpha\rho_S(a)|
\le C_{d,\alpha,e,\eta}
\max_{0\le k\le|\alpha|}\|\theta^{(k)}\|_\infty,
\label{E45}
\end{equation}
where \(C_{d,\alpha,e,\eta}>0\) depends only on \(d,\alpha,e,\eta\) and is independent of \(L_S\). Since \(q_S\ge2\) on \(\operatorname{supp}S\) and \(\theta=0\) on \((-\infty,0]\), the cutoffs in \eqref{E101} satisfy
\begin{equation*}
\begin{aligned}
&\rho_S=1\ \text{on the open set }\{q_S>1\}\supset\operatorname{supp}S,\\
&\operatorname{supp}\rho_S\subset\{q_S\ge0\}
 \subset\{a:e\cdot a\ge(\eta/2)|a|-L_S\}.
\end{aligned}
\tag{\getrefnumber{E45}a}\label{E102}
\end{equation*}
For \(a\in\operatorname{supp}\rho_U\) and \(b\in\operatorname{supp}\rho_V\), adding the two bounds in \eqref{E102} gives
\begin{equation}
\frac\eta2(|a|+|b|)
\le e\cdot(a+b)+L_U+L_V
\le |e|\,|a+b|+L_U+L_V. \label{eq:cone-control}
\end{equation}

With these cutoffs, the function in the asserted pairing is
\begin{equation}
j_\varphi(a,b)=\rho_U(a)\rho_V(b)e^{i\sigma(a,b)}\varphi(a+b). \label{eq:schwartz-pullback}
\end{equation}
For \(\alpha\in\mathbb N_0^{2d}\), let \(\partial_{a,b}^\alpha\) denote differentiation in the real coordinates of \((a,b)\). Leibniz' rule and the bilinearity of \(\sigma\) in \eqref{E55} express
\(\partial_{a,b}^\alpha j_\varphi\) as a finite sum of terms of the form
$$
(\partial^\gamma\rho_U)(a)(\partial^\delta\rho_V)(b)
p(a,b)e^{i\sigma(a,b)}(\partial^\beta\varphi)(a+b),
$$
where \(\gamma,\delta,\beta\in\mathbb N_0^d\),
\(|\gamma|,|\delta|,|\beta|\le|\alpha|\), and \(p\) is a polynomial
of degree at most \(|\alpha|\), independent of \(\varphi\).
The inclusions
\(\operatorname{supp}(\partial^\gamma\rho_U)\subset\operatorname{supp}\rho_U\)
and
\(\operatorname{supp}(\partial^\delta\rho_V)\subset\operatorname{supp}\rho_V\)
place the support of each term in
\(\operatorname{supp}\rho_U\times\operatorname{supp}\rho_V\).
On \(\operatorname{supp}\rho_U\times\operatorname{supp}\rho_V\), \eqref{eq:cone-control} bounds the weight in \((a,b)\) in terms of the weight in \(a+b\), since
$$
\begin{aligned}
\langle(a,b)\rangle
&\le1+|a|+|b|\le1+\frac2\eta\bigl(|e|\,|a+b|+L_U+L_V\bigr)\le C_0\langle a+b\rangle,
\end{aligned}
$$
where \(C_0=1+(2/\eta)(|e|+L_U+L_V)\).
The degree bound on \(p\) and this comparison imply for each \(N\in\mathbb N_0\) that
$$
\begin{aligned}
\langle(a,b)\rangle^N|p(a,b)|
&\le C_\alpha\langle(a,b)\rangle^{N+|\alpha|}\le C_\alpha C_0^{N+|\alpha|}\langle a+b\rangle^{N+|\alpha|}.
\end{aligned}
$$
Here \(C_\alpha>0\) depends only on \(d,\alpha\) and works for all the finitely many polynomial factors. Applying \eqref{E45} to the cutoff derivatives, summing the terms, and taking the supremum yields
\begin{equation}
\begin{aligned}
&\sup_{a,b}\langle(a,b)\rangle^N
|\partial_{a,b}^\alpha j_\varphi(a,b)|\le C_{N,\alpha}
\max_{\substack{\beta\in\mathbb N_0^d\\|\beta|\le|\alpha|}}\sup_z
\langle z\rangle^{N+|\alpha|}|\partial^\beta\varphi(z)|.
\end{aligned}
\label{E4}
\end{equation}
The constant \(C_{N,\alpha}>0\) may depend on \(N,\alpha,d,e,\eta,M_U,M_V,\theta\), but is independent of \(\varphi\). In the notation of \eqref{E51}, the right-hand side of \eqref{E4} is \(C_{N,\alpha}p_{N+|\alpha|,|\alpha|}(\varphi)\). By this seminorm bound and the continuity of \(U\otimes V\), both maps in
$$
\mathcal S(\mathbb R^d)
\xrightarrow{\ \varphi\mapsto j_\varphi\ }
\mathcal S(\mathbb R^{2d})
\xrightarrow{\ \Psi\mapsto\langle U\otimes V,\Psi\rangle\ }
\mathbb C
$$
are continuous and linear. Their composition defines a tempered distribution by
\begin{equation}
\langle U\star V,\varphi\rangle
=\langle U\otimes V,j_\varphi\rangle. \label{eq:distribution-product}
\end{equation}
For \(\varphi\in C_c^\infty(\mathbb R^d)\), the product defining \(j_\varphi\) has support in
$$
(\operatorname{supp}\rho_U\times\operatorname{supp}\rho_V)
\cap\{(a,b):a+b\in\operatorname{supp}\varphi\}.
$$
This set is closed, and \eqref{eq:cone-control} bounds \(|a|+|b|\) there, so \(j_\varphi\in C_c^\infty(\mathbb R^{2d})\). By \eqref{E102} and the choice of \(\chi_\varphi\) in \eqref{E23}, both cutoffs equal one near the compact fiber \((\operatorname{supp}U\times\operatorname{supp}V)\cap\{(a,b):a+b\in\operatorname{supp}\varphi\}\). At every other point of \(\operatorname{supp}U\times\operatorname{supp}V\), the function \(\varphi(a+b)\) vanishes locally. By the support identity \eqref{E121}, the tests in \eqref{eq:distribution-product} and \eqref{E23} therefore agree near \(\operatorname{supp}(U\otimes V)\), so their pairings coincide. By the density in \eqref{E18}, the tempered extension is unique and hence independent of the cutoffs.

The properness verified at the start of the proof also makes \(\operatorname{supp}U+\operatorname{supp}V\) closed, by the compact-image argument following \eqref{E38}. The support theorem used for \eqref{E22} applies to distributions \cite[Theorem 1, p. 44]{LoupiasMiracleSole1967}, so the local product \eqref{E23} has support in \(\operatorname{supp}U+\operatorname{supp}V\). Its tempered extension has the same support because it agrees on every compactly supported smooth test. For \(a\in\operatorname{supp}U\) and \(b\in\operatorname{supp}V\), adding the bounds in \eqref{eq:cone} and using the triangle inequality yields
$$
e\cdot(a+b)\ge\eta(|a|+|b|)-M_U-M_V
\ge\eta|a+b|-M_U-M_V.
$$
Thus each \(a+b\) in the support sum satisfies the cone inequality, so $$\operatorname{supp}U+\operatorname{supp}V\subset\{z:e\cdot z\ge\eta|z|-M_U-M_V\},$$ proving \eqref{eq:product-support}.

For \(U,V,T\in\mathcal C_{e,\eta}\), the support estimate \eqref{eq:product-support} places \(U\star V\) and \(V\star T\) in the same class. Both iterated products are therefore tempered. To compare their restrictions to compactly supported tests, we verify properness of total addition on the three supports. Fix a nonempty compact set \(L\subset\mathbb R^d\). We claim that
\begin{equation}
K_L=\{(a,b,c)\in\operatorname{supp}U\times\operatorname{supp}V\times\operatorname{supp}T:a+b+c\in L\}\ \text{is compact}.
\label{E50}
\end{equation}
Indeed, \(K_L\) is closed because the supports and \(L\) are closed and addition is continuous. Let \(M_T\ge0\) be a support constant for \(T\) as in \eqref{eq:cone}. For \((a,b,c)\in K_L\), summing the three cone inequalities gives
$$
\begin{aligned}
\eta(|a|+|b|+|c|)
&\le e\cdot(a+b+c)+M_U+M_V+M_T\\
&\le |e|\max_{z\in L}|z|+M_U+M_V+M_T<\infty.
\end{aligned}
$$
This bound proves the compactness in \eqref{E50}, so total addition on the three supports satisfies the regularity at infinity required in \cite[Definition 1 and Corollary (c), pp. 43--44]{LoupiasMiracleSole1967}. Taking \(-\sigma\) as the symplectic form in that definition gives the phase in \eqref{E23}. The two bracketings have the same phase because bilinearity in \eqref{E55} implies
$$
\sigma(a,b)+\sigma(a+b,c)
=\sigma(a,b)+\sigma(a,c)+\sigma(b,c)
=\sigma(b,c)+\sigma(a,b+c).
$$
Under this properness condition the cited corollary identifies the local products, so $$(U\star V)\star T=U\star(V\star T)$$ in \(\mathcal D'\). As the two tempered distributions agree on \(C_c^\infty(\mathbb R^d)\), the density in \eqref{E18} extends the identity to all Schwartz tests. For a finite product, \eqref{eq:product-support} keeps every intermediate product in \(\mathcal C_{e,\eta}\). The three-factor identity therefore applies to consecutive blocks of factors and reduces any bracketing to the left-associated one.
\end{proof}

To localize the triple product, we shall pair compactly supported distributions with smooth functions.

\begin{pairingdefinition}[Action of a compactly supported distribution]
Let \(r\ge1\) be an integer and \(S\in\mathcal D'(\mathbb R^r)\) have compact support. Choose \(\rho\in C_c^\infty(\mathbb R^r)\) equal to one near \(\operatorname{supp}S\). Following \cite[Section 2.3, p. 44]{Hormander2003}, we define the action of \(S\) on \(f\in C^\infty(\mathbb R^r)\) by
\begin{equation}
\langle S,f\rangle=\langle S,\rho f\rangle.
\label{E56}
\end{equation}
\end{pairingdefinition}
For fixed \(\rho\), multiplication by \(\rho\) maps \(C^\infty(\mathbb R^r)\) continuously into the test functions supported in \(\operatorname{supp}\rho\), so the right side is well defined and continuous in \(f\). If \(\widetilde\rho\) is another such cutoff, then \((\rho-\widetilde\rho)f\) vanishes near \(\operatorname{supp}S\). Thus \eqref{E56} is independent of the cutoff and agrees with the original action on \(C_c^\infty(\mathbb R^r)\). If \(S\) is tempered, its original action and \(f\mapsto\langle S,\rho f\rangle\) are continuous on \(\mathcal S(\mathbb R^r)\) and agree on that dense subspace. By \eqref{E18} the two actions agree for every Schwartz test as well.

\begin{corollary}\label{cor:local-pairing}
Let \(U,V,T\in\mathcal C_{e,\eta}\), let \(L\subset\mathbb R^d\) be compact, and let \(K_L\) be the compact inverse image in \eqref{E50}. For \((a,b,c)\in(\mathbb R^d)^3\), define the coordinate projections by \(\pi_1(a,b,c)=a\), \(\pi_2(a,b,c)=b\), and \(\pi_3(a,b,c)=c\). Choose \(\chi_j\in C_c^\infty(\mathbb R^d)\) with \(\chi_j=1\) near \(\pi_j(K_L)\) for \(j=1,2,3\), and set
$$
P=\chi_1U,\qquad Q=\chi_2V,\qquad R=\chi_3T.
$$
These cutoffs are fixed for the chosen \(U,V,T,L\). For every \(\varphi\in C_c^\infty(\mathbb R^d)\) with \(\operatorname{supp}\varphi\subset L\), we have
\begin{equation}
\begin{aligned}
&\quad\langle U\star V\star T,\varphi\rangle=\langle(P\star Q)\star R,\varphi\rangle\\
&=\big\langle P(a)\otimes Q(b)\otimes R(c),
e^{i\sigma(a,b)+i\sigma(a+b,c)}\varphi(a+b+c)\big\rangle.
\end{aligned}
\label{E39}
\end{equation}
The last pairing is interpreted by \eqref{E56}.
\end{corollary}

\begin{proof}
We compare the original and truncated triple products on tests supported in \(L\) by expressing both pairings in terms of the fixed tensor \(U\otimes V\otimes T\). The compactness of \(K_L\) in \eqref{E50} makes each projection \(\pi_j(K_L)\) compact, so the same cutoffs \(\chi_j\) serve for all these tests. For such a \(\varphi\), define \(\Psi_\varphi\in C^\infty(\mathbb R^{3d})\) by
$$
\Psi_\varphi(a,b,c)=e^{i\sigma(a,b)+i\sigma(a+b,c)}\varphi(a+b+c).
$$
Iterating \eqref{E81} gives
$$
P\otimes Q\otimes R
=\bigl(\chi_1(a)\chi_2(b)\chi_3(c)\bigr)(U\otimes V\otimes T)
$$
in \(\mathcal S'(\mathbb R^{3d})\). To use this identity with \(\Psi_\varphi\), choose \(\omega\in C_c^\infty(\mathbb R^{3d})\) equal to one near \(\operatorname{supp}\chi_1\times\operatorname{supp}\chi_2\times\operatorname{supp}\chi_3\) and test it on \(\omega\Psi_\varphi\). On the left, \eqref{E121} and \eqref{E56} allow us to omit \(\omega\). On the right, \(\omega(\chi_1\otimes\chi_2\otimes\chi_3)=\chi_1\otimes\chi_2\otimes\chi_3\). We obtain
$$
\langle P\otimes Q\otimes R,\Psi_\varphi\rangle
=\big\langle U\otimes V\otimes T,\chi_1(a)\chi_2(b)\chi_3(c)\Psi_\varphi\big\rangle.
$$
For the truncated factors, \eqref{eq:product-support} places \(\operatorname{supp}(P\star Q)\) in the compact set \(\operatorname{supp}P+\operatorname{supp}Q\). In both applications of \eqref{E23} we may therefore take the cutoff equal to one near the entire product of the relevant supports and use \eqref{E56} to omit it from the pairing. Pairing first with \(R\) in the outer convolution leaves the test
$$
F(s)=\big\langle R(c),e^{i\sigma(s,c)}\varphi(s+c)\big\rangle
$$
for \(P\star Q\). With a fixed cutoff in \(c\) as in \eqref{E56}, we may differentiate this pairing in \(s\) by \cite[Theorem 2.1.3]{Hormander2003}. If \(s\notin\{z-c:z\in\operatorname{supp}\varphi,\ c\in\operatorname{supp}R\}\), then \(s+\operatorname{supp}R\) is disjoint from \(\operatorname{supp}\varphi\). Thus \(c\mapsto\varphi(s+c)\) vanishes near \(\operatorname{supp}R\), so \(F(s)=0\). This difference set is compact as the image of two compact supports under subtraction, and \eqref{E37} places \(\operatorname{supp}F\) in it. We may now test \(P\star Q\) against \(F\in C_c^\infty(\mathbb R^d)\) in \eqref{E23}, which pairs \(P(a)\otimes Q(b)\) with \(e^{i\sigma(a,b)}F(a+b)\). Substitution of \(F(a+b)\) and the iterated pairing \eqref{E36} yield
$$
\langle(P\star Q)\star R,\varphi\rangle
=\langle P\otimes Q\otimes R,\Psi_\varphi\rangle.
$$
For the original factors, choose \(\zeta\in C_c^\infty(\mathbb R^{2d})\) equal to one near \(\{(a+b,c):(a,b,c)\in K_L\}\). If \(s\in\operatorname{supp}(U\star V)\), \(c\in\operatorname{supp}T\), and \(s+c\in\operatorname{supp}\varphi\), then by \eqref{eq:product-support} we may write \(s=a+b\) with \(a\in\operatorname{supp}U\) and \(b\in\operatorname{supp}V\). Hence \((a,b,c)\in K_L\), so \(\zeta=1\) near \((s,c)\) as required by \eqref{E23}. To choose one inner cutoff for all \(c\), project \(\operatorname{supp}\zeta\) onto the sum variable. The resulting compact set \(L_\zeta=\{z\in\mathbb R^d:(z,c)\in\operatorname{supp}\zeta\text{ for some }c\in\mathbb R^d\}\) contains every \(a+b\) for which \(\zeta(a+b,c)\ne0\). By \eqref{E38}, we may therefore choose \(\rho\in C_c^\infty(\mathbb R^{2d})\) with
$$
\rho=1\quad\text{near }\{(a,b)\in\operatorname{supp}U\times\operatorname{supp}V:a+b\in L_\zeta\}.
$$
With \(\zeta\) as the outer cutoff in \eqref{E23}, pairing first with \(T\) leaves the function
$$
H(s)=\big\langle T(c),\zeta(s,c)e^{i\sigma(s,c)}\varphi(s+c)\big\rangle
$$
as the test for \(U\star V\). The \(c\)-supports lie in the compact projection of \(\operatorname{supp}\zeta\), so \(H\) is smooth by \cite[Theorem 2.1.3]{Hormander2003}. Since \(\operatorname{supp}H\subset L_\zeta\), the chosen \(\rho\) is admissible for the inner convolution. Substituting \(H(a+b)\) in \eqref{E23} and applying \eqref{E36} to the compactly supported test \(\rho(a,b)\zeta(a+b,c)\Psi_\varphi(a,b,c)\), we obtain
\begin{equation*}
\begin{aligned}
\langle U\star V\star T,\varphi\rangle
&=\big\langle U\otimes V\otimes T,\rho(a,b)\zeta(a+b,c)\Psi_\varphi\big\rangle\\
&=\big\langle U\otimes V\otimes T,\chi_1(a)\chi_2(b)\chi_3(c)\Psi_\varphi\big\rangle.
\end{aligned}
\tag{\getrefnumber{E39}a}\label{E103}
\end{equation*}
To justify the replacement of cutoffs in \eqref{E103}, let \((a,b,c)\in K_L\). The outer cutoff \(\zeta\) equals one near \((a+b,c)\), so \(a+b\in L_\zeta\) and the inner cutoff \(\rho\) equals one near \((a,b)\). Their product is therefore one near \(K_L\), as is \(\chi_1(a)\chi_2(b)\chi_3(c)\) by the choice of the coordinate cutoffs. At every other point of the product of the three supports, \(a+b+c\notin L\), so \(\Psi_\varphi\) vanishes on a neighborhood of that point. By \eqref{E121} applied twice, the two tests therefore agree near \(\operatorname{supp}(U\otimes V\otimes T)\). The cutoff identity obtained by testing with \(\omega\Psi_\varphi\) identifies the last pairing in \eqref{E103} with \(\langle P\otimes Q\otimes R,\Psi_\varphi\rangle\). The calculation with \(F(a+b)\) identified this same pairing with \(\langle(P\star Q)\star R,\varphi\rangle\), proving \eqref{E39}.
\end{proof}

The cutoffs in Corollary \ref{cor:local-pairing} work for all tests supported in a given compact set. For the triple product \(V\star g\star U\) in Subsection \ref{sec:weyl-smoothing} and the two-factor composition in Proposition \ref{prop:bounded-composition}, we need a single sequence of truncations that converges against every Schwartz test.

\begin{lemma}\label{L5}

Let \(U,V\in\mathcal C_{e,\eta}\). If \(\chi\in C_c^\infty(\mathbb R^d)\) equals one near zero, set
$$
U_N=\chi(\cdot/N)U\quad\text{and}\quad V_N=\chi(\cdot/N)V
$$
for \(N=1,2,\ldots\), with multiplication by the cutoff defined by \eqref{E35}. Then
\postdisplaypenalty=10000\relax
\begin{equation}
U_N\star V_N\longrightarrow U\star V
\label{eq:cutoff-convergence}
\end{equation}
weakly in \(\mathcal S'\) as \(N\to\infty\).
More generally, fix an integer \(m\ge1\) and \(T_1,\ldots,T_m\in\mathcal C_{e,\eta}\). If \(T_{j,N}=\chi(\cdot/N)T_j\) for \(j=1,\ldots,m\), then
\begin{equation}
T_{1,N}\star\cdots\star T_{m,N}
\longrightarrow T_1\star\cdots\star T_m
\label{eq:finite-cutoff}
\end{equation}
weakly in \(\mathcal S'\) as \(N\to\infty\).

\end{lemma}

\begin{proof}

Choose support constants \(M_U,M_V\) in \eqref{eq:cone} and take the cutoffs \(\rho_S\) in \eqref{E101} for \(S\in\{U,V\}\), independently of \(N\). For \(\varphi\in\mathcal S(\mathbb R^d)\), define \(j_\varphi\) by \eqref{eq:schwartz-pullback}, so \(j_\varphi\in\mathcal S(\mathbb R^{2d})\) by \eqref{E4}. Since \(\operatorname{supp}U_N\subset\operatorname{supp}U\) and \(\operatorname{supp}V_N\subset\operatorname{supp}V\), the same cutoffs apply to \(U_N,V_N\). Formula \eqref{eq:distribution-product} pairs \(U_N\otimes V_N\) with \(j_\varphi\). Transferring the two factors \(\chi(a/N)\) and \(\chi(b/N)\) to the test by \eqref{E81}, we obtain
\begin{equation}
\begin{aligned}
\langle U_N\star V_N,\varphi\rangle
&=\langle U_N\otimes V_N,j_\varphi\rangle=\langle U\otimes V,\chi(a/N)\chi(b/N)j_\varphi(a,b)\rangle.
\end{aligned}
\label{E52}
\end{equation}
To pass to the limit through \(U\otimes V\) in \eqref{E52}, we estimate the cutoff on Schwartz tests in arbitrary dimension, so that the same bound also applies to a finite product. Let \(\lambda\ge1\) be an integer, let \(\psi\in\mathcal S(\mathbb R^\lambda)\), and choose \(\vartheta\in C_c^\infty(\mathbb R^\lambda)\) with \(\vartheta=1\) on \(B(0,r)\) for some \(r>0\). Write \(d_N(w)=\vartheta(w/N)-1\) for \(w\in\mathbb R^\lambda\). All derivatives of \(d_N\), including \(d_N\) itself, vanish on \(B(0,rN)\). Fix \(k\in\mathbb N_0\). For \(\gamma\in\mathbb N_0^\lambda\) with \(0<|\gamma|\le k\),
$$
\partial^\gamma d_N(w)=N^{-|\gamma|}(\partial^\gamma\vartheta)(w/N).
$$
Since \(N\ge1\) and \(|d_N|\le1+\|\vartheta\|_\infty\), all derivatives through order \(k\) satisfy \(\|\partial^\gamma d_N\|_\infty\le C_{k,\vartheta}\), with \(C_{k,\vartheta}>0\) independent of \(N\). For \(|\alpha|\le k\), we apply this bound to the Leibniz expansion
$$
\begin{aligned}
|\partial^\alpha(d_N\psi)(w)|
&\le\sum_{\beta\le\alpha}\binom{\alpha}{\beta}
 |\partial^{\alpha-\beta}d_N(w)|\,|\partial^\beta\psi(w)|\\
&\le 2^k C_{k,\vartheta}\max_{|\beta|\le k}|\partial^\beta\psi(w)|.
\end{aligned}
$$
For \(M\in\mathbb N_0\), the weight satisfies \(\langle w\rangle^M\le(rN)^{-1}\langle w\rangle^{M+1}\) on \(|w|\ge rN\), and \(\partial^\alpha(d_N\psi)(w)=0\) on \(|w|<rN\). Taking the weighted supremum in \(w\) and the maximum over \(|\alpha|\le k\) as in \eqref{E51}, we obtain
\begin{equation}
p_{M,k}\bigl((\vartheta(\cdot/N)-1)\psi\bigr)
\le \frac{2^k C_{k,\vartheta}}{rN}\,p_{M+1,k}(\psi).
\label{E5}
\end{equation}
By \eqref{E17}, the cutoff in \eqref{E52} is \((\chi\otimes\chi)(a/N,b/N)\). Taking \(\lambda=2d\), \(\vartheta=\chi\otimes\chi\), and \(\psi=j_\varphi\) in \eqref{E5} shows that the test on the right of \eqref{E52} converges to \(j_\varphi\) in \(\mathcal S(\mathbb R^{2d})\) as \(N\to\infty\). Continuity of \(U\otimes V\) and the identity \(\langle U\otimes V,j_\varphi\rangle=\langle U\star V,\varphi\rangle\) from \eqref{eq:distribution-product} now prove \eqref{eq:cutoff-convergence}.

For a single factor, \eqref{E35} gives \(\langle T_{1,N}-T_1,\varphi\rangle=\langle T_1,(\chi(\cdot/N)-1)\varphi\rangle\), which tends to zero by \eqref{E5} with \(\lambda=d\) and \(\vartheta=\chi\). When \(m>2\), the product of the first \(m-1\) truncated factors need not equal the cutoff of their product, as a direct iteration of \eqref{eq:cutoff-convergence} would require. Instead we keep \(T_1\otimes\cdots\otimes T_m\) fixed, as in \eqref{E52}, and construct a single Schwartz test to which \eqref{E5} applies in all \(m\) variables. For each \(j\), take \(S=T_j\) in \eqref{E101} with its support constant \(M_{T_j}\) from \eqref{eq:cone}. The corresponding replacement for \eqref{eq:schwartz-pullback} is the function \(j_\varphi:(\mathbb R^d)^m\to\mathbb C\), defined for \(\varphi\in\mathcal S(\mathbb R^d)\) by
\begin{equation*}
j_\varphi(a_1,\ldots,a_m)=
\prod_{j=1}^m\rho_{T_j}(a_j)
\exp\!\bigg(i\sum_{j<k}\sigma(a_j,a_k)\bigg)
\varphi\!\bigg(\sum_{j=1}^m a_j\bigg).
\tag{\getrefnumber{E5}a}\label{E104}
\end{equation*}
For \(a_j\in\operatorname{supp}\rho_{T_j}\), \(j=1,\ldots,m\), summing \eqref{E102} over the \(m\) factors gives
\begin{equation*}
\begin{aligned}
\frac\eta2\sum_{j=1}^m|a_j|
&\le e\cdot\sum_{j=1}^m a_j+\sum_{j=1}^m L_{T_j}\le |e|\bigg|\sum_{j=1}^m a_j\bigg|+\sum_{j=1}^m L_{T_j}.
\end{aligned}
\tag{\getrefnumber{E5}b}\label{E105}
\end{equation*}
Each derivative of order \(k\) of the exponential in \eqref{E104} is that exponential multiplied by a polynomial of degree at most \(k\). Estimate \eqref{E105} bounds the resulting polynomial factors and the Schwartz weight in terms of \(\langle\sum_j a_j\rangle\), while \eqref{E45} bounds the cutoff derivatives. With \(N,k\in\mathbb N_0\), the Leibniz estimate in \eqref{E4} thus takes the form
$$
p_{N,k}^{(a_1,\ldots,a_m)}(j_\varphi)
\le C_{N,k}\,p_{N+k,k}(\varphi).
$$
For the fixed factors and cutoffs the constant is independent of \(\varphi\). Thus \(\varphi\mapsto j_\varphi\) is continuous from \(\mathcal S(\mathbb R^d)\) to \(\mathcal S(\mathbb R^{md})\). Its composition with the tempered tensor \(T_1\otimes\cdots\otimes T_m\) defines a tempered distribution in the sum variable.

We identify this distribution with \(T_1\star\cdots\star T_m\) by induction on \(m\), initially testing against \(\varphi\in C_c^\infty(\mathbb R^d)\). For one factor the identity follows from \(\rho_{T_1}=1\) near \(\operatorname{supp}T_1\). When \(m\ge2\), write \(s=\sum_{j<m}a_j\). The tuples with \(a_j\in\operatorname{supp}T_j\) and \(s+a_m\in\operatorname{supp}\varphi\) form a closed set, bounded by \eqref{E105} since \(\operatorname{supp}T_j\subset\operatorname{supp}\rho_{T_j}\). Choose \(\zeta\in C_c^\infty(\mathbb R^{2d})\) equal to one near the image of this set under \((a_1,\ldots,a_m)\mapsto(s,a_m)\). Repeated use of \eqref{eq:product-support} places \(\operatorname{supp}(T_1\star\cdots\star T_{m-1})\) in \(\sum_{j<m}\operatorname{supp}T_j\). Thus every \(s\in\operatorname{supp}(T_1\star\cdots\star T_{m-1})\) can be written as \(\sum_{j<m}a_j\) with \(a_j\in\operatorname{supp}T_j\). If also \(a_m\in\operatorname{supp}T_m\) and \(s+a_m\in\operatorname{supp}\varphi\), this tuple belongs to the compact set above, so \(\zeta=1\) near the outer support fiber in \eqref{E23}. To apply the induction hypothesis to the first \(m-1\) factors, pair first in \(a_m\) and set
$$
\omega(s)=\big\langle T_m(a_m),\zeta(s,a_m)e^{i\sigma(s,a_m)}\varphi(s+a_m)\big\rangle.
$$
As for the inner test in \eqref{E103}, the two projections of \(\operatorname{supp}\zeta\) control different variables. The \(s\)-projection contains \(\operatorname{supp}\omega\), and the \(a_m\)-projection is a common compact support for the tests paired with \(T_m\). Hence \(\omega\in C_c^\infty(\mathbb R^d)\) by \cite[Theorem 2.1.3]{Hormander2003}. With the outer cutoff \(\zeta\), \eqref{E23} and \eqref{E36} give
$$
\langle(T_1\star\cdots\star T_{m-1})\star T_m,\varphi\rangle
=\langle T_1\star\cdots\star T_{m-1},\omega\rangle.
$$
Applying the induction hypothesis to the right-hand side introduces the cutoffs \(\prod_{j<m}\rho_{T_j}(a_j)\). On the support of \(\zeta(s,a_m)\), both \(s\) and \(a_m\) are bounded. Estimate \eqref{E105} for the first \(m-1\) cutoffs then bounds \(\sum_{j<m}|a_j|\). Thus \eqref{E36} expresses the iterated product on \(\varphi\) as the action of \(T_1\otimes\cdots\otimes T_m\) on the compactly supported smooth test
$$
\bigg(\prod_{j<m}\rho_{T_j}(a_j)\bigg)\zeta(s,a_m)
\exp\!\bigg(i\sum_{j<k<m}\sigma(a_j,a_k)+i\sigma(s,a_m)\bigg)
\varphi(s+a_m).
$$
On \(\prod_{j=1}^m\operatorname{supp}T_j\), the factors \(\zeta(s,a_m)\) and \(\rho_{T_m}(a_m)\) both equal one near the tuples with \(s+a_m\in\operatorname{supp}\varphi\), while \(\varphi(s+a_m)\) vanishes in a neighborhood of every other tuple. After replacing \(\zeta(s,a_m)\) by \(\rho_{T_m}(a_m)\), the test is still compactly supported, since \(\sum_j a_j\in\operatorname{supp}\varphi\) on its support and \eqref{E105} bounds all \(a_j\) there. The two tests agree near the tensor support, so the replacement preserves the pairing as in \eqref{E103}. By \eqref{E55}, the outer phase contributes \(\sigma(s,a_m)=\sum_{j<m}\sigma(a_j,a_m)\), supplying exactly the pairs with second index \(m\) in \eqref{E104}. This proves
\begin{equation}
\langle T_1\star\cdots\star T_m,\varphi\rangle
=\langle T_1\otimes\cdots\otimes T_m,j_\varphi\rangle
\label{E92}
\end{equation}
for compactly supported smooth tests. The iterated product is tempered by Proposition \ref{prop:tempered-product}, and the pairing on the right is continuous by the seminorm estimate following \eqref{E105}. Hence \eqref{E92} extends by the density in \eqref{E18} to every \(\varphi\in\mathcal S(\mathbb R^d)\).

Since \(\operatorname{supp}T_{j,N}\subset\operatorname{supp}T_j\), formula \eqref{E92} applies to the original and truncated factors with the same \(\rho_{T_j}\). Transfer the factors \(\chi(a_j/N)\) to \(j_\varphi\) by \eqref{E81}. With \(\lambda=md\) and \(\vartheta=\chi^{\otimes m}\), \eqref{E5} gives convergence of the test difference in \(\mathcal S(\mathbb R^{md})\), so continuity of \(T_1\otimes\cdots\otimes T_m\) yields
$$
\begin{aligned}
&\quad\langle T_{1,N}\star\cdots\star T_{m,N}-T_1\star\cdots\star T_m,\varphi\rangle\\
&=\bigg\langle T_1\otimes\cdots\otimes T_m,
\bigg(\prod_{j=1}^m\chi(a_j/N)-1\bigg)j_\varphi\bigg\rangle
\longrightarrow0.
\end{aligned}
$$
Since \(\varphi\in\mathcal S(\mathbb R^d)\) is arbitrary, \eqref{eq:finite-cutoff} holds in the weak topology of \(\mathcal S'(\mathbb R^d)\), which completes the proof.
\end{proof}

\subsection{Weyl maps and compactly supported products}\label{sec:4}

Grossmann, Loupias and Stein relate operator composition to twisted convolution by taking the symplectic Fourier transform of the symbols \cite[Section 2, (7)--(10)]{GrossmannLoupiasStein1968}. For extensions to distributional symbols and multiplier spaces, see \cite{Maillard1986,DuboisVioletteEtAl2001,Soloviev2012}. We fix the normalization of the Weyl operators to match the phase in \eqref{eq:twisted-product} (see \cite[Sections 1.3--1.4 and Chapter 2]{Folland1989}).

The Weyl operator \(W_a:L^2(\mathbb R^n)\to L^2(\mathbb R^n)\) associated with \(a=(s,t)\in\mathbb R^{2n}\), with \(s,t\in\mathbb R^n\), is defined by
\begin{equation}
(W_a u)(x)=e^{is\cdot t}e^{-i\sqrt2\,t\cdot x}u(x-\sqrt2\,s)
\label{E24}
\end{equation}
for \(x\in\mathbb R^n\). This is the operator \(\rho(p,q,0)\) in \cite[(1.25)]{Folland1989} with \(p=-\sqrt2\,s\) and \(q=-t/(\sqrt2\pi)\).
We write \(A^*\) for the Hilbert-space adjoint of a bounded operator \(A\) between Hilbert spaces. For \(a=(s,t)\) and \(b=(s',t')\), applying \eqref{E24} twice produces the translation in \(W_{a+b}\). Relative to the exponent in \(W_{a+b}\), the additional phase is \(i\) times
$$
s\cdot t+s'\cdot t'+2s\cdot t'-(s+s')\cdot(t+t')
=s\cdot t'-t\cdot s'=\sigma(a,b)
$$
by \eqref{E55}. Each \(W_a\) is unitary, being a translation followed by multiplication by a function of modulus one. Taking \(b=-a\) in the product formula therefore identifies its adjoint, so
\begin{equation}
W_aW_b=e^{i\sigma(a,b)}W_{a+b},\qquad W_a^*=W_{-a}.
\label{eq:weyl-relations}
\end{equation}
The identity \(\langle w,\overline v\rangle=\langle w,v\rangle_2\) for \(w\in L^2(\mathbb R^n)\) and \(v\in\mathcal S(\mathbb R^n)\) follows from \eqref{E44}, since \(w\overline v\in L^1(\mathbb R^n)\) by Cauchy--Schwarz.

For \(u,v\in\mathcal S(\mathbb R^n)\), formula \eqref{E24} expresses the Weyl coefficient as
\begin{equation}
\langle W_{(s,t)}u,v\rangle_2
=e^{is\cdot t}\int_{\mathbb R^n}u(x-\sqrt2\,s)\overline{v(x)}
e^{-i\sqrt2\,t\cdot x}\,dx.
\label{E25}
\end{equation}
To determine the Schwartz regularity of the Weyl coefficient in \((s,t)\), we use two continuous automorphisms of \(\mathcal S(\mathbb R^{p+q})\). For integers \(p,q\ge1\), an invertible linear map \(J:\mathbb R^{p+q}\to\mathbb R^{p+q}\), and \(H\in\mathcal S(\mathbb R^{p+q})\), these are given by
\begin{equation*}
\begin{aligned}
(J^*H)(x,y)&=H(J(x,y)),\quad(\mathcal F_y H)(x,\xi)&=\int_{\mathbb R^q}H(x,y)e^{-iy\cdot\xi}\,dy,
\end{aligned}
\tag{\getrefnumber{E25}a}\label{E107}
\end{equation*}
for $ x\in\mathbb R^p, y,\xi\in\mathbb R^q.$
Invertible linear changes of variables preserve the Schwartz topology \cite[p. 268]{Treves1967}. The partial Fourier transform \(\mathcal F_y\) in the convention \eqref{E7} acts by the Fourier isomorphism on the \(y\)-factor of the Schwartz tensor product, at frequency \(\xi/(2\pi)\) in the convention of \cite[Theorems 25.1 and 51.6]{Treves1967}. In \eqref{E25}, apply \eqref{E107} to \(H=u\otimes\overline v\) with \(J(s,x)=(x-\sqrt2s,x)\), then transform in \(x\) at \(\sqrt2t\). The tensor-product estimate \eqref{E118} shows that the resulting function is Schwartz in \((s,t)\) and depends continuously on \(u,\overline v\). Since every derivative of \(e^{is\cdot t}\) is a polynomial times this phase, the multiplier criterion \eqref{E89} applied to \eqref{E25} gives the continuous bilinear map
\begin{equation}
\mathcal S(\mathbb R^n)\times\mathcal S(\mathbb R^n)\longrightarrow\mathcal S(\mathbb R^{2n})
\colon\quad (u,\overline v)\longmapsto\bigl[a\mapsto\langle W_a u,v\rangle_2\bigr].
\label{E88}
\end{equation} For fixed \(U\in\mathcal S'(\mathbb R^{2n})\) and \(u\in\mathcal S(\mathbb R^n)\), composing the continuous coefficient map in \eqref{E88} with \(U\) shows that the map
$$
\varphi\longmapsto
\big\langle U(a),\langle W_a u,\overline\varphi\rangle_2\big\rangle,
\qquad \varphi\in\mathcal S(\mathbb R^n),
$$
is a continuous linear functional on \(\mathcal S(\mathbb R^n)\), hence a tempered distribution. In this way we extend the integrated Schr\"odinger representation \cite[(1.28) and Theorem (1.30)]{Folland1989} to tempered distributions in the normalization \eqref{E24}, using our complex-linear pairing.

\begin{definition}[Weyl map and spreading distribution]\label{def:weyl-map}
The Weyl map associated with \(U\in\mathcal S'(\mathbb R^{2n})\) is the operator \(W(U):\mathcal S(\mathbb R^n)\to\mathcal S'(\mathbb R^n)\) defined, for \(u,v\in\mathcal S(\mathbb R^n)\), by
\begin{equation}
\langle W(U)u,\overline v\rangle
=\big\langle U(a),\langle W_a u,v\rangle_2\big\rangle.
\label{eq:weyl-map}
\end{equation}
Here \(a\in\mathbb R^{2n}\). The distribution \(U\) is called the spreading distribution of \(W(U)\).
\end{definition}

By \eqref{E88} and \eqref{eq:weyl-map}, the form \(b_U(u,\varphi)=\langle W(U)u,\varphi\rangle\) is jointly continuous and bilinear on \(\mathcal S(\mathbb R^n)\times\mathcal S(\mathbb R^n)\). For any such form \(b\), continuity at \((0,0)\) and the increasing seminorms in \eqref{E51} give \(C>0\) and \(N,k\in\mathbb N_0\), depending only on \(b\), such that for every bounded \(D\subset\mathcal S(\mathbb R^n)\) we have
\begin{equation*}
\sup_{\varphi\in D}|b(u,\varphi)|
\le C p_{N,k}(u)\sup_{\varphi\in D}p_{N,k}(\varphi)
\qquad (u\in\mathcal S(\mathbb R^n)).
\tag{\getrefnumber{eq:weyl-map}a}\label{E129}
\end{equation*}
Since \(p_{N,k}\) is bounded on \(D\), \eqref{E129} with \(b=b_U\) shows that \(W(U)u\to0\) uniformly on each bounded set of tests as \(u\to0\) in \(\mathcal S(\mathbb R^n)\). This proves continuity of \(W(U)\) into the strong dual. The following kernel formula recovers \(U\) from \(W(U)\).

\begin{lemma}\label{L3}

For \(U\in\mathcal S'(\mathbb R^{2n})\), the distribution kernel \(K_U\in\mathcal S'(\mathbb R^n\times\mathbb R^n)\) of \(W(U)\) is given by
\begin{equation}
K_U(x,y)=2^{-n/2}\int_{\mathbb R^n}
U\!\left(\frac{x-y}{\sqrt2},t\right)
e^{-it\cdot(x+y)/\sqrt2}\,dt.
\label{eq:weyl-kernel}
\end{equation}
The integral is interpreted in the sense of distributions. The map \(U\mapsto K_U\) is an isomorphism from \(\mathcal S'(\mathbb R^{2n})\) onto \(\mathcal S'(\mathbb R^n\times\mathbb R^n)\), restricting to an isomorphism between the corresponding Schwartz spaces. In particular, \(U\mapsto W(U)\) is injective.

\end{lemma}

\begin{proof}

The Schwartz kernel and Weyl correspondences extend to tempered distributions \cite[Propositions 3.7, 3.18, 5.1 and 5.4]{KeylKiukasWerner2016}. In our pairing convention the kernel is characterized by \(\langle K_U,\overline v\otimes u\rangle=\langle W(U)u,\overline v\rangle\) for \(u,v\in\mathcal S(\mathbb R^n)\). To compute it, first take \(U\in\mathcal S(\mathbb R^{2n})\), identified with its regular distribution by \eqref{E18}. We write the pairing with \(U\) in \eqref{eq:weyl-map} as an integral by \eqref{E44} and insert the coefficient formula \eqref{E25}. Absolute convergence follows from Cauchy--Schwarz in \(x\), which gives
$$
\int_{\mathbb R^{2n}}|U(s,t)|
\int_{\mathbb R^n}|u(x-\sqrt2s)|\,|v(x)|\,dx\,ds\,dt
\le\|U\|_1\|u\|_2\|v\|_2<\infty.
$$
We may therefore apply Fubini's theorem and substitute \(y=x-\sqrt2s\). With \(ds=2^{-n/2}\,dy\) and \(is\cdot t-i\sqrt2t\cdot x=-it\cdot(x+y)/\sqrt2\), the pairing becomes
\begin{equation*}
\begin{aligned}
\langle W(U)u,\overline v\rangle
&=2^{-n/2}\iiint_{\mathbb R^{3n}}
U\!\left(\frac{x-y}{\sqrt2},t\right)
e^{-it\cdot(x+y)/\sqrt2}\\
&\quad \cdot u(y)\overline{v(x)}\,dt\,dy\,dx\\
&=\iint_{\mathbb R^{2n}}K_U(x,y)u(y)\overline{v(x)}\,dy\,dx,
\end{aligned}
\tag{\getrefnumber{eq:weyl-kernel}a}\label{E115}
\end{equation*}
where \(K_U\) is the function in \eqref{eq:weyl-kernel}.

For \(f\in\mathcal S(\mathbb R^n\times\mathbb R^n)\), write \(\mathcal F_t f\) for the partial Fourier transform \eqref{E107} in the second variable, so that
\begin{equation}
(\mathcal F_tf)(s,q)
=\int_{\mathbb R^n}f(s,t)e^{-it\cdot q}\,dt,\quad s,q\in\mathbb R^n.
\label{E53}
\end{equation}
Using \eqref{E53} and the coordinates \(s=(x-y)/\sqrt2\), \(q=(x+y)/\sqrt2\), we rewrite \eqref{eq:weyl-kernel} as
$$
K_U\!\left(\frac{q+s}{\sqrt2},\frac{q-s}{\sqrt2}\right)
=2^{-n/2}(\mathcal F_tU)(s,q).
$$
By \eqref{E107} with \(s,t\in\mathbb R^n\) and transformation in \(t\), \(\mathcal F_t\) is an isomorphism of \(\mathcal S(\mathbb R^{2n})\). Define \(J:\mathbb R^{2n}\to\mathbb R^{2n}\) by \(J(s,q)=((q+s)/\sqrt2,(q-s)/\sqrt2)\), and let \(B:\mathcal S(\mathbb R^{2n})\to\mathcal S(\mathbb R^{2n})\) be given by \((B\varphi)(s,q)=\varphi(J(s,q))\). This \(J\) is orthogonal with absolute Jacobian one, and \(B=J^*\) is an isomorphism by \eqref{E107}. For \(U\in\mathcal S(\mathbb R^{2n})\), the coordinate formula for \(K_U\) thus gives
\begin{equation*}
K_U=2^{-n/2}B^{-1}\mathcal F_tU,\qquad
U=2^{n/2}\mathcal F_t^{-1}BK_U.
\tag{\getrefnumber{E53}a}\label{E124}
\end{equation*}
Both maps are continuous in the Schwartz topology by \eqref{E107}.

To extend the kernel correspondence \eqref{E124} to tempered distributions, we transpose the partial Fourier transform and the change of coordinates on the test space. Let \(A\) be a continuous linear isomorphism of \(\mathcal S(\mathbb R^{2n})\). The transpose \(A^t:\mathcal S'(\mathbb R^{2n})\to\mathcal S'(\mathbb R^{2n})\) and its inverse are given by
\begin{equation}
\langle A^tT,\varphi\rangle=\langle T,A\varphi\rangle,
\qquad (A^t)^{-1}=(A^{-1})^t
\label{E74}
\end{equation}
for \(T\in\mathcal S'(\mathbb R^{2n})\) and \(\varphi\in\mathcal S(\mathbb R^{2n})\). Since \(T\circ A\) is a continuous linear functional, \(A^tT\) is tempered. For a bounded set \(D\) of Schwartz tests, the seminorm \(\nu_D\) of the strong dual satisfies
\begin{equation*}
\nu_D(T):=\sup_{\varphi\in D}|\langle T,\varphi\rangle|,
\qquad \nu_D(A^tT)=\nu_{A(D)}(T).
\tag{\getrefnumber{E74}a}\label{E128}
\end{equation*}
Since \(A\) and \(A^{-1}\) take bounded sets to bounded sets, \eqref{E128} proves continuity of both transposes in \eqref{E74}. This construction uses the complex-linear pairing \eqref{E43}. For \(U\in\mathcal S'(\mathbb R^{2n})\) and \(\psi\in\mathcal S(\mathbb R^{2n})\), apply \eqref{E74} to Fourier transformation of \(\psi(s,q)\) in \(q\), with output variable \(t\). In the convention \eqref{E15}, this defines
\begin{equation}
\langle\mathcal F_tU,\psi\rangle
=\bigg\langle U(s,t),\int_{\mathbb R^n}e^{-it\cdot q}\psi(s,q)\,dq\bigg\rangle.
\label{E54}
\end{equation}

For Schwartz \(U\) and \(\psi\), \eqref{E44} turns the right side of \eqref{E54} into the integral of \(U(s,t)e^{-it\cdot q}\psi(s,q)\). Schwartz decay makes this integral absolutely convergent, so integrating first in \(t\) recovers the pairing of the function in \eqref{E53} with \(\psi\). By \eqref{E74}, this extension is an isomorphism of \(\mathcal S'(\mathbb R^{2n})\). For \(U\in\mathcal S'(\mathbb R^{2n})\), set \(K_U=2^{-n/2}B^t\mathcal F_tU\). For every \(\varphi\in\mathcal S(\mathbb R^{2n})\), \eqref{E74} gives
\begin{equation*}
\langle K_U,\varphi\rangle
=2^{-n/2}\bigg\langle\mathcal F_tU(s,q),
\varphi\!\bigg(\frac{q+s}{\sqrt2},\frac{q-s}{\sqrt2}\bigg)\bigg\rangle.
\tag{\getrefnumber{E54}a}\label{E110}
\end{equation*}
Since \(J\) has absolute Jacobian one, \(B^tf=f\circ J^{-1}\) on Schwartz functions. Thus \eqref{E110} recovers the pointwise formula \eqref{E124} for Schwartz symbols and specifies the distributional meaning of \eqref{eq:weyl-kernel} for general \(U\). Using the inverse formula in \eqref{E74}, we recover \(U\) from \(K_U\) by
$$
U=2^{n/2}\mathcal F_t^{-1}(B^{-1})^tK_U
$$
and hence \(U\mapsto K_U\) is an isomorphism from \(\mathcal S'(\mathbb R^{2n})\) onto \(\mathcal S'(\mathbb R^n\times\mathbb R^n)\). To extend the kernel identity \eqref{E115} to tempered \(U\), take \(u,v\in\mathcal S(\mathbb R^n)\) and set \(\varphi(x,y)=u(y)\overline{v(x)}\) in \eqref{E110}. By \eqref{E54}, the test paired with \(U(s,t)\), including the factor \(2^{-n/2}\), is
$$
2^{-n/2}\int_{\mathbb R^n}e^{-it\cdot q}\,
u\!\left(\frac{q-s}{\sqrt2}\right)
\overline{v\!\left(\frac{q+s}{\sqrt2}\right)}\,dq.
$$
With \(q=\sqrt2x-s\), the arguments of \(u\) and \(v\) become \(x-\sqrt2s\) and \(x\), respectively, and \(e^{-it\cdot q}=e^{is\cdot t}e^{-i\sqrt2t\cdot x}\). The Jacobian \(dq=2^{n/2}\,dx\) cancels the prefactor, so \eqref{E25} identifies this test as the coefficient \((s,t)\mapsto\langle W_{(s,t)}u,v\rangle_2\). Substitution in \eqref{E110} yields
\begin{equation*}
\langle K_U(x,y),u(y)\overline{v(x)}\rangle
=\big\langle U(s,t),\langle W_{(s,t)}u,v\rangle_2\big\rangle
=\langle W(U)u,\overline v\rangle.
\tag{\getrefnumber{E54}b}\label{E116}
\end{equation*}
The last equality is \eqref{eq:weyl-map}, so \(K_U\) is the kernel of \(W(U)\). Uniqueness in the tempered Schwartz kernel theorem \cite[Corollary to Theorem 51.6]{Treves1967}, together with the injectivity of \(U\mapsto K_U\), gives
$$
W(U)=0\quad\Longrightarrow\quad K_U=0\quad\Longrightarrow\quad U=0.
$$
We conclude that
\begin{equation}
U\mapsto W(U)\quad\hbox{is injective}.
\label{E40}
\end{equation}
This completes the proof.
\end{proof}

For a composition \(W(P)W(Q)\) to act on \(\mathcal S(\mathbb R^n)\), the output of \(W(Q)\) must lie in the domain of \(W(P)\). Compact support of \(Q\) ensures that \(W(Q)\) preserves \(\mathcal S(\mathbb R^n)\).

\begin{lemma}\label{L4}

Let \(P,Q\in\mathcal S'(\mathbb R^{2n})\) be compactly supported. Then \(W(P)\) and \(W(Q)\) map \(\mathcal S(\mathbb R^n)\) continuously to itself and satisfy
\begin{equation}
W(P\star Q)u=W(P)W(Q)u
\label{eq:compact-composition}
\end{equation}
for every \(u\in\mathcal S(\mathbb R^n)\).

\end{lemma}

\begin{proof}

We first show that the tempered distribution \(W(P)u\) given by Definition \ref{def:weyl-map} is a Schwartz function and estimate its seminorms in terms of those of \(u\). Fix a compact neighborhood \(K\) of \(\operatorname{supp}P\). Pairing with \(P\) requires control of derivatives in \(a\) on \(K\), while the Schwartz seminorms of the output involve weights and derivatives in \(x\). For \(a=(s,t)\in\mathbb R^{2n}\), \(u\in\mathcal S(\mathbb R^n)\), and \(j\in\{1,\ldots,n\}\), differentiating \eqref{E24} gives the identities
\begin{equation}
\begin{aligned}
\partial_{s_j}(W_a u)
&=it_jW_a u-\sqrt2W_a(\partial_{x_j}u),\\
\partial_{t_j}(W_a u)
&=i(s_j-\sqrt2x_j)W_a u,\\
\partial_{x_j}(W_a u)
&=-i\sqrt2t_jW_a u+W_a(\partial_{x_j}u).
\end{aligned}
\label{E84}
\end{equation}

We use \(\alpha\in\mathbb N_0^{2n}\) for parameter derivatives and \(\beta\in\mathbb N_0^n\) for spatial derivatives. Fix \(N,k,m\in\mathbb N_0\). For \(|\alpha|\le m\) and \(|\beta|\le k\), we iterate \eqref{E84} to expand each mixed derivative in the form
$$
\partial_x^\beta\partial_a^\alpha(W_a u)(x)
=\sum_{r=1}^{J}q_{a,r}(x)(W_a(\partial^{\gamma_r}u))(x),
$$
where \(J\in\mathbb N\) is the number of terms, \(\gamma_r\in\mathbb N_0^n\), \(|\gamma_r|\le k+m\), and \(\deg q_{a,r}\le m\). Only a \(t\)-derivative in \eqref{E84} can raise the degree in \(x\), and each differentiation introduces at most one further derivative of \(u\). These facts give the bounds \(m\) and \(k+m\), respectively. The dependence of \(J,\gamma_r,q_{a,r}\) on \((\alpha,\beta)\) is suppressed. Write \(q_{a,r}(x)=\sum_{|\nu|\le m}b_{r,\nu}(a)x^\nu\), with \(\nu\in\mathbb N_0^n\) and coefficients polynomial in \(a\). These coefficients are bounded on \(K\), and only finitely many occur as \(\alpha,\beta\) vary. We may therefore choose \(C_1\) such that
\begin{equation}
\sup_{a\in K}\sum_{r=1}^{J}\sum_{|\nu|\le m}|b_{r,\nu}(a)|\le C_1,
\qquad C_1=C_1(n,K,k,m)>0,
\label{E83}
\end{equation}
uniformly for \(|\alpha|\le m\) and \(|\beta|\le k\). Choose also \(C_K\ge1\), depending only on \(K\), such that \(\sqrt2|s|\le C_K-1\) for \(a=(s,t)\in K\). The bounds \(|x^\nu|\le\langle x\rangle^m\), \(\langle x-\sqrt2s\rangle\ge1\), and the triangle inequality in \(\mathbb R^{n+1}\) give
\begin{equation}
\begin{aligned}
&|q_{a,r}(x)|\le\biggl(\sum_{|\nu|\le m}|b_{r,\nu}(a)|\biggr)\langle x\rangle^m,\\
&\langle x\rangle\le\langle x-\sqrt2s\rangle+\sqrt2|s|
\le C_K\langle x-\sqrt2s\rangle.
\end{aligned}
\label{E82}
\end{equation}

Set \(C_{K,N,k,m}=C_1C_K^{N+m}\). In the derivative expansion the phase in \eqref{E24} has modulus one, and \eqref{E83} bounds the coefficient sum uniformly on \(K\). Using \eqref{E82} to bound the monomials and transfer the weight to \(x-\sqrt2s\), we obtain from \eqref{E51} that
$$
\begin{aligned}
\langle x\rangle^N|\partial_x^\beta\partial_a^\alpha(W_a u)(x)|
&\le C_1\langle x\rangle^{N+m}
\max_{|\gamma|\le k+m}|\partial^\gamma u(x-\sqrt2s)|\\
&\le C_{K,N,k,m}\langle x-\sqrt2s\rangle^{N+m}
\max_{|\gamma|\le k+m}|\partial^\gamma u(x-\sqrt2s)|\\
&\le C_{K,N,k,m}p_{N+m,k+m}(u)
\end{aligned}
$$
for \(x\in\mathbb R^n\), \(a\in K\), \(|\alpha|\le m\), and \(|\beta|\le k\). Taking the suprema in \(x,a\) and the maxima in \(\alpha,\beta\) yields
\begin{equation}
\max_{|\alpha|\le m}\sup_{a\in K}
p_{N,k}(\partial_a^\alpha W_a u)
\le C_{K,N,k,m}p_{N+m,k+m}(u).
\label{E26}
\end{equation}

By restriction, \(P\in\mathcal S'(\mathbb R^{2n})\) also belongs to \(\mathcal D'(\mathbb R^{2n})\), with the same compact support. We combine its finite-order estimate with \eqref{E26} to construct a Schwartz representative of \(W(P)u\). Write \(\operatorname{int}K\) for the interior of \(K\) and choose \(\rho\in C_c^\infty(\operatorname{int}K)\) equal to one near \(\operatorname{supp}P\). Formula \eqref{E56} defines the action of \(P\) on \(C^\infty(\mathbb R^{2n})\). With \(\rho\) fixed, \cite[(2.3.1), p. 44]{Hormander2003} gives \(m\in\mathbb N_0\) and \(C_2>0\), depending on \(P,\rho\), such that every \(\psi\in C^\infty(\mathbb R^{2n})\) satisfies
\begin{equation}
\begin{aligned}
|\langle P,\psi\rangle|
&\le C_2\sum_{|\alpha|\le m}\sup_{a\in\operatorname{supp}\rho}
|\partial_a^\alpha\psi(a)|\le C_P\max_{|\alpha|\le m}\sup_{a\in K}
|\partial_a^\alpha\psi(a)|.
\end{aligned}
\label{E57}
\end{equation}
In the second inequality, we use \(\operatorname{supp}\rho\subset K\) and absorb the finite number of terms into \(C_P>0\).
We now fix the order \(m\) supplied by this estimate. For \(u\in\mathcal S(\mathbb R^n)\), \eqref{E24} shows that \((a,x)\mapsto(W_a u)(x)\) is jointly smooth on \(\mathbb R^{2n}\times\mathbb R^n\). Thus, \eqref{E56} defines \(f_u:\mathbb R^n\to\mathbb C\) by
\begin{equation}
f_u(x)=\langle P(a),(W_a u)(x)\rangle=\langle P(a),\rho(a)(W_a u)(x)\rangle.
\label{E58}
\end{equation}
Set \(\phi_u(a,x)=\rho(a)(W_a u)(x)\) for \((a,x)\in\mathbb R^{2n}\times\mathbb R^n\). Then
\begin{equation}
\phi_u\in C^\infty(\mathbb R^{2n}\times\mathbb R^n),\qquad
\operatorname{supp}\phi_u(\cdot,x)\subset\operatorname{supp}\rho
\quad(x\in\mathbb R^n).
\label{E91}
\end{equation}
To differentiate \eqref{E58}, let \(T\in\mathcal D'(\mathbb R^p)\) and \(F\in C^\infty(\mathbb R^p\times O)\), where \(O\subset\mathbb R^q\) is open and \(p,q\ge1\) are integers. Suppose that each \(z_0\in O\) has an open neighborhood \(O_0\subset O\) and a compact set \(K_0\subset\mathbb R^p\) such that \(\operatorname{supp}F(\cdot,z)\subset K_0\) for all \(z\in O_0\). On each \(O_0\), the same compact set \(K_0\) works for every parameter, so \cite[Theorem 2.1.3]{Hormander2003} applies there. Since these neighborhoods cover \(O\), the function \(z\mapsto\langle T,F(\cdot,z)\rangle\) is smooth on \(O\), with
\begin{equation*}
\partial_z^\beta\langle T,F(\cdot,z)\rangle
=\langle T,\partial_z^\beta F(\cdot,z)\rangle,
\qquad z\in O,\quad\beta\in\mathbb N_0^q.
\tag{\getrefnumber{E91}a}\label{E108}
\end{equation*}
By \eqref{E91}, formula \eqref{E108} applies to \(T=P\) and \(F(a,x)=\phi_u(a,x)\) on \(O=\mathbb R^n\), with \(\operatorname{supp}\rho\) as the fixed compact set. It follows that \(f_u\) is smooth. Differentiation in \(x\) leaves the cutoff \(\rho\) unchanged, so
$$
\begin{aligned}
\partial_x^\beta f_u(x)
=\langle P(a),\rho(a)\partial_x^\beta(W_a u)(x)\rangle=\langle P(a),\partial_x^\beta(W_a u)(x)\rangle
\end{aligned}
$$
for \(x\in\mathbb R^n\) and \(\beta\in\mathbb N_0^n\). The last equality uses \eqref{E56}. Substituting this derivative formula into the seminorm \eqref{E51} for \(N,k\in\mathbb N_0\), we obtain
\begin{equation}
\begin{aligned}
p_{N,k}(f_u)
&=\max_{|\beta|\le k}\sup_{x\in\mathbb R^n}\langle x\rangle^N
\bigl|\langle P(a),\partial_x^\beta(W_a u)(x)\rangle\bigr|\\
&\le C_P\max_{|\alpha|\le m}\sup_{a\in K}
p_{N,k}(\partial_a^\alpha W_a u)\\
&\le C_PC_{K,N,k,m}p_{N+m,k+m}(u)<\infty,
\end{aligned}
\label{E85}
\end{equation}
where the inequalities follow successively from \eqref{E57} and \eqref{E26} with this fixed \(m\). Since \(f_u\) is smooth, these bounds and \eqref{E51} show that \(f_u\in\mathcal S(\mathbb R^n)\).

To identify \(f_u\) with \(W(P)u\), fix \(v\in\mathcal S(\mathbb R^n)\). By \eqref{E58} and \eqref{eq:weyl-map}, it suffices to pass \(P\) through the integral of \((W_a u)(x)\overline{v(x)}\) in \(x\). On a fixed ball this follows by approximating the integral by Riemann sums. To remove the truncation, \eqref{E57} requires a tail estimate uniform on \(K\) through \(m\) derivatives in \(a\). For this estimate, take \(N=n+1\) and \(k=0\) in \eqref{E26} to obtain
\begin{equation}
\begin{aligned}
&\quad\max_{|\alpha|\le m}\sup_{a\in K}
|\partial_a^\alpha\bigl((W_a u)(x)\overline{v(x)}\bigr)|\\
&\le C_{K,n+1,0,m}p_{n+1+m,m}(u)\|v\|_\infty
\langle x\rangle^{-n-1}.
\end{aligned}
\label{E61}
\end{equation}
For \(R>0\) and \(a\in\mathbb R^{2n}\), set
$$
\begin{aligned}
&H_R(a)=\int_{|x|\le R}(W_a u)(x)\overline{v(x)}\,dx,\\
&H(a)=\int_{\mathbb R^n}(W_a u)(x)\overline{v(x)}\,dx
=\langle W_a u,v\rangle_2.
\end{aligned}
$$
As the coefficient function in \eqref{E88}, \(H\) belongs to \(\mathcal S(\mathbb R^{2n})\) and is the test function paired with \(P\) in \eqref{eq:weyl-map}. The bounds in \eqref{E26}, applied at every parameter derivative order on compact parameter sets, justify differentiating both integrals in \(a\). In particular, \(H_R\) is smooth. Applying \eqref{E61} to the integral over \(|x|>R\), we can estimate the derivatives of \(H-H_R\) by
$$
\begin{aligned}
&\quad\max_{|\alpha|\le m}\sup_{a\in K}
|\partial_a^\alpha(H-H_R)(a)|\le\int_{|x|>R}\max_{|\alpha|\le m}\sup_{a\in K}
\left|\partial_a^\alpha\bigl((W_a u)(x)\overline{v(x)}\bigr)\right|\,dx\\
&\le C_{K,n+1,0,m}p_{n+1+m,m}(u)\|v\|_\infty
\int_{|x|>R}\langle x\rangle^{-n-1}\,dx\longrightarrow0
\end{aligned}
$$
as \(R\to\infty\). With the pairings against \(P\) interpreted by \eqref{E56}, estimate \eqref{E57} yields
\begin{equation}\label{E90}
|\langle P,H_R\rangle-\langle P,H\rangle|
\le C_P\max_{|\alpha|\le m}\sup_{a\in K}
|\partial_a^\alpha(H_R-H)(a)|\rightarrow0
\end{equation}
as $R\to\infty.$

It remains to identify \(\langle P,H_R\rangle\) for fixed \(R>0\) before taking the limit in \eqref{E90}. Put \(B_R=\{x\in\mathbb R^n:|x|\le R\}\) and take finite measurable partitions \(\{D_{N,j}\}_{j=1}^{J_N}\) of \(B_R\). Writing \(\operatorname{diam}D_{N,j}=\sup\{|x-y|:x,y\in D_{N,j}\}\) for the diameter of a cell, we choose the partitions so that their mesh \(\delta_N=\max_j\operatorname{diam}D_{N,j}\) tends to zero as \(N\to\infty\). In each cell \(D_{N,j}\) of Lebesgue measure \(|D_{N,j}|\), fix a point \(x_{N,j}\). To pass the compact integral through \(P\), we use the convergence
$$
S_N=\sum_{j=1}^{J_N}|D_{N,j}|\phi_u(\cdot,x_{N,j})\overline{v(x_{N,j})}
\longrightarrow\rho H_R\quad\text{in }\mathcal D(\mathbb R^{2n}).
$$
Indeed, put \(F(a,x)=\phi_u(a,x)\overline{v(x)}\). For each \(\alpha\in\mathbb N_0^{2n}\), joint smoothness and the common support in \eqref{E91} give a finite bound \(C_{\alpha,R}\) for \(|\nabla_x\partial_a^\alpha F|\) on \(\mathbb R^{2n}\times B_R\). For \(x\in D_{N,j}\), the segment from \(x_{N,j}\) to \(x\) lies in \(B_R\) and its length satisfies \(|x-x_{N,j}|\le\operatorname{diam}D_{N,j}\le\delta_N\). We apply the mean value estimate along each segment and integrate over its cell to bound the error by
\begin{equation*}
\sup_a|\partial_a^\alpha(S_N-\rho H_R)(a)|
\le\sum_{j=1}^{J_N}|D_{N,j}|C_{\alpha,R}\delta_N
=|B_R|C_{\alpha,R}\delta_N\rightarrow0.
\tag{\getrefnumber{E90}a}\label{E94}
\end{equation*}
Estimate \eqref{E94} holds for every derivative, and \(S_N-\rho H_R\) is supported in the fixed compact set \(\operatorname{supp}\rho\). Thus \(S_N\to\rho H_R\) in \(\mathcal D(\mathbb R^{2n})\).

To pass from test functions to scalar integrals, let \(G:B\to\mathcal D(\mathbb R^q)\) be continuous on a closed ball \(B\subset\mathbb R^p\), where \(p,q\ge1\) are integers. Take tagged partitions \(\{D_{N,j}\}_j\) of \(B\) with \(x_{N,j}\in D_{N,j}\) and mesh tending to zero as \(N\to\infty\). If \(I_N=\sum_j|D_{N,j}|G(x_{N,j})\) converges to \(I\) in \(\mathcal D(\mathbb R^q)\), then every \(T\in\mathcal D'(\mathbb R^q)\) satisfies
\begin{equation*}
\begin{aligned}
&\langle T,I_N\rangle=\sum_j|D_{N,j}|\langle T,G(x_{N,j})\rangle,\\
&\langle T,I\rangle=\lim_{N\to\infty}\langle T,I_N\rangle
=\int_B\langle T,G(x)\rangle\,dx.
\end{aligned}
\tag{\getrefnumber{E90}b}\label{E95}
\end{equation*}
The finite sum identity follows from linearity. Its left side tends to \(\langle T,I\rangle\) by continuity of \(T\), and its right side is a Riemann sum for the continuous scalar function \(x\mapsto\langle T,G(x)\rangle\) on \(B\).

By \eqref{E91}, joint smoothness and common compact support make \(F(\cdot,x)\) continuous in \(x\in B_R\) with values in \(\mathcal D(\mathbb R^{2n})\). Thus \eqref{E94} allows us to take \(T=P\), \(G(x)=F(\cdot,x)\) and \(I=\rho H_R\) in \eqref{E95}. Using \eqref{E56} for \(H_R\) and \eqref{E58} for \(f_u\), we obtain
$$
\begin{aligned}
\langle P,H_R\rangle
&=\langle P,\rho H_R\rangle
=\lim_{N\to\infty}\langle P,S_N\rangle=\int_{|x|\le R}\langle P,\phi_u(\cdot,x)\rangle\overline{v(x)}\,dx\\
&=\int_{|x|\le R}f_u(x)\overline{v(x)}\,dx.
\end{aligned}
$$
We now let \(R\to\infty\) in this identity for \(B_R\). On the integral side, \(f_u\overline v\in L^1(\mathbb R^n)\) and \eqref{E44} identifies the limit. On the distribution side, \eqref{E90} gives \(\langle P,H\rangle\). Substituting the coefficient \(H\) in \eqref{eq:weyl-map}, we obtain
$$
\begin{aligned}
\langle f_u,\overline v\rangle
&=\lim_{R\to\infty}\int_{|x|\le R}f_u(x)\overline{v(x)}\,dx=\lim_{R\to\infty}\langle P,H_R\rangle=\langle P,H\rangle
=\langle W(P)u,\overline v\rangle.
\end{aligned}
$$
Since \(v\) is arbitrary, \(W(P)u=f_u\) as tempered distributions, and \eqref{E85} now reads
\begin{equation}
p_{N,k}(W(P)u)\le C_PC_{K,N,k,m}p_{N+m,k+m}(u)
\label{E6}
\end{equation}
for all \(N,k\in\mathbb N_0\) and \(u\in\mathcal S\), proving continuity of \(W(P):\mathcal S\to\mathcal S\).

To prove \eqref{eq:compact-composition}, fix \(u,v\in\mathcal S(\mathbb R^n)\). Choose \(R>0\) so that \(|a|\le R\) on \(\operatorname{supp}P\cup\operatorname{supp}Q\). Then \(0\ge |a|-R\) on both supports, so \(P,Q\in\mathcal C_{0,1}\) by \eqref{eq:cone}. By \eqref{eq:product-support}, \(\operatorname{supp}(P\star Q)\) lies in the compact sum \(\operatorname{supp}P+\operatorname{supp}Q\). Thus \eqref{E6} applies to \(Q\) and \(P\star Q\) as well as to \(P\), and both \(W(P)W(Q)u\) and \(W(P\star Q)u\) belong to \(\mathcal S(\mathbb R^n)\). For each fixed \(a\in\mathbb R^{2n}\), \eqref{eq:weyl-relations} identifies \(W_a^*v\) with \(W_{-a}v\), which is Schwartz by \eqref{E24}. Taking \(W_a^*v\) as the second vector in \eqref{eq:weyl-map} for \(Q\), we rewrite the coefficient \(b\mapsto\langle W_bu,W_a^*v\rangle_2\), \(b\in\mathbb R^{2n}\), by the adjoint identity to obtain
\begin{equation}
\begin{aligned}
\langle W_aW(Q)u,v\rangle_2
&=\langle W(Q)u,W_a^*v\rangle_2=\big\langle Q(b),\langle W_aW_bu,v\rangle_2\big\rangle.
\end{aligned}
\label{E93}
\end{equation}

Write \(c_{u,v}(z)=\langle W_z u,v\rangle_2\), which is Schwartz by \eqref{E88}. For the Schwartz function \(W(P)W(Q)u\), the pairing against \(\overline v\) equals the \(L^2\) inner product with \(v\) by \eqref{E44}. Taking \(W(Q)u\) as input in \eqref{eq:weyl-map}, we may substitute \eqref{E93} into the pairing with \(P\) to obtain
$$
\begin{aligned}
\langle W(P)W(Q)u,v\rangle_2
&=\langle W(P)W(Q)u,\overline v\rangle=\big\langle P(a),\langle W_aW(Q)u,v\rangle_2\big\rangle\\
&=\Big\langle P(a),\big\langle Q(b),\langle W_aW_bu,v\rangle_2\big\rangle\Big\rangle\\
&=\big\langle P(a)\otimes Q(b),e^{i\sigma(a,b)}c_{u,v}(a+b)\big\rangle.
\end{aligned}
$$
The phase in the last equality comes from \eqref{eq:weyl-relations}. To apply \eqref{E36} to the smooth function \(e^{i\sigma(a,b)}c_{u,v}(a+b)\), choose \(\chi_P,\chi_Q\in C_c^\infty(\mathbb R^{2n})\) equal to one near \(\operatorname{supp}P\) and \(\operatorname{supp}Q\), respectively. By \eqref{E56} we may insert \(\chi_P(a)\chi_Q(b)\) in both pairings, so \eqref{E36} applies to a compactly supported test. The cutoffs in \eqref{eq:distribution-product} also equal one near \(\operatorname{supp}P\times\operatorname{supp}Q\), identifying the tensor pairing with \(\langle P\star Q,c_{u,v}\rangle\). Since \(W(P\star Q)u\) is Schwartz, we can use \eqref{E44} to write its defining pairing \eqref{eq:weyl-map} in the form
$$
\langle P\star Q,c_{u,v}\rangle
=\langle W(P\star Q)u,v\rangle_2.
$$
Since \(v\) was arbitrary, the two Schwartz functions \(W(P)W(Q)u\) and \(W(P\star Q)u\) agree, proving \eqref{eq:compact-composition}.
\end{proof}

\subsection{Weyl smoothing and the proof of Theorem \ref*{thm:distribution}}\label{sec:weyl-smoothing}

The order of the factors in \(h=V\star g\star U\) places \(U\star V\) in the middle of \(h\star h\). To apply Proposition \ref{prop:smooth-square}, we first establish smoothness and the support bound for \(h\), then show that for each pair of nonzero factors \(U,V\) one can choose \(g\) so that \(h\ne0\).

\begin{lemma}\label{lem:smooth-regularization}

Let \(U,V\in\mathcal C_{e,\eta}\) and \(g\in C_c^\infty(\mathbb R^{2n})\). Then \(h=V\star g\star U\) is a smooth function in \(\mathcal C_{e,\eta}\).

\end{lemma}

\begin{proof}

Recall that \(d=2n\), so \(g\in C_c^\infty(\mathbb R^d)\). Choose \(L>0\) with \(\operatorname{supp}g\subset\{a\in\mathbb R^d:|a|\le L\}\), and put \(M_g=(\eta+|e|)L\ge0\). For \(a\in\operatorname{supp}g\), we have
\begin{equation*}
e\cdot a\ge-|e|\,|a|
=\eta|a|-(\eta+|e|)|a|
\ge\eta|a|-M_g.
\tag{\getrefnumber{E93}a}\label{E109}
\end{equation*}
Since \(g\in\mathcal S'(\mathbb R^d)\) by \eqref{E18}, this estimate verifies \(g\in\mathcal C_{e,\eta}\). Applying Proposition \ref{prop:tempered-product} first to \(V,g\) and then to \(V\star g,U\), and using its associativity, we obtain \(h_0=V\star g\star U\in\mathcal C_{e,\eta}\).

Use the support constants \(M_U,M_V\) from \eqref{eq:cone} and \(M_g\) from \eqref{E109} to choose \(\rho_U,\rho_V,\rho_g\) as in \eqref{E101}. In \eqref{E92} for \((V,g,U)\), \eqref{E44} expresses the action of \(g\) as integration in \(r\). Its cutoff drops out because \eqref{E102} gives \(\rho_g=1\) near \(\operatorname{supp}g\), hence \(\rho_g g=g\). Writing \(r=z-y-x\) in \(g\) and \(y+r=z-x\) in the phase of \eqref{E104} suggests a candidate for the smooth representative of \(h_0\). Its test in \((y,x)\) is the function \(\varphi_z:\mathbb R^d\times\mathbb R^d\to\mathbb C\) given by
$$
\begin{aligned}
\varphi_z(y,x)
&=\rho_V(y)\rho_U(x)g(z-y-x) e^{i\sigma(y,z-y-x)+i\sigma(z-x,x)}.
\end{aligned}
$$
For fixed \(z\), the factor \(g(z-y-x)\) confines \(x+y\) to the compact set \(\{z-r:r\in\operatorname{supp}g\}\). On the two cutoff supports, \eqref{eq:cone-control} then bounds both \(x\) and \(y\). Hence \(\varphi_z\in C_c^\infty(\mathbb R^{2d})\), and we can define
\begin{equation}
h(z)=\langle V(y)\otimes U(x),\varphi_z(y,x)\rangle.
\label{eq:smooth-formula}
\end{equation}
To differentiate this pairing, we need a common compact support as \(z\) varies on a ball. On \(\operatorname{supp}\rho_V\times\operatorname{supp}\rho_U\), the condition \(r=z-y-x\in\operatorname{supp}g\) implies \(|r|\le L\). Applying \eqref{eq:cone-control} with \(x+y=z-r\) gives \(\frac\eta2(|x|+|y|)\le |e|(|z|+L)+L_U+L_V\), so there is a constant \(C_g>0\), depending only on \(e,\eta,M_U,M_V,L\), such that
\begin{equation}
|x|+|y|\le C_g(1+|z|).
\label{eq:smoothing-support}
\end{equation}
Since \(\operatorname{supp}\partial^\alpha g\subset\operatorname{supp}g\), the same bound holds on the support of every \(z\)-derivative of \(\varphi_z\). Fix \(z_0\in\mathbb R^d\). For \(z\in B(z_0,1)\), all these supports lie in the compact set
$$
\{(y,x)\in\mathbb R^{2d}:|x|+|y|\le C_g(2+|z_0|)\}.
$$
Thus \eqref{E108} applies on \(B(z_0,1)\) with \(T=V\otimes U\) acting in \((y,x)\) on \(F((y,x),z)=\varphi_z(y,x)\). These balls cover \(\mathbb R^d\), so the function \(h\) defined by \eqref{eq:smooth-formula} is smooth everywhere and satisfies
$$
\partial_z^\alpha h(z)
=\langle V(y)\otimes U(x),\partial_z^\alpha\varphi_z(y,x)\rangle,
\quad\alpha\in\mathbb N_0^d.
$$
To see that integration against this smooth function defines a tempered distribution, choose \(M,k\in\mathbb N_0\) and \(C>0\) such that \(|\langle V\otimes U,\Phi\rangle|\le C p_{M,k}^{(y,x)}(\Phi)\). On the support of \(\varphi_z\) and its derivatives, \eqref{eq:smoothing-support} bounds \(\langle(y,x)\rangle\) by a constant times \(\langle z\rangle\). With \(r=z-y-x\), bilinearity and alternation in \eqref{E55} give \(\sigma(y,r)+\sigma(y+r,x)=\sigma(y,z)-\sigma(y,x)+\sigma(z,x)\). Each differentiation of the exponential factor in \((y,x)\) therefore contributes at most one power of \(\langle z\rangle\) on this support. The derivatives of \(g\) and of the cutoffs are bounded by compact support and \eqref{E45}, respectively. Hence for some \(C'>0\) independent of \(z\), \eqref{eq:smooth-formula} yields
$$
|h(z)|\le C p_{M,k}^{(y,x)}(\varphi_z)\le C'\langle z\rangle^{M+k}.
$$
In particular, the integral on each Schwartz test \(\psi\in\mathcal S(\mathbb R^d)\) satisfies
\begin{equation*}
\int_{\mathbb R^d}|h(z)\psi(z)|\,dz
\le C'\|\langle\cdot\rangle^{-d-1}\|_1\,p_{M+k+d+1,0}(\psi).
\tag{\getrefnumber{eq:smoothing-support}a}\label{E130}
\end{equation*}

We identify \(h\) with \(h_0\) by reducing their actions on each \(\psi\in C_c^\infty(\mathbb R^d)\) to the same pairing with \(V\otimes U\). Fix such a \(\psi\) and choose \(R_\psi>0\) so that the closed ball \(B_\psi=\{z\in\mathbb R^d:|z|\le R_\psi\}\) contains \(\operatorname{supp}\psi\). By \eqref{eq:smoothing-support}, all functions \(\psi(z)\varphi_z\) are supported in the compact set
$$
K_\psi=\bigl\{(y,x)\in\mathbb R^{2d}:|x|+|y|\le C_g(1+R_\psi)\bigr\}.
$$
Since \(\Psi_\psi(z,y,x)=\psi(z)\varphi_z(y,x)\) is smooth with support in \(\operatorname{supp}\psi\times K_\psi\), its integral in \(z\) defines \(f_\psi\in C_c^\infty(\mathbb R^{2d})\) by
$$
f_\psi(y,x)=\int_{\mathbb R^d}\psi(z)\varphi_z(y,x)\,dz.
$$
For a tagged partition of \(B_\psi\) of mesh \(\delta\), the estimate \eqref{E94} in the integration variable \(z\) bounds the error in each derivative \(\partial_{(y,x)}^\alpha\), \(\alpha\in\mathbb N_0^{2d}\), by \(\delta |B_\psi|\|\nabla_z\partial_{(y,x)}^\alpha\Psi_\psi\|_\infty\). As \(\delta\to0\), the bound tends to zero for each \(\alpha\). Because all sums are supported in \(K_\psi\), they converge to \(f_\psi\) in \(\mathcal D(\mathbb R^{2d})\). We may therefore apply \eqref{E95} with \(G(z)=\Psi_\psi(z,\cdot,\cdot)\) and \(T=V\otimes U\). Since \(\langle T,G(z)\rangle=\psi(z)h(z)\) by \eqref{eq:smooth-formula}, the integral representing \(\langle h,\psi\rangle\) in \eqref{E44} satisfies
\begin{equation}
\begin{aligned}
\int_{\mathbb R^d}h(z)\psi(z)\,dz
&=\int_{\mathbb R^d}\langle V\otimes U,\Psi_\psi(z,\cdot,\cdot)\rangle\,dz=\langle V\otimes U,f_\psi\rangle.
\end{aligned}
\label{E62}
\end{equation}
Changing variables from \(z\) to \(r=z-y-x\) in the definition of \(f_\psi\), we obtain
\begin{equation}
\begin{aligned}
f_\psi(y,x)
=\rho_V(y)\rho_U(x)\int_{\mathbb R^d}&g(r)e^{i\sigma(y,r)+i\sigma(y+r,x)}\\
&\quad\cdot
\psi(y+r+x)\,dr.
\end{aligned}
\label{E63}
\end{equation}
Apply \eqref{E104} to \((T_1,T_2,T_3)=(V,g,U)\) in the variables \((a_1,a_2,a_3)=(y,r,x)\). By bilinearity in \eqref{E55}, the sum of the three phase terms is \(\sigma(y,r)+\sigma(y,x)+\sigma(r,x)=\sigma(y,r)+\sigma(y+r,x)\). Thus the Schwartz test for \(\psi\in C_c^\infty(\mathbb R^d)\) is
$$
J_\psi(y,r,x)=\rho_V(y)\rho_g(r)\rho_U(x)
 e^{i\sigma(y,r)+i\sigma(y+r,x)}\psi(y+r+x).
$$
Pairing in \(r\) by \eqref{E44} removes the middle cutoff because \(\rho_g g=g\). The remaining integral is precisely the function of \((y,x)\) in \eqref{E63}, so
\begin{equation*}
\langle g,J_\psi(y,\cdot,x)\rangle
 =\int_{\mathbb R^d}g(r)J_\psi(y,r,x)\,dr=f_\psi(y,x).
\tag{\getrefnumber{E63}a}\label{E106}
\end{equation*}
In the tensor pairing \eqref{E92}, we may act first in \(r\) by \eqref{E36}. With the inner pairing evaluated in \eqref{E106}, the identity becomes
$$
\langle h_0,\psi\rangle
=\langle V\otimes g\otimes U,J_\psi\rangle
=\langle V\otimes U,f_\psi\rangle
=\int_{\mathbb R^d} h(z)\psi(z)\,dz,
$$
where the last equality is \eqref{E62}. By \eqref{E130}, integration against \(h\) is continuous on \(\mathcal S(\mathbb R^d)\). Its agreement with \(h_0\) on compactly supported smooth tests therefore extends by density in \eqref{E18} to every Schwartz test.

Formula \eqref{E37} now identifies the two supports. Applying \eqref{eq:product-support} first to \((V,g)\) and then to \((V\star g,U)\), we obtain
$$
\operatorname{supp}h=\operatorname{supp}h_0
\subset\operatorname{supp}V+\operatorname{supp}g+\operatorname{supp}U.
$$
For \(z=y+r+x\) in this sum, adding the cone bounds for \(y,r,x\) and using \(|z|\le |y|+|r|+|x|\) yields
\begin{equation}
e\cdot z\ge\eta|z|-M_h,\qquad
z\in\operatorname{supp}h,\qquad M_h=M_V+M_g+M_U,
\label{eq:smoothing-cone}
\end{equation}
which completes the proof.
\end{proof}

For nonzero \(U,V\in\mathcal C_{e,\eta}\), it remains to choose \(g\) so that the smooth product \(V\star g\star U\) is nonzero. To use the Weyl maps, the middle factor must take the tempered distribution \(W(U)u\) into the Schwartz domain of \(W(V)\). We obtain this extension from the Schwartz kernel of \(W(g)\), using the seminorm notation introduced after \eqref{E51}.

\begin{lemma}\label{L6}
Let \(g\in\mathcal S(\mathbb R^{2n})\), and write \(\kappa_g(r,s)=K_g(r,s)\) for the kernel in \eqref{eq:weyl-kernel}, with \(r,s\in\mathbb R^n\). The formula
\begin{equation}
(S_gw)(r)=w\bigl(\kappa_g(r,\cdot)\bigr)
\label{E41}
\end{equation}
for \(w\in\mathcal S'(\mathbb R^n)\) and \(r\in\mathbb R^n\) defines a linear map \(S_g:\mathcal S'(\mathbb R^n)\to\mathcal S(\mathbb R^n)\) that agrees with \(W(g)\) on \(\mathcal S(\mathbb R^n)\). For each fixed \(w\), there exist \(M,k\in\mathbb N_0\) and \(C_w>0\), depending only on \(w\), such that
\begin{equation}
p_{N,m}^{(r)}(S_gw)\le C_w p_{N+M,m+k}^{(r,s)}(\kappa_g)
\label{eq:fixed-distribution-smoothing}
\end{equation}
for all \(N,m\in\mathbb N_0\) and \(g\in\mathcal S(\mathbb R^{2n})\).
\end{lemma}

\begin{proof}
By Lemma \ref{L3}, \(\kappa_g\in\mathcal S(\mathbb R^n\times\mathbb R^n)\), so \(\kappa_g(r,\cdot)\in\mathcal S(\mathbb R^n)\) for every \(r\in\mathbb R^n\). The pairing in \eqref{E41} therefore defines \((S_gw)(r)\) and is linear in \(w\). For \(w\in\mathcal S(\mathbb R^n)\), the inclusion in \eqref{E18} allows us to apply \eqref{E41} to \(w\), and \eqref{E44} then gives \((S_gw)(r)=\int_{\mathbb R^n}\kappa_g(r,s)w(s)\,ds\). Taking \(U=g\) and \(u=w\) in \eqref{E115} and integrating first in \(s\), we obtain
$$
\langle W(g)w,\overline v\rangle
=\int_{\mathbb R^n}(S_gw)(r)\overline{v(r)}\,dr
\qquad(v\in\mathcal S(\mathbb R^n)).
$$
The double integral is absolutely convergent since \(\kappa_g,w,v\) are Schwartz. Thus, \(S_gw\) represents \(W(g)w\).

For general \(w\in\mathcal S'(\mathbb R^n)\), continuity in the topology defined by \eqref{E51} provides \(M,k\in\mathbb N_0\) and \(C_w>0\), depending only on \(w\), so that \(|w(\psi)|\le C_w p_{M,k}^{(s)}(\psi)\) for every \(\psi\in\mathcal S(\mathbb R^n)\). In the notation of \eqref{E127}, \(S_gw=F_{\kappa_g}\) with \(w\) acting in \(s\). To pass differentiation in \(r\) through \(w\) in \eqref{E41}, we estimate the remainder of the kernel difference quotient in \(p_{M,k}^{(s)}\). For the \(j\)-th coordinate vector \(e_j\), \(r\in\mathbb R^n\), and real \(h\ne0\), linearity gives
$$
\frac{(S_gw)(r+h e_j)-(S_gw)(r)}{h}
-w\bigl(\partial_{r_j}\kappa_g(r,\cdot)\bigr)
=w\bigl(R_{j,h}(r,\cdot)\bigr),
$$
where, for \(s\in\mathbb R^n\), Taylor's formula with integral remainder reads
$$
\begin{aligned}
R_{j,h}(r,s)
&:=\frac{\kappa_g(r+h e_j,s)-\kappa_g(r,s)}{h}
-\partial_{r_j}\kappa_g(r,s)\\
&=h\int_0^1(1-t)\,\partial_{r_j}^2\kappa_g(r+t h e_j,s)\,dt.
\end{aligned}
$$
Taking \(s\)-derivatives through order \(k\) under the integral and using \(\langle s\rangle\le\langle(r+t h e_j,s)\rangle\) in \eqref{E51}, we obtain
$$
p_{M,k}^{(s)}\bigl(R_{j,h}(r,\cdot)\bigr)
\le \frac{|h|}{2}\,p_{M,k+2}^{(r,s)}(\kappa_g)
\longrightarrow 0
$$
as \(h\to0\), uniformly in \(r\), for every \(M,k\in\mathbb N_0\).

In particular, applying the bound chosen for \(w\) to the kernel remainder gives
$$
|w(R_{j,h}(r,\cdot))|
\le \tfrac12 C_w|h|\,p_{M,k+2}^{(r,s)}(\kappa_g)\longrightarrow0
\qquad(h\to0).
$$
Taking the limit in the difference quotient for \(S_gw\) therefore gives
$$
\partial_{r_j}(S_gw)(r)=w\bigl(\partial_{r_j}\kappa_g(r,\cdot)\bigr).
$$
For every \(\alpha\in\mathbb N_0^n\), the kernel \(\partial_r^\alpha\kappa_g\) is again Schwartz. Applying the first-order calculation to this kernel yields the recursion
$$
\partial_{r_j}\!\left[w\bigl(\partial_r^\alpha\kappa_g(r,\cdot)\bigr)\right]
=w\bigl(\partial_r^{\alpha+e_j}\kappa_g(r,\cdot)\bigr).
$$
Starting from \eqref{E41} and iterating this identity, we obtain
$$
\partial_r^\alpha(S_gw)(r)=w\bigl(\partial_r^\alpha\kappa_g(r,\cdot)\bigr)
\qquad(\alpha\in\mathbb N_0^n).
$$
For arbitrary \(M,k\in\mathbb N_0\), integrating the first \(r\)-derivatives along the segment from \(r\) to \(r'\) gives
$$
\begin{aligned}
&p_{M,k}^{(s)}\bigl(\partial_r^\alpha\kappa_g(r',\cdot)
                 -\partial_r^\alpha\kappa_g(r,\cdot)\bigr)\\
&\qquad\le \sqrt n\,|r'-r|\,p_{M,k+|\alpha|+1}^{(r,s)}(\kappa_g).
\end{aligned}
$$
The derivative slices are therefore continuous in the Schwartz topology, so all the scalar derivatives of \(S_gw\) are continuous. The bound for \(w\) controls \(|\partial_r^\alpha(S_gw)(r)|\) by \(C_wp_{M,k}^{(s)}(\partial_r^\alpha\kappa_g(r,\cdot))\). Multiplying by \(\langle r\rangle^N\) before taking the supremum in \(r\) and the maximum over \(|\alpha|\le m\) in \eqref{E51}, we obtain
$$
\begin{aligned}
p_{N,m}^{(r)}(S_gw)
&\le C_w\max_{\substack{|\alpha|\le m\\|\beta|\le k}}
\sup_{r,s}\langle r\rangle^N\langle s\rangle^M
|\partial_r^\alpha\partial_s^\beta\kappa_g(r,s)|.
\end{aligned}
$$
Since \(\langle r\rangle^N\langle s\rangle^M\le\langle(r,s)\rangle^{N+M}\) and \(|\alpha|+|\beta|\le m+k\), the right-hand side is bounded by \(C_w p_{N+M,m+k}^{(r,s)}(\kappa_g)\). This proves \eqref{eq:fixed-distribution-smoothing} for all \(N,m\in\mathbb N_0\) and hence \(S_gw\in\mathcal S(\mathbb R^n)\). The proof is complete.
\end{proof}

By Lemma \ref{L6}, \(S_g\) takes \(W(U)u\) into the Schwartz space on which \(W(V)\) acts. Expressing \(\langle W(V)S_gW(U)u,\overline v\rangle\) as a pairing with \(U\otimes V\) will allow us to pass from the compact-support formula \eqref{eq:compact-composition} to factors satisfying the common cone condition.

\begin{proposition}\label{lem:sandwich}

Let \(U,V\in\mathcal S'(\mathbb R^{2n})\), and write \(A=W(U)\), \(B=W(V)\). For \(g\in\mathcal S(\mathbb R^{2n})\) and \(u,v\in\mathcal S(\mathbb R^n)\), the function
$$
f_g(x,y)=\langle W_yS_gW_xu,v\rangle_2
$$
with \(x,y\in\mathbb R^{2n}\) belongs to \(\mathcal S(\mathbb R^{4n})\), continuously in \(g,u,\overline v\), and
\begin{equation}
\langle BS_gAu,\overline v\rangle
=\langle U(x)\otimes V(y),f_g(x,y)\rangle.
\label{eq:mixed-pairing}
\end{equation}
If, in addition, \(U,V\in\mathcal C_{e,\eta}\) and \(g\in C_c^\infty(\mathbb R^{2n})\), then
\begin{equation}
W(V\star g\star U)=W(V)S_gW(U)=BS_gA
\label{eq:sandwich}
\end{equation}
as maps \(\mathcal S(\mathbb R^n)\to\mathcal S'(\mathbb R^n)\).

\end{proposition}

\begin{proof}

To construct the Schwartz test \(f_g\) in \eqref{eq:mixed-pairing}, we keep the parameter \(x\) in the smoothed vector \(S_gW_xu\). This will also allow us to pair the vector with \(U\) in \(x\) before taking its Weyl coefficient in \(y\). Write \(x=(a,b)\), \(y=(c,\zeta)\), with \(a,b,c,\zeta\in\mathbb R^n\), and let \(\kappa=\kappa_g\). For fixed \(x\), \(W_xu\) is Schwartz by \eqref{E24}, so the pairing in \eqref{E41} has the integral representation \eqref{E44}. We denote the resulting function by \(q_g:\mathbb R^{2n}\times\mathbb R^n\to\mathbb C\) and substitute the formula for \(W_xu\) to obtain
\begin{equation*}
\begin{aligned}
q_g(x,r)&=(S_gW_xu)(r)=\int_{\mathbb R^n}\kappa_g(r,s)(W_xu)(s)\,ds\\
&=e^{ia\cdot b}\int_{\mathbb R^n}
\kappa_g(r,s)u(s-\sqrt2a)e^{-i\sqrt2b\cdot s}\,ds.
\end{aligned}
\tag{\getrefnumber{eq:sandwich}a}\label{E117}
\end{equation*}
Set \(G(r,a,s)=\kappa_g(r,s)u(s-\sqrt2a)\) for \(r,a,s\in\mathbb R^n\). By \eqref{E107}, the invertible substitution \((r,a,s)\mapsto(r,s,s-\sqrt2a)\) takes \(\kappa_g\otimes u\) to \(G\in\mathcal S(\mathbb R^{3n})\). In the Fourier notation of \eqref{E107}, with \((r,a)\) fixed, \eqref{E117} reads
$$
q_g((a,b),r)=e^{ia\cdot b}(\mathcal F_sG)(r,a,\sqrt2b).
$$
The partial Fourier transform in \(s\) and the frequency dilation \(b\mapsto\sqrt2 b\) preserve the Schwartz space by \eqref{E107}, and multiplication by \(e^{ia\cdot b}\) is continuous by \eqref{E89}. Hence \(q_g\in\mathcal S(\mathbb R^{2n}\times\mathbb R^n)\). For fixed \(x\), the Weyl coefficient of the Schwartz vector \(q_g(x,\cdot)\) at \(y=(c,\zeta)\) is \(f_g(x,y)\). By \eqref{E25} it has the form
$$
f_g(x,y)=e^{ic\cdot\zeta}\int_{\mathbb R^n}
e^{-i\sqrt2\,\zeta\cdot R}q_g(x,R-\sqrt2c)\overline{v(R)}\,dR.
$$
Substitute \eqref{E117} for \(q_g(x,R-\sqrt2c)\) in this integral. By translation invariance in \(R\), the resulting integrand has \(L^1(dR\,ds)\) norm at most \(\|\kappa\|_1\|u\|_\infty\|v\|_\infty\). Fubini's theorem therefore allows us to integrate first in \(R\) and write
$$
\begin{aligned}
f_g(x,y)
=&e^{i(a\cdot b+c\cdot\zeta)}
\iint_{\mathbb R^n\times\mathbb R^n}\kappa(R-\sqrt2c,s)\\
&\quad\cdot u(s-\sqrt2a)\overline{v(R)}
e^{-i\sqrt2(\zeta\cdot R+b\cdot s)}\,dR\,ds.
\end{aligned}
$$ To leave the kernel in the unshifted variables \((r,s)\), put \(R=r+\sqrt2c\). This substitution has unit Jacobian and changes the phase factor to \(e^{i(a\cdot b+c\cdot\zeta)}e^{-2ic\cdot\zeta}=e^{i(a\cdot b-c\cdot\zeta)}\), so the coefficient becomes 
\begin{equation}
\begin{aligned}
f_g(x,y)
=&e^{i(a\cdot b-c\cdot\zeta)}
\iint\kappa(r,s)u(s-\sqrt2a)\overline{v(r+\sqrt2c)} \\
&\quad \cdot e^{-i\sqrt2(\zeta\cdot r+b\cdot s)}\,dr\,ds.
\end{aligned}
\label{eq:mixed-coefficient}
\end{equation}
Applying the pointwise tensor-product formula \eqref{E17} twice, we write the amplitude in \eqref{eq:mixed-coefficient} as the pullback
$$
H(a,c,r,s)=(\kappa\otimes u\otimes\overline v)(r,s,s-\sqrt2a,r+\sqrt2c).
$$
Since \(\kappa,u,\overline v\) are Schwartz and the linear change of variables on the right is invertible, \(H\in\mathcal S(\mathbb R^{4n})\) by \eqref{E107}.
With \((a,c)\) fixed, the partial Fourier transform in \((r,s)\) from \eqref{E107} lets us write \eqref{eq:mixed-coefficient} as
$$
f_g\bigl((a,b),(c,\zeta)\bigr)
=e^{i(a\cdot b-c\cdot\zeta)}
(\mathcal F_{r,s}H)(a,c,\sqrt2\,\zeta,\sqrt2\,b).
$$
Both the partial Fourier transform and the invertible linear substitution \((a,b,c,\zeta)\mapsto(a,c,\sqrt2\,\zeta,\sqrt2\,b)\) preserve the Schwartz topology by \eqref{E107}. Every derivative of the exponential factor is a polynomial times that factor, so multiplication by it is continuous by \eqref{E89}. Thus, \(f_g\in\mathcal S(\mathbb R^{4n})\). Its dependence on \(g,u,\overline v\) is jointly continuous because \(g\mapsto\kappa_g\) is continuous by Lemma \ref{L3}, the tensor product is jointly continuous by \eqref{E118}, and the subsequent operations are continuous linear maps.

We now evaluate the tensor pairing with \(f_g\) in \eqref{eq:mixed-pairing}. To act first in \(x\), we identify \(r\mapsto\langle U,q_g(\cdot,r)\rangle\) with \(S_gAu\). The joint Schwartz regularity proved from \eqref{E117} allows us to apply \eqref{E127} in \(x\), so this partial pairing is Schwartz in \(r\). For each \(r\), \eqref{E117} also reads \(q_g(x,r)=\langle W_xu,\overline{\kappa_g(r,\cdot)}\rangle_2\). Evaluating \(S_gAu\) by \eqref{E41} and taking \(\overline{\kappa_g(r,\cdot)}\) as the second vector in \eqref{eq:weyl-map}, we find
\begin{equation}
\begin{aligned}
(S_gAu)(r)
&=(Au)\bigl(\kappa_g(r,\cdot)\bigr)=\langle U(x),q_g(x,r)\rangle.
\end{aligned}
\label{E64}
\end{equation}

Fix \(y\in\mathbb R^{2n}\) and put \(v_y=W_y^*v=W_{-y}v\in\mathcal S(\mathbb R^n)\) by \eqref{E24} and \eqref{eq:weyl-relations}. Moving \(W_y\) to the second vector in the definition of \(f_g\) and using \eqref{E117}, we obtain
$$
f_g(x,y)=\langle S_gW_xu,v_y\rangle_2
=\int_{\mathbb R^n}q_g(x,r)\overline{v_y(r)}\,dr.
$$
We will evaluate \(U\otimes\overline{v_y}\) on the joint Schwartz test \(q_g\) in the two orders allowed by \eqref{E36}. To identify its partial pairing in \(r\) with the displayed integral in the Schwartz topology of \(x\), denote the integral over \(|r|\le R\) by \(F_{R,y}\), where \(R>0\), and estimate the tail in the seminorms \eqref{E51}. Each slice at \(r\in\mathbb R^n\) satisfies
$$
p_{M,k}^{(x)}\bigl(q_g(\cdot,r)\bigr)
=\max_{|\alpha|\le k}\sup_{x\in\mathbb R^{2n}}
\langle x\rangle^M|\partial_x^\alpha q_g(x,r)|
\le p_{M,k}^{(x,r)}(q_g).
$$
Indeed, the joint seminorm on the right includes the derivatives indexed by \((\alpha,0)\) in \((x,r)\), and \(\langle x\rangle\le\langle(x,r)\rangle\). Since \(v_y\in L^1(\mathbb R^n)\), this uniform bound permits differentiation in \(x\) under the tail integral. Taking the weighted suprema in \(x\) yields
$$
\begin{aligned}
p_{M,k}^{(x)}\bigl(f_g(\cdot,y)-F_{R,y}\bigr)
&=p_{M,k}^{(x)}\!\bigg(\int_{|r|>R}q_g(\cdot,r)\overline{v_y(r)}\,dr\bigg)\\
&\le\int_{|r|>R}p_{M,k}^{(x)}\bigl(q_g(\cdot,r)\bigr)|v_y(r)|\,dr\\
&\le p_{M,k}^{(x,r)}(q_g)\int_{|r|>R}|v_y(r)|\,dr\longrightarrow0
\end{aligned}
$$
as \(R\to\infty\). Then \(F_{R,y}\to f_g(\cdot,y)\) in \(\mathcal S(\mathbb R^{2n})\). For each \(R>0\), the truncated function \(\mathbf 1_{\{|r|\le R\}}\overline{v_y}\) is integrable. Applying \eqref{E36} to its tensor product with \(U\), with Schwartz test \(q_g\), and using \eqref{E44} for the truncated function, we obtain
$$
\langle U,F_{R,y}\rangle
=\int_{|r|\le R}\langle U,q_g(\cdot,r)\rangle\overline{v_y(r)}\,dr.
$$
The left side converges to \(\langle U,f_g(\cdot,y)\rangle\) by continuity of \(U\). On the right, \eqref{E64} identifies the integrand with \((S_gAu)(r)\overline{v_y(r)}\), which is integrable since both factors are Schwartz by Lemma \ref{L6} and the choice of \(v_y\). Passing to the limit and recalling \(v_y=W_y^*v\), we find
\begin{equation}
\begin{aligned}
\langle W_yS_gAu,v\rangle_2
&=\int_{\mathbb R^n}\langle U(x),q_g(x,r)\rangle\overline{v_y(r)}\,dr\\
&=\langle U(x),\int_{\mathbb R^n}q_g(x,r)\overline{v_y(r)}\,dr\rangle
=\langle U(x),f_g(x,y)\rangle.
\end{aligned}
\label{E65}
\end{equation}

Lemma \ref{L6} places \(S_gAu\) in \(\mathcal S(\mathbb R^n)\), so its Weyl coefficient in \eqref{E65} is a Schwartz test in \(y\) by \eqref{E88}. Pairing \eqref{E65} with \(V\) and applying \eqref{eq:weyl-map} to \(B=W(V)\), we obtain
$$
\begin{aligned}
\langle BS_gAu,\overline v\rangle
&=\Big\langle V(y),\big\langle U(x),f_g(x,y)\big\rangle\Big\rangle=\langle U(x)\otimes V(y),f_g(x,y)\rangle,
\end{aligned}
$$
where \eqref{E36} applies to the Schwartz test \(f_g\) in the last equality. This proves \eqref{eq:mixed-pairing}.

Now suppose \(U,V\in\mathcal C_{e,\eta}\) and \(g\in C_c^\infty(\mathbb R^{2n})\). For the compact approximants, the Weyl composition formula \eqref{eq:compact-composition} will identify the two sides of \eqref{eq:sandwich}. We use \eqref{eq:mixed-pairing} to pass to the limit on the operator side. Choose \(L>0\) so that \(|z|\le L\) on \(\operatorname{supp}g\). By \eqref{E109}, \(g\in\mathcal C_{e,\eta}\) with \(M_g=(\eta+|e|)L\). Take \(\chi\in C_c^\infty(\mathbb R^{2n})\) equal to one on the closed unit ball and set \(\chi_N(z)=\chi(z/N)\), \(U_N=\chi_NU\), and \(V_N=\chi_NV\) for \(N=1,2,\ldots\). Fix \(u,v\in\mathcal S(\mathbb R^n)\). By \eqref{E88}, the coefficient \(c_{u,v}(z)=\langle W_z u,v\rangle_2\) belongs to \(\mathcal S(\mathbb R^{2n})\).

For each \(N\), the factors \(V_N,g,U_N\) are compactly supported, and \eqref{eq:product-support} places \(\operatorname{supp}(V_N\star g)\) in the compact sum \(\operatorname{supp}V_N+\operatorname{supp}g\). For compact factors, Proposition \ref{prop:tempered-product} identifies the product in \(\mathcal C_{0,1}\) used in \eqref{eq:compact-composition} and the product in \(\mathcal C_{e,\eta}\) as the unique tempered extension of the same local pairing \eqref{E23}. The Weyl maps of \(V_N,g,U_N\), and \(V_N\star g\) preserve \(\mathcal S(\mathbb R^n)\) by \eqref{E6}. Applying \eqref{eq:compact-composition} first to \(V_N,g\) and then to \(V_N\star g,U_N\), with \(W(g)=S_g\) on \(\mathcal S(\mathbb R^n)\) by Lemma \ref{L6}, we obtain
$$
W(V_N\star g\star U_N)u=W(V_N)S_gW(U_N)u.
$$
The Weyl pairing \eqref{eq:weyl-map} for the left side and the mixed pairing \eqref{eq:mixed-pairing} for the right side now give
\begin{equation}
\begin{aligned}
\langle V_N\star g\star U_N,c_{u,v}\rangle
&=\langle W(V_N)S_gW(U_N)u,\overline v\rangle=\langle U_N\otimes V_N,f_g\rangle\\
&=\langle U\otimes V,(\chi_N\otimes\chi_N)f_g\rangle,
\end{aligned}
\label{E66}
\end{equation}
where the last equality moves the cutoffs from \(U_N,V_N\) onto \(f_g\) by \eqref{E81}.

On the left of \eqref{E66} the test \(c_{u,v}\) is fixed. Since \(\chi_Ng=g\) for \(N\ge L\), weak convergence in \eqref{eq:finite-cutoff} for \((V,g,U)\), tested on this coefficient, yields
\begin{equation*}
\langle V_N\star g\star U_N,c_{u,v}\rangle
\longrightarrow \langle V\star g\star U,c_{u,v}\rangle
\qquad(N\to\infty).
\tag{\getrefnumber{E66}a}\label{E119}
\end{equation*}
On the right of \eqref{E66}, \(U\otimes V\) is fixed, so we estimate the change in its test function. With \(\lambda=4n\) and \(\vartheta=\chi\otimes\chi\), we have \(\vartheta=1\) on the unit ball. Applying \eqref{E5} to \(f_g\in\mathcal S(\mathbb R^{4n})\), we obtain
$$
p_{M,k}^{(x,y)}\bigl((\chi_N\otimes\chi_N-1)f_g\bigr)
\le \frac{2^k C_{k,\vartheta}}{N}\,p_{M+1,k}^{(x,y)}(f_g)
\longrightarrow 0
\qquad(N\to\infty)
$$
for all \(M,k\in\mathbb N_0\). Thus \((\chi_N\otimes\chi_N)f_g\to f_g\) in \(\mathcal S(\mathbb R^{4n})\) by \eqref{E51}. By the tensor product theorem recalled after \eqref{E36}, \(U\otimes V\) acts continuously on this space. Hence
\begin{equation*}
\langle U\otimes V,(\chi_N\otimes\chi_N)f_g\rangle
\longrightarrow \langle U\otimes V,f_g\rangle
\qquad(N\to\infty).
\tag{\getrefnumber{E66}b}\label{E120}
\end{equation*}
By \eqref{E66}, the scalar pairings in \eqref{E119} and \eqref{E120} agree for every \(N\). Equality of their limits gives the middle equality in
$$
\begin{aligned}
\langle W(V\star g\star U)u,\overline v\rangle
&=\langle V\star g\star U,c_{u,v}\rangle=\langle U\otimes V,f_g\rangle
=\langle BS_gAu,\overline v\rangle.
\end{aligned}
$$
The first equality is \eqref{eq:weyl-map} with \(V\star g\star U\) acting on \(c_{u,v}(z)=\langle W_z u,v\rangle_2\), and the last is \eqref{eq:mixed-pairing} for \(U,V\). Since \(u,v\in\mathcal S(\mathbb R^n)\) are arbitrary, \eqref{eq:sandwich} holds as an identity of maps \(\mathcal S(\mathbb R^n)\to\mathcal S'(\mathbb R^n)\). The proof is complete.
\end{proof}

For \(U,V\in\mathcal C_{e,\eta}\) and \(g\in C_c^\infty(\mathbb R^{2n})\), \eqref{eq:sandwich} reduces nonvanishing of \(V\star g\star U\) to that of \(W(V)S_gW(U)\). The operator construction works for any two nonzero tempered distributions \(U,V\). A rank-one map can send a nonzero value of \(W(U)\) to a Schwartz function outside the kernel of \(W(V)\). To replace that map by \(S_g\) with compactly supported \(g\), it suffices to preserve nonvanishing of the pairing \(\langle W(V)S_gW(U)u,\overline v\rangle\) for one pair \(u,v\in\mathcal S(\mathbb R^n)\).

Exact realization of the rank-one map with such a \(g\) is precluded by \cite[Theorem 2.2]{NarayananRatnakumar2010}. Indeed, for \(J(s,t)=-\sqrt2\,t-i\sqrt2\,s\), the representation in \cite[Section 1]{NarayananRatnakumar2010} satisfies \(\pi_1(J(s,t),0)=W_{(s,t)}\) by \eqref{E24}. Since \(|\det J|=2^n\), \(W(g)\) is the Weyl transform in that normalization of \(2^{-n}g\circ J^{-1}\), a nonzero compactly supported \(L^1\) function whenever \(g\in C_c^\infty\) is nonzero. Its Weyl transform has infinite rank.

We approximate the Schwartz spreading function of the rank-one map by compactly supported smooth functions. To justify this step, suppose \(g_j\to g\) in \(\mathcal S(\mathbb R^{2n})\). By \eqref{E124}, \(\kappa_{g_j}-\kappa_g\to0\) in the same space. Applying \eqref{eq:fixed-distribution-smoothing} to the difference gives
\begin{equation}
S_{g_j}w\longrightarrow S_gw
\quad\text{in }\mathcal S(\mathbb R^n)
\label{E87}
\end{equation}
as \(j\to\infty\) for each fixed \(w\in\mathcal S'(\mathbb R^n)\).

\begin{lemma}\label{lem:nonzero-sandwich}

If \(U,V\in\mathcal S'(\mathbb R^{2n})\) are nonzero, there is \(g\in C_c^\infty(\mathbb R^{2n})\) such that
$$
W(V)S_gW(U)\ne0
$$
as a map \(\mathcal S(\mathbb R^n)\to\mathcal S'(\mathbb R^n)\).

\end{lemma}

\begin{proof}

Let \(A=W(U)\) and \(B=W(V)\). Since both maps are nonzero by \eqref{E40}, choose \(u,p\in\mathcal S(\mathbb R^n)\) with \(\alpha=Au\ne0\) and \(Bp\ne0\). These are tempered distributions, so there exist \(q,v\in\mathcal S(\mathbb R^n)\) such that
\begin{equation*}
\alpha(q)\ne0\quad\text{and}\quad \langle Bp,\overline v\rangle\ne0.
\tag{\getrefnumber{E87}a}\label{E123}
\end{equation*}
To send \(\alpha=Au\) to a vector on which \(B\) is nonzero, define \(S:\mathcal S'(\mathbb R^n)\to\mathcal S(\mathbb R^n)\) by \(Sf=f(q)p\). Then \(S\alpha=\alpha(q)p\), so \(BS\alpha=\alpha(q)Bp\ne0\).

To construct \(g_0\) with \(S_{g_0}=S\), we first compute the kernel of the restriction of \(S\) to \(\mathcal S(\mathbb R^n)\). For \(h\in\mathcal S(\mathbb R^n)\), the embedding \eqref{E18} and the integral pairing \eqref{E44} give
$$
(Sh)(r)=\langle h,q\rangle p(r)=\int_{\mathbb R^n}p(r)q(s)h(s)\,ds.
$$
By \eqref{E17} the kernel is \(p\otimes q\), which is Schwartz by \eqref{E118}. Applying the inverse map in \eqref{E124} to \(p\otimes q\) gives the unique \(g_0\in\mathcal S(\mathbb R^{2n})\) with \(K_{g_0}=p\otimes q\). To act on \(\alpha=Au\), we use the extension \(S_{g_0}\) of \(W(g_0)\) from Lemma \ref{L6}. For fixed \(r\), its kernel slice \(K_{g_0}(r,\cdot)=p(r)q\) belongs to \(\mathcal S(\mathbb R^n)\). Applying \eqref{E41} to any \(w\in\mathcal S'(\mathbb R^n)\) and using linearity of \(w\), we obtain
\begin{equation}
(S_{g_0}w)(r)=w\bigl(p(r)q(\cdot)\bigr)
=w(q)p(r)=(Sw)(r).
\label{E67}
\end{equation}

By \eqref{E67} and \eqref{E123}, the coefficient to be preserved under compact approximation is \(\langle BS_{g_0}Au,\overline v\rangle=\alpha(q)\langle Bp,\overline v\rangle\ne0\). Choose \(\chi\in C_c^\infty(\mathbb R^{2n})\) equal to one near zero and put \(g_j=\chi(\cdot/j)g_0\in C_c^\infty(\mathbb R^{2n})\) for positive integers \(j\). With \(\lambda=2n\), \(\vartheta=\chi\), and \(\psi=g_0\), estimate \eqref{E5} implies \(g_j\to g_0\) in \(\mathcal S(\mathbb R^{2n})\) as \(j\to\infty\). Applying \eqref{E87} to the fixed distribution \(\alpha\) and evaluating the limit by \eqref{E67}, we find that
$$
S_{g_j}\alpha\longrightarrow S_{g_0}\alpha=\alpha(q)p
\quad\text{in }\mathcal S(\mathbb R^n).
$$
Set \(d_j=S_{g_j}\alpha-\alpha(q)p\). Since \(d_j\to0\) in \(\mathcal S(\mathbb R^n)\), \eqref{E88} with \(v\) fixed implies that the coefficient \(a\mapsto\langle W_a d_j,v\rangle_2\) tends to zero in \(\mathcal S(\mathbb R^{2n})\). With \(\alpha=Au\) and \(B=W(V)\), formula \eqref{eq:weyl-map} and the continuity of \(V\) on this space yield
$$
\begin{aligned}
\langle BS_{g_j}Au,\overline v\rangle-\alpha(q)\langle Bp,\overline v\rangle
&=\langle Bd_j,\overline v\rangle=\big\langle V(a),\langle W_a d_j,v\rangle_2\big\rangle
\longrightarrow0
\end{aligned}
$$ as $j\to\infty.$
By \eqref{E123}, \(\alpha(q)\langle Bp,\overline v\rangle\ne0\), so \(\langle BS_{g_j}Au,\overline v\rangle\ne0\) for all sufficiently large \(j\). For any such \(j\), the compactly supported function \(g=g_j\) satisfies \(W(V)S_gW(U)\ne0\). This completes the proof.
\end{proof}

\begin{proof}[Proof of Theorem \ref{thm:distribution}]
For the algebra structure, fix \(U_1,\ldots,U_m\in\mathcal C_{e,\eta}\) and take
\(X=\sum_{i=1}^m\lambda_iU_i\) and \(Y=\sum_{j=1}^m\mu_jU_j\) with complex coefficients.
The finite union \(F=\bigcup_{j=1}^m\operatorname{supp}U_j\) is closed. For \(\varphi\in C_c^\infty(\mathbb R^d\setminus F)\), each \(\langle U_j,\varphi\rangle=0\), so \(\langle X,\varphi\rangle=\langle Y,\varphi\rangle=0\).
By Definition \ref{def:distribution-support} and \eqref{eq:cone} with \(M=\max_{1\le j\le m}M_{U_j}\), it follows that
$$
\operatorname{supp}X\cup\operatorname{supp}Y
\subset\bigcup_{j=1}^m\operatorname{supp}U_j
\subset\{a\in\mathbb R^d:e\cdot a\ge\eta|a|-M\}.
$$
In particular, \(X,Y\in\mathcal C_{e,\eta}\), so \(\mathcal C_{e,\eta}\) is a complex vector space.
Using the same \(M\) in \eqref{E101} gives one cutoff for every element of this span.
For each \(\varphi\in\mathcal S(\mathbb R^d)\), the test \(j_\varphi\) in \eqref{eq:schwartz-pullback} is therefore common to all pairs.
Expanding the tensor product in \eqref{eq:distribution-product} by \eqref{E36}, we obtain
$$
\langle X\star Y,\varphi\rangle
=\sum_{i,j=1}^m\lambda_i\mu_j\langle U_i\otimes U_j,j_\varphi\rangle
=\sum_{i,j=1}^m\lambda_i\mu_j\langle U_i\star U_j,\varphi\rangle.
$$
As \(\varphi\) is arbitrary, \(X\star Y=\sum_{i,j=1}^m\lambda_i\mu_j(U_i\star U_j)\), proving bilinearity.
For \(U,V\in\mathcal C_{e,\eta}\), \eqref{eq:product-support} places \(U\star V\) in the same class with support constant \(M_U+M_V\).
Proposition \ref{prop:tempered-product} also gives \((U\star V)\star T=U\star(V\star T)\) for every \(T\in\mathcal C_{e,\eta}\).

The point mass \(\delta_0\) belongs to \(\mathcal C_{e,\eta}\) with \(M_{\delta_0}=0\).
For \(\varphi\in C_c^\infty(\mathbb R^d)\), the cutoff in \eqref{E23} may be taken equal to one near both \(\{0\}\times\operatorname{supp}\varphi\) and \(\operatorname{supp}\varphi\times\{0\}\). Since \(\sigma(0,a)=\sigma(a,0)=0\), evaluating the point mass in either variable yields
$$
\langle\delta_0\star U,\varphi\rangle=\langle U,\varphi\rangle
=\langle U\star\delta_0,\varphi\rangle
\qquad(\varphi\in C_c^\infty(\mathbb R^d)).
$$
By the density in \eqref{E18}, the tempered extensions agree on \(\mathcal S(\mathbb R^d)\).
Thus, \(\delta_0\star U=U=U\star\delta_0\), completing the proof that \(\mathcal C_{e,\eta}\) is a unital associative complex algebra.

To exclude zero divisors, suppose \(U,V\in\mathcal C_{e,\eta}\setminus\{0\}\) and \(U\star V=0\).
We seek a nonzero smooth function of the form \(h=V\star g\star U\), so that \(U\star V\) occurs in the middle of its square.
By Lemma \ref{lem:nonzero-sandwich} there exists \(g\in C_c^\infty(\mathbb R^d)\) with \(W(V)S_gW(U)\ne0\).
Estimate \eqref{E109} places \(g\) in \(\mathcal C_{e,\eta}\), and Lemma \ref{lem:smooth-regularization} identifies \(h\) with a smooth function in this class.
For this choice of \(g\), the composition formula \eqref{eq:sandwich} yields
\begin{equation}
W(h)=W(V)S_gW(U)\ne0.
\label{E68}
\end{equation}
In particular, \(h\not\equiv0\). By associativity in Proposition \ref{prop:tempered-product}, the adjacent factors \(U,V\) in
\((V\star g\star U)\star(V\star g\star U)\) can be grouped together.
The assumption \(U\star V=0\) and bilinearity then imply in \(\mathcal S'(\mathbb R^d)\) that
\begin{equation}
h\star h
=V\star g\star(U\star V)\star g\star U=0.
\label{E48}
\end{equation}

To apply Proposition \ref{prop:smooth-square}, we identify the distributional square in \eqref{E48} with the locally defined smooth square.
For \(\ell(z)=-e\cdot z\), the bound \eqref{eq:smoothing-cone} reads
\(\ell(z)\le M_h-\eta|z|\) on \(\operatorname{supp}h\), so \eqref{E80} verifies condition \((Q)\) for this support.
Lemma \ref{L2} therefore applies with \(f=g=h\). The integral in \eqref{eq:twisted-product} converges absolutely for each \(z\in\mathbb R^d\) and defines the smooth function
$$
k(z)=\int_{\mathbb R^d}h(x)h(z-x)e^{i\sigma(x,z-x)}\,dx.
$$
Let \(T\in\mathcal D'(\mathbb R^d)\) be the pushforward of \(e^{i\sigma(x,y)}(h\otimes h)\) under addition \((x,y)\mapsto x+y\). Testing the two representations of \(T\) on \(\varphi\in C_c^\infty(\mathbb R^d)\) relates the function \(k\) to the vanishing distributional square in \eqref{E48} through the identities
$$
\int_{\mathbb R^d}k(z)\varphi(z)\,dz
\overset{\text{\eqref{E86}}}{=}\langle T,\varphi\rangle
\overset{\text{Prop. \ref{prop:tempered-product}}}{=}\langle h\star h,\varphi\rangle
\overset{\text{\eqref{E48}}}{=}0.
$$
Thus, \(k\) represents the zero distribution, and its continuity and \eqref{E37} imply \(k\equiv0\).
Proposition \ref{prop:smooth-square} now forces \(h\equiv0\), contradicting \eqref{E68}.
Hence \(U\star V=0\) in \(\mathcal C_{e,\eta}\) implies \(U=0\) or \(V=0\), which completes the proof.
\end{proof}

\begin{remark}[An unbounded Weyl map]
The cone class in Theorem \ref{thm:distribution} contains distributions whose Weyl maps have no bounded \(L^2\) extension. Fix \(j\in\{1,\ldots,n\}\) and take \(U=\partial_{s_j}\delta_0\), so \(U(\varphi)=-\partial_{s_j}\varphi(0)\) for \(\varphi\in\mathcal S(\mathbb R^{2n})\). In \eqref{eq:weyl-map}, pairing with \(\partial_{s_j}\delta_0\) differentiates the Weyl coefficient at the origin with a minus sign. For \(u,v\in\mathcal S(\mathbb R^n)\), estimate \eqref{E26} with \(N=k=0\) and \(m=1\) bounds the differentiated integrand in \eqref{E25} near \(a=0\) by a constant times \(|v(x)|\). Differentiating under this integral and using \eqref{E84} at the origin, we find
$$
\langle W(U)u,\overline v\rangle
=-\left.\partial_{s_j}\langle W_{(s,t)}u,v\rangle_2\right|_{(0,0)}
=\langle\sqrt2\,\partial_{x_j}u,v\rangle_2.
$$
Since \(v\) is arbitrary, the regular distribution pairing \eqref{E44} yields
$$
W(U)u=-\left.\partial_{s_j}W_{(s,t)}u\right|_{(0,0)}
=\sqrt2\,\partial_{x_j}u.
$$
Since \(\operatorname{supp}U=\{0\}\), the cone bound \eqref{eq:cone} holds with \(M_U=0\) for every admissible \(e,\eta\). To test whether \(W(U)\) has a bounded \(L^2\) extension, we apply the derivative to functions of increasing frequency. Fix \(\phi\in C_c^\infty(\mathbb R^n)\) with \(\phi\not\equiv0\), and set \(u_k(x)=e^{ikx_j}\phi(x)\) for \(x\in\mathbb R^n\) and positive integers \(k\). Since \(\partial_{x_j}u_k=e^{ikx_j}(ik\phi+\partial_{x_j}\phi)\), the modulus of the exponential and the reverse triangle inequality give
$$
\|u_k\|_2=\|\phi\|_2,
\qquad\|\partial_{x_j}u_k\|_2
\ge k\|\phi\|_2-\|\partial_{x_j}\phi\|_2\longrightarrow\infty
$$
as \(k\to\infty\). Thus, \(\|W(U)u_k\|_2/\|u_k\|_2=\sqrt2\,\|\partial_{x_j}u_k\|_2/\|\phi\|_2\to\infty\), which rules out a bounded extension.
\end{remark}

\begin{example}[Proper addition and a vanishing twisted product]\label{Eg1}
In dimension \(n=1\), write \(z=(s,t)\in\mathbb R^2\) and set
$$
U=-\frac12(\partial_s^2+\partial_t^2)\delta_0-\delta_0,
\qquad
V(s,t)=e^{-(s^2+t^2)/2}.
$$
Then \(U\) and \(V\) are nonzero tempered distributions, addition is proper on \(\operatorname{supp}U\times\operatorname{supp}V\), and \(U\star V=0\).
\end{example}

\begin{proof}
The value of \(\langle U,\psi\rangle\) depends only on \(\psi(0)\), \(\partial_s^2\psi(0)\), and \(\partial_t^2\psi(0)\), so \eqref{E51} gives \(U\in\mathcal S'(\mathbb R^2)\), with \(\operatorname{supp}U\subset\{0\}\). If \(\psi\in C_c^\infty(\mathbb R^2)\) equals one near zero, then \(\langle U,\psi\rangle=-1\), so \(\operatorname{supp}U=\{0\}\). The strictly positive Schwartz function \(V\) defines a tempered distribution by \eqref{E18}, with \(\operatorname{supp}V=\mathbb R^2\) by \eqref{E37}. On \(\{0\}\times\mathbb R^2\), the inverse image under addition of any compact set \(L\subset\mathbb R^2\) is \(\{0\}\times L\). Thus \eqref{E34} holds and \eqref{E23} defines \(U\star V\). The full support of \(V\) precludes the common cone condition in Theorem \ref{thm:distribution}. Indeed, for any \(e\in\mathbb R^2\) there is a unit vector \(\omega\) with \(e\cdot\omega\le0\), so \(e\cdot(R\omega)\le0<\eta R-M_V\) for all sufficiently large \(R\).

To compute the contribution of \(\Delta\delta_0\), where \(\Delta=\partial_s^2+\partial_t^2\), fix \(\varphi\in C_c^\infty(\mathbb R^2)\). Let \(J_\varphi\) denote the test function obtained in \eqref{E23} by pairing first with \(V\) in the second variable \(y\), as in \eqref{E36}. Writing \(w=(a,b)\) and substituting \(z=w+y=(s,t)\) in the integral given by \eqref{E44}, we have \(\sigma(w,z-w)=at-bs\) by \eqref{E55}. Since \(\chi_\varphi=1\) near \(\{0\}\times\operatorname{supp}\varphi\) and \(\operatorname{supp}\varphi\) is compact, \(\chi_\varphi(w,z-w)=1\) for all \(z\in\operatorname{supp}\varphi\) whenever \(|w|\) is sufficiently small. On this neighborhood of zero, \(J_\varphi\) therefore has the form
$$
J_\varphi(w)=\int_{\mathbb R^2}e^{i(at-bs)}V(s-a,t-b)\varphi(s,t)\,ds\,dt.
$$
By the definition of distributional derivatives, \(\Delta\delta_0\) acts by the Laplacian in \(w\) at zero. The smooth integrand has fixed compact support in \(z\), so differentiation under the integral gives
$$
\langle(\Delta\delta_0)\star V,\varphi\rangle
=\langle\Delta\delta_0,J_\varphi\rangle
=(\Delta_w J_\varphi)(0)
=\int_{\mathbb R^2}F(z)\varphi(z)\,dz,
$$
where the product rule at fixed \(z=(s,t)\) yields
$$
\begin{aligned}
F(s,t)
&=\left.(\partial_a^2+\partial_b^2)
\bigl[e^{i(at-bs)}V(s-a,t-b)\bigr]\right|_{a=b=0}\\
&=\bigl((\partial_s-it)^2+(\partial_t+is)^2\bigr)V\\
&=\bigl(\partial_s^2+\partial_t^2-s^2-t^2
+2i(s\partial_t-t\partial_s)\bigr)V.
\end{aligned}
$$
For the Gaussian \(V\), the identities \(\partial_sV=-sV\) and \(\partial_tV=-tV\) give \(s\partial_tV-t\partial_sV=(-st+ts)V=0\) and \(\Delta V=(s^2+t^2-2)V\). Thus \(F=-2V\), so \((\Delta\delta_0)\star V=-2V\).

Evaluating \(J_\varphi\) itself at zero gives \(J_\varphi(0)=\langle V,\varphi\rangle\), and hence \(\delta_0\star V=V\) by \eqref{E23}. Combining the two contributions to \(U\star V\) now yields
\begin{equation*}
U\star V=-\tfrac12\bigl((\Delta\delta_0)\star V\bigr)-\delta_0\star V
=-\tfrac12(-2V)-V=0.\qedhere
\end{equation*}
\end{proof}

\section{Bounded Weyl operators and the parabolic algebra}\label{sec:5}

\subsection{Bounded Weyl operators}

To apply Theorem \ref{thm:distribution} to products of bounded operators, we first identify the spreading distribution of each bounded operator \(A\) on \(L^2(\mathbb R^n)\). Let \(i:\mathcal S(\mathbb R^n)\hookrightarrow L^2(\mathbb R^n)\) denote the continuous inclusion. For \(f\in L^2(\mathbb R^n)\), write \(Jf\) for the tempered distribution in \eqref{E44}, so that
$$
\langle Jf,q\rangle=\int_{\mathbb R^n}f(x)q(x)\,dx
=\langle f,\overline q\rangle_2,\qquad q\in\mathcal S(\mathbb R^n).
$$
On each bounded set \(D\subset\mathcal S(\mathbb R^n)\), Cauchy--Schwarz gives \(\nu_D(Jf)\le\|f\|_2\sup_{q\in D}\|q\|_2\) in the notation of \eqref{E128}. The supremum is finite by continuity of \(i\), so \(J:L^2(\mathbb R^n)\hookrightarrow\mathcal S'(\mathbb R^n)\) is continuous. If \(Jf=0\), then \(\langle f,v\rangle_2=0\) for every \(v\in\mathcal S(\mathbb R^n)\), and density in \(L^2(\mathbb R^n)\) forces \(f=0\). With this identification, \(T=J\circ A\circ i\) is a continuous linear map from \(\mathcal S(\mathbb R^n)\) to \(\mathcal S'(\mathbb R^n)\).

Applying the tempered Schwartz kernel theorem \cite[Corollary to Theorem 51.6, (51.7)]{Treves1967} to \(T\) with \(m=n\), we obtain its unique kernel \(K_A\in\mathcal S'(\mathbb R^{2n})\). Evaluating \(Tu=J(Au)\) in the kernel identity gives
\begin{equation}
\langle K_A(x,y),u(y)q(x)\rangle
=\langle Tu,q\rangle=\langle Au,\overline q\rangle_2,
\qquad u,q\in\mathcal S(\mathbb R^n).
\label{E78}
\end{equation}
Lemma \ref{L3} determines the unique \(U_A\in\mathcal S'(\mathbb R^{2n})\) with \(K_{U_A}=K_A\). Taking \(q=\overline v\) in \eqref{E78} and using \eqref{E116}, we obtain \(\langle W(U_A)u,\overline v\rangle=\langle Tu,\overline v\rangle\) for all \(u,v\in\mathcal S(\mathbb R^n)\). Then \(W(U_A)=T\), which identifies \(W(U_A)\) with \(A|_{\mathcal S(\mathbb R^n)}\) under \(J\). By density of \(\mathcal S(\mathbb R^n)\) in \(L^2(\mathbb R^n)\), \(A\) is the unique bounded extension to \(L^2(\mathbb R^n)\).

The spreading distribution \(U_A\) of a bounded operator may have noncompact support. To approximate \(A\) by operators to which Lemma \ref{L4} applies, we express the Weyl maps of compactly supported cutoffs of \(U_A\) as averages of Weyl conjugates of \(A\). These averages are instances of the function--operator convolution of quantum harmonic analysis \cite{LuefSkrettingland2018, Werner1984}.

\begin{lemma}\label{lem:bounded-cutoffs}

Let \(U\in\mathcal S'(\mathbb R^d)\), and suppose that \(W(U)\) extends to a bounded operator \(A\) on \(L^2(\mathbb R^n)\). For \(\chi\in C_c^\infty(\mathbb R^d)\) equal to one near zero, set
$$
A_N=W(\chi(\cdot/N)U)
$$
for \(N=1,2,\ldots\). These maps have bounded extensions to \(L^2(\mathbb R^n)\), denoted by the same symbols, with \(\sup_N\|A_N\|<\infty\), and satisfy
$$
A_Nu\longrightarrow Au
$$
in \(L^2(\mathbb R^n)\) as \(N\to\infty\) for every \(u\in L^2(\mathbb R^n)\).

\end{lemma}

\begin{proof}

Write \(d=2n\). To realize the cutoff at \(a=(s,\tau)\in\mathbb R^d\), we average conjugates of \(A\) over \(t=(t_1,t_2)\in\mathbb R^d\). With \(L_\sigma(s,\tau)=(2\tau,-2s)\), formula \eqref{E55} reads
$$
2\sigma(t,a)=2(t_1\cdot\tau-t_2\cdot s)=t\cdot L_\sigma a.
$$
To recover \(\chi(a)\) by inversion at \(L_\sigma a\), take \(F=\chi\circ L_\sigma^{-1}\), so \(F(L_\sigma a)=\chi(a)\). The linear isomorphism \(L_\sigma\) preserves smoothness and compact support under pullback, hence \(F\in C_c^\infty(\mathbb R^d)\subset\mathcal S(\mathbb R^d)\). Since the Fourier transform preserves \(\mathcal S(\mathbb R^d)\), the function \(k=(2\pi)^{-d}\widehat F\) is Schwartz. Fourier inversion at \(L_\sigma a\) in the convention \eqref{E7} reads
$$
\chi(a)=F(L_\sigma a)
=(2\pi)^{-d}\int_{\mathbb R^d}\widehat F(t)e^{it\cdot L_\sigma a}\,dt
=\int_{\mathbb R^d}k(t)e^{2i\sigma(t,a)}\,dt.
$$
At \(a=0\) this identity gives \(\int k=\chi(0)=1\). For \(k_N(t)=N^d k(Nt)\), the substitution \(r=Nt\) gives \(\int k_N=\int k\) and \(\|k_N\|_1=\|k\|_1\). Applying the same substitution to the representation of \(\chi(a/N)\) and using \(\sigma(Nt,a/N)=\sigma(t,a)\), we obtain
\begin{equation}
\chi(a/N)=\int_{\mathbb R^d}k_N(t)e^{2i\sigma(t,a)}\,dt.
\label{E27}
\end{equation}
To realize the multiplier in \eqref{E27} by conjugation, write \(W_t^*=W_{-t}\) and apply \eqref{eq:weyl-relations} twice. The Weyl parameter is \(t+a-t=a\), and the combined exponent is \(i[\sigma(t,a)+\sigma(t+a,-t)]=2i\sigma(t,a)\) by \eqref{E55}. Then for \(t,a\in\mathbb R^d\) we have
\begin{equation*}
W_tW_aW_t^*=e^{2i\sigma(t,a)}W_a.
\tag{\getrefnumber{E27}a}\label{E111}
\end{equation*}

In view of \eqref{E27} and \eqref{E111}, for \(u\in L^2(\mathbb R^n)\) we set
\begin{equation}
\widetilde A_Nu
=\int_{\mathbb R^d}k_N(t)W_tAW_t^*u\,dt. \label{eq:conjugation-average}
\end{equation}
To justify this vector integral, fix \(t_0\in\mathbb R^d\). Unitarity and strong continuity of the translations and modulations in \eqref{E24}, together with \(W_t^*=W_{-t}\), imply
\begin{equation*}
\begin{aligned}
\|W_tAW_t^*u-W_{t_0}AW_{t_0}^*u\|_2
&\le \|A\|\,\|(W_t^*-W_{t_0}^*)u\|_2\\
&\quad+\|(W_t-W_{t_0})AW_{t_0}^*u\|_2\longrightarrow0
\end{aligned}
\tag{\getrefnumber{eq:conjugation-average}a}\label{E112}
\end{equation*}
as \(t\to t_0\). Moreover, unitarity bounds the integral of the norm by
\begin{equation}
\int_{\mathbb R^d}|k_N(t)|\,\|W_tAW_t^*u\|_2\,dt
\le\|k_N\|_1\|A\|\|u\|_2
=\|k\|_1\|A\|\|u\|_2.
\label{E69}
\end{equation}
For each \(u\), the integrand in \eqref{eq:conjugation-average} is \(L^2\)-continuous by \eqref{E112} and continuity of \(k_N\). The norm bound \eqref{E69} ensures Bochner integrability and yields \(\|\widetilde A_Nu\|_2\le\|k\|_1\|A\|\|u\|_2\), uniformly in \(N\).

To transfer the uniform operator bound to \(A_N\), we compare the pairings of \(\widetilde A_N\) and \(A_N\) on \(\mathcal S(\mathbb R^n)\). Fix \(u,v\in\mathcal S(\mathbb R^n)\) and set \(c(a)=\langle W_a u,v\rangle_2\). The adjoint \(W_t^*=W_{-t}\) in \eqref{eq:weyl-relations} preserves \(\mathcal S\) by \eqref{E24}. Moving \(W_t\) to the second vector gives \(\langle AW_t^*u,W_t^*v\rangle_2\). Since \(A\) extends \(W(U)\), \eqref{eq:weyl-map} applies to these two Schwartz vectors. Moving the adjoint back inside the Weyl coefficient and applying \eqref{E111}, we obtain
\begin{equation}
\begin{aligned}
\langle W_tAW_t^*u,v\rangle_2
&=\big\langle U(a),\langle W_tW_aW_t^*u,v\rangle_2\big\rangle\\
&=\big\langle U(a),e^{2i\sigma(t,a)}\langle W_a u,v\rangle_2\big\rangle.
\end{aligned}
\label{E28}
\end{equation}
Applying the bounded functional \(w\mapsto\langle w,v\rangle_2\) to the Bochner integral \eqref{eq:conjugation-average} and substituting \eqref{E28} yields the scalar integral below. Its integrand pairs \(U\) with \(k_N(t)e^{2i\sigma(t,a)}c(a)\) in \(a\). The coefficient \(c\) is Schwartz by \eqref{E88}, so \(k_N\otimes c\) is Schwartz by \eqref{E118}. Multiplication by \(e^{2i\sigma(t,a)}\) preserves this space by \eqref{E89}. We may therefore evaluate \(1(t)\otimes U(a)\) on this joint test in either order by \eqref{E36}, obtaining
\begin{equation*}
\begin{aligned}
&\langle\widetilde A_Nu,v\rangle_2
=\int_{\mathbb R^d}k_N(t)\big\langle U(a),e^{2i\sigma(t,a)}c(a)\big\rangle\,dt=\langle U,m_N c\rangle,\\
&\quad m_N(a)=\int_{\mathbb R^d}k_N(t)e^{2i\sigma(t,a)}\,dt.
\end{aligned}
\tag{\getrefnumber{E28}a}\label{E126}
\end{equation*}
By \eqref{E27}, \(m_N(a)=\chi(a/N)\). Substituting in \eqref{E126}, moving the cutoff onto \(U\) by \eqref{E35}, and applying \eqref{eq:weyl-map}, we find
$$
\langle\widetilde A_Nu,v\rangle_2
=\langle U,\chi(\cdot/N)c\rangle
=\langle\chi(\cdot/N)U,c\rangle=\langle A_Nu,\overline v\rangle.
$$
As \(v\) ranges over \(\mathcal S(\mathbb R^n)\), the last identity says \(J\widetilde A_Nu=A_Nu\) for every Schwartz \(u\), where \(J:L^2\hookrightarrow\mathcal S'\). Density of \(\mathcal S(\mathbb R^n)\) in \(L^2(\mathbb R^n)\) makes \(\widetilde A_N\) the unique bounded extension. We denote it by \(A_N\) in the strong convergence argument.

For strong convergence, fix \(u\in L^2(\mathbb R^n)\). Rescaling the integration variable in \eqref{eq:conjugation-average} and using \(\int k=1\) to subtract \(Au\) inside the integral, we obtain
\begin{equation}
\|(A_N-A)u\|_2
\le\int|k(t)|\,
\|(W_{t/N}AW_{t/N}^*-A)u\|_2\,dt.
\label{E70}
\end{equation}
\postdisplaypenalty=10000\relax
For each fixed \(t\in\mathbb R^d\), strong continuity at \(t_0=0\) in \eqref{E112} and \(W_0=I\) imply \(W_{t/N}AW_{t/N}^*u\to Au\) in \(L^2(\mathbb R^n)\). The integrand in \eqref{E70} is bounded by \(2\|A\|\|u\|_2|k(t)|\in L^1(\mathbb R^d)\), independently of \(N\). Dominated convergence in \eqref{E70} now gives \(A_Nu\to Au\) in \(L^2(\mathbb R^n)\) for every \(u\). The proof is complete.
\end{proof}

\begin{remark}[Normalization of the conjugation average]
The average in \eqref{eq:conjugation-average} is the regularization of Georgescu and Iftimovici in our Weyl coordinates. To see the normalization, set \(\Lambda(s,\tau)=(-\sqrt2s,-\sqrt2\tau)\). With \(W_{\rm GI}\) and \(\sigma_{\rm GI}\) denoting their Schr\"odinger representation and symplectic form, \cite[(1.36)--(1.37)]{GeorgescuIftimovici2023} and \eqref{E24}, \eqref{E55} give \(W_{\rm GI}(\Lambda a)=W_a\) and \(\sigma_{\rm GI}(\Lambda a,\Lambda t)=-2\sigma(a,t)\). Under this change of coordinates, their coefficient formula \cite[(3.3)]{GeorgescuIftimovici2023} and \eqref{eq:weyl-map} identify their distribution as \(A^\#=\Lambda_*U\), where \(\langle \Lambda_*U,\Phi\rangle=\langle U,\Phi\circ \Lambda\rangle\) for \(\Phi\in\mathcal S(\mathbb R^d)\).

For the averaging measure, \(|\det \Lambda|=2^n\), so the substitution \(\xi=\Lambda t\) gives \(d\bar\xi=(2\pi)^{-n}d\xi=\pi^{-n}dt\). Define \(\theta\in\mathcal S(\mathbb R^d)\) by \(\theta(\Lambda t)=\pi^n k(t)\) and set \(\theta_\varepsilon(\xi)=\varepsilon^{-d}\theta(\xi/\varepsilon)\) for \(\varepsilon=N^{-1}\). Then
$$
\theta_\varepsilon(\Lambda t)\,\pi^{-n}dt=k_N(t)\,dt,
\qquad \int_{\mathbb R^d}\theta_\varepsilon(\xi)\,d\bar\xi=1,
$$
so the regularization of \(A\) in \cite[(3.4)]{GeorgescuIftimovici2023} equals \(\widetilde A_N\). In their Fourier convention \cite[(3.2)]{GeorgescuIftimovici2023}, alternation of \(\sigma\) and \eqref{E27} identify the multiplier as
$$
\widehat{\theta_\varepsilon}(\Lambda a)
=\int_{\mathbb R^d}e^{-2i\sigma(a,t)}k_N(t)\,dt
=\chi(a/N).
$$
For this normalized Schwartz weight, \cite[Lemma 3.2(iii)]{GeorgescuIftimovici2023} gives \((\widetilde A_N)^\#=\widehat{\theta_\varepsilon}A^\#\). Testing against \(\Phi\in\mathcal S(\mathbb R^d)\) and using \(A^\#=\Lambda_*U\) and \eqref{E35}, we obtain
$$
\langle (\widetilde A_N)^\#,\Phi\rangle
=\langle U,(\widehat{\theta_\varepsilon}\circ \Lambda)(\Phi\circ \Lambda)\rangle
=\langle \Lambda_*(\chi(\cdot/N)U),\Phi\rangle.
$$
\end{remark}

For bounded Weyl operators whose spreading distributions satisfy a common cone condition, we apply \eqref{eq:compact-composition} to compactly supported cutoffs of both distributions. On the distribution side the cone condition gives weak convergence of the twisted products by \eqref{eq:cutoff-convergence}. On the operator side Lemma \ref{lem:bounded-cutoffs} gives strong convergence of the factors with uniform norm bounds, so their compositions converge strongly as well.

\begin{proposition}\label{prop:bounded-composition}

Let \(U,V\in\mathcal C_{e,\eta}\), and suppose that \(W(U)\) and \(W(V)\) have bounded extensions \(A,B\) on \(L^2(\mathbb R^n)\). Then
\begin{equation}
W(U\star V)=AB \label{eq:bounded-composition}
\end{equation}
as maps \(\mathcal S(\mathbb R^n)\to\mathcal S'(\mathbb R^n)\), and the left side has the bounded extension \(AB\).

\end{proposition}

\begin{proof}

Choose \(\chi\) as in Lemma \ref{lem:bounded-cutoffs}, and set \(U_N=\chi(\cdot/N)U\), \(V_N=\chi(\cdot/N)V\), \(A_N=W(U_N)\), and \(B_N=W(V_N)\). By that lemma, \(A_N\) and \(B_N\) have uniformly bounded extensions to \(L^2(\mathbb R^n)\), denoted by the same symbols, and converge strongly to \(A\) and \(B\), respectively. The uniform bound on \(A_N\) allows us to pass to the product, since for every \(u\in L^2(\mathbb R^n)\)
\begin{equation}
\|(A_NB_N-AB)u\|_2
\le\|A_N\|\,\|(B_N-B)u\|_2
+\|(A_N-A)Bu\|_2\longrightarrow0
\label{E71}
\end{equation}
as \(N\to\infty\).

It remains to identify the strong limit \(AB\) in \eqref{E71} with the Weyl map of \(U\star V\). Fix \(u,v\in\mathcal S(\mathbb R^n)\) and set \(c_{u,v}(a)=\langle W_a u,v\rangle_2\). This is a Schwartz test by \eqref{E88}, so we may use it in the weak convergence \eqref{eq:cutoff-convergence}. At each \(N\), compact support of \(U_N,V_N\) ensures that \(A_N,B_N\) preserve \(\mathcal S\) by Lemma \ref{L4}, and \eqref{eq:compact-composition} identifies \(W(U_N\star V_N)\) there with \(A_NB_N\). Consequently \(\langle U_N\star V_N,c_{u,v}\rangle=\langle A_NB_Nu,v\rangle_2\) for every \(N\) by \eqref{eq:weyl-map}. Taking the weak distributional limit on the left and the strong operator limit on the right, we find
$$
\begin{aligned}
\langle W(U\star V)u,\overline v\rangle
&=\langle U\star V,c_{u,v}\rangle
=\lim_{N\to\infty}\langle U_N\star V_N,c_{u,v}\rangle\\
&=\lim_{N\to\infty}\langle A_NB_Nu,v\rangle_2
=\langle ABu,v\rangle_2.
\end{aligned}
$$
Indeed, \(|\langle (A_NB_N-AB)u,v\rangle_2|\le\|(A_NB_N-AB)u\|_2\|v\|_2\to0\) by \eqref{E71}. Since \(u,v\in\mathcal S(\mathbb R^n)\) are arbitrary, the identification \eqref{E44} proves \eqref{eq:bounded-composition} as maps \(\mathcal S\to\mathcal S'\), with \(AB\) as the bounded extension. This completes the proof.
\end{proof}

\begin{corollary}\label{cor:bounded-zero-product}

Let \(U,V\in\mathcal C_{e,\eta}\). If \(W(U)\) and \(W(V)\) have bounded extensions \(A,B\) on \(L^2(\mathbb R^n)\), then \(AB=0\) implies \(A=0\) or \(B=0\).

\end{corollary}

\begin{proof}

If \(AB=0\), Proposition \ref{prop:bounded-composition} gives \(W(U\star V)=0\) on \(\mathcal S(\mathbb R^n)\). Injectivity in \eqref{E40} yields \(U\star V=0\), and Theorem \ref{thm:distribution} applies to the two factors in \(\mathcal C_{e,\eta}\). Hence, \(U=0\) or \(V=0\). The corresponding bounded extension \(A\) or \(B\) vanishes first on the dense subspace \(\mathcal S(\mathbb R^n)\) and then on \(L^2(\mathbb R^n)\) by continuity, which completes the proof.
\end{proof}

\subsection{The parabolic algebra}\label{sec:6}

The Volterra and analytic nests, studied on \(L^2(\mathbb R)\) in \cite{KatavolosPower1997} and on \(L^p(\mathbb R)\), \(1<p<\infty\), in \cite{Kastis2017}, give the support information needed for Corollary \ref{cor:bounded-zero-product}. Invariance of both nests under each \(A\in\mathcal A_p\) places \(\operatorname{supp}U_A\) in a fixed quadrant.

Let \(\mathcal F:L^2(\mathbb R)\to L^2(\mathbb R)\) denote the unitary normalization of the Fourier transform in \eqref{E7}, given by
$$
(\mathcal Fu)(\xi)=(2\pi)^{-1/2}
\int_{\mathbb R}u(x)e^{-ix\xi}\,dx
$$
for \(u\in\mathcal S(\mathbb R)\) and \(\xi\in\mathbb R\).
For \(r\in\mathbb R\), define the closed subspaces $$E_r=\{u\in L^2(\mathbb R):u=0\text{ a.e. on }(-\infty,r)\}$$ and \(\widehat E_r=\mathcal F^{-1}E_r\), so that \(u\in\widehat E_r\) precisely when \(\mathcal Fu\) vanishes almost everywhere on \((-\infty,r)\). With \(\{0\}\) and \(L^2(\mathbb R)\) adjoined to each family, the \(E_r\) form the Volterra nest and the \(\widehat E_r\) form the analytic nest \cite[Section 2]{KastisPower2024}.

A bounded operator \(A\) on \(L^2(\mathbb R)\) preserves a closed subspace \(E\) if \(AE\subset E\). For \(\mu,\lambda\ge0\), the relations \(D_\mu E_r\subset E_{r+\mu}\subset E_r\) and \(M_\lambda E_r=E_r\) show that both generators preserve every \(E_r\). For \(u\in\mathcal S(\mathbb R)\), a change of variable in \eqref{E7} for translations and a frequency shift for modulations give
$$
\begin{aligned}
(\mathcal FD_\mu\mathcal F^{-1}u)(\xi)&=e^{-i\mu\xi}u(\xi),\quad
(\mathcal FM_\lambda\mathcal F^{-1}u)(\xi)&=u(\xi-\lambda)
\end{aligned}
$$
for every \(\xi\in\mathbb R\). Both identities extend to \(L^2(\mathbb R)\) by density, since the operators on either side are unitary. The conjugated operators preserve every \(E_r\), so \(D_\mu\widehat E_r\subset\widehat E_r\) and \(M_\lambda\widehat E_r\subset\widehat E_r\).

To pass from the generators to \(\mathcal A_p\), fix a closed subspace \(E\) and write \(E^\perp\) for its orthogonal complement. Since \(E=(E^\perp)^\perp\), the inclusion \(AE\subset E\) is equivalent to
$$
\langle Au,v\rangle_2=0
$$
for every \(u\in E\) and \(v\in E^\perp\).
For each such pair \(u,v\), the functional \(A\mapsto\langle Au,v\rangle_2\) is continuous in the weak operator topology, so the algebra of operators preserving \(E\) is weak-operator closed. When \(E=E_r\) or \(E=\widehat E_r\), it contains all products \(D_\mu M_\lambda\) with \(\mu,\lambda\ge0\) and hence their span \(\mathcal A_{\mathrm{alg}}\). By \eqref{E1}, every \(A\in\mathcal A_p\) therefore satisfies
\begin{equation*}
A E_r\subset E_r,\qquad A\widehat E_r\subset\widehat E_r
\qquad(r\in\mathbb R).
\tag{\getrefnumber{E71}a}\label{E113}
\end{equation*}

\begin{lemma}\label{lem:quadrant}

For \(A\in\mathcal A_p\), its tempered spreading distribution \(U_A\) satisfies
\begin{equation}
\operatorname{supp}U_A
\subset Q:=\{(s,t)\in\mathbb R^2:s\ge0,\ t\le0\}. \label{eq:quadrant}
\end{equation}

\end{lemma}

\begin{proof}

We first prove a half-plane support criterion for any bounded operator \(B\) on \(L^2(\mathbb R)\) satisfying \(BE_r\subset E_r\) for every \(r\in\mathbb R\). In the kernel coordinates of \eqref{eq:weyl-kernel}, \(s=(x-y)/\sqrt2\), so \(s<0\) exactly when \(x<y\). Fix \(r\) and take \(\varphi\in C_c^\infty((r,\infty))\) and \(\psi\in C_c^\infty((-\infty,r))\), extended by zero to \(\mathbb R\). Since \(B\varphi\in E_r\) and \(\psi\in E_r^\perp\), the kernel pairing \eqref{E78} satisfies
$$
\langle K_B(x,y),\varphi(y)\overline{\psi(x)}\rangle
=\langle B\varphi,\psi\rangle_2=0.
$$
Put \(R_r=(-\infty,r)\times(r,\infty)\). The preceding identity shows that \(K_B\) annihilates every product test on \(R_r\). Their linear span is dense in \(C_c^\infty(R_r)\) \cite[Theorem 39.2]{Treves1967}, so continuity of \(K_B\) implies \(K_B|_{R_r}=0\).

To extend this vanishing to the open half-plane \(\{x<y\}\), fix \(\Phi\in C_c^\infty(\{x<y\})\). Each point \((x,y)\) of its support lies in \(R_{(x+y)/2}\), so compactness provides a finite cover by rectangles \(R_{r_1},\ldots,R_{r_m}\). Choose a smooth partition of unity \(\{\chi_j\}_{j=1}^m\) subordinate to this cover on a neighborhood of \(\operatorname{supp}\Phi\). Then \(\Phi=\sum_{j=1}^m\chi_j\Phi\), with \(\chi_j\Phi\in C_c^\infty(R_{r_j})\). Since \(K_B|_{R_{r_j}}=0\) for every \(j\), linearity gives
$$
\langle K_B,\Phi\rangle
=\sum_{j=1}^m\langle K_B,\chi_j\Phi\rangle=0.
$$
Since \(\Phi\) was arbitrary, \(K_B\) vanishes on \(\{x<y\}\). By Definition \ref{def:distribution-support}, its support satisfies
\begin{equation}
\operatorname{supp}K_B\subset\{(x,y):x\ge y\}.
\label{E72}
\end{equation}

The partial Fourier transform in \eqref{E110} leaves \(s\) unchanged. To carry the support bound \eqref{E72} from \(K_B\) to \(U_B\), we therefore localize in \(s\). Write $$J(s,q)=((q+s)/\sqrt2,(q-s)/\sqrt2)$$ as in \eqref{E110}. For \(\theta\in C_c^\infty((-\infty,0))\) and \(\Phi\in C_c^\infty(\mathbb R^2)\), we define $$(\theta\Phi)(s,q)=\theta(s)\Phi(s,q).$$ Multiplication by \(\theta(s)\) commutes with \(\mathcal F_t\) by \eqref{E53} and \eqref{E54}. Since \(x-y=\sqrt2s\) under \(J\), the test \((\theta\Phi)\circ J^{-1}\) is supported in \(\{x<y\}\). Then \eqref{E110} with \(n=1\) and the support bound \eqref{E72} yield
$$
\begin{aligned}
\langle\mathcal F_t(\theta U_B),\Phi\rangle
&=\langle\mathcal F_tU_B,\theta\Phi\rangle=\sqrt2\,\langle K_B,(\theta\Phi)\circ J^{-1}\rangle=0.
\end{aligned}
$$
By density in \eqref{E18}, \(\mathcal F_t(\theta U_B)=0\) in \(\mathcal S'(\mathbb R^2)\), and injectivity of \(\mathcal F_t\) by \eqref{E74} implies \(\theta U_B=0\). For any \(\varphi\in C_c^\infty(\{s<0\})\), choose \(\theta=1\) near the \(s\)-projection of \(\operatorname{supp}\varphi\). Then \(\langle U_B,\varphi\rangle=\langle\theta U_B,\varphi\rangle=0\), proving the implication
\begin{equation}
BE_r\subset E_r\ \ (r\in\mathbb R)\quad\Longrightarrow\quad\operatorname{supp}U_B\subset\{s\ge0\}.
\label{E29}
\end{equation}

By the first inclusion in \eqref{E113}, \eqref{E29} with \(B=A\) gives \(\operatorname{supp}U_A\subset\{s\ge0\}\). For the other coordinate we use the same criterion with \(B=\mathcal FA\mathcal F^{-1}\). This operator is bounded by unitarity, and the second inclusion in \eqref{E113} implies
\begin{equation}
(\mathcal FA\mathcal F^{-1})E_r
=\mathcal F(A\widehat E_r)
\subset\mathcal F\widehat E_r=E_r
\label{E79}
\end{equation}
for every \(r\in\mathbb R\). Hence \eqref{E29} gives \(\operatorname{supp}U_B\subset\{s\ge0\}\) for \(B=\mathcal FA\mathcal F^{-1}\). To read this bound in the original coordinates, we compute the action of Fourier conjugation on \(W_{(s,t)}\). For \(u\in\mathcal S(\mathbb R)\) and \(\xi\in\mathbb R\), substituting \(x=y+\sqrt2s\) in the Fourier transform of \eqref{E24} gives
$$
\begin{aligned}
(\mathcal FW_{(s,t)}u)(\xi)
&=e^{-ist}e^{-i\sqrt2s\xi}(\mathcal Fu)(\xi+\sqrt2t)=(W_{(-t,s)}\mathcal Fu)(\xi).
\end{aligned}
$$
Since \(\mathcal F\) and the Weyl operators are unitary, this identity extends from \(\mathcal S\) to \(L^2\) and yields
\begin{equation}
\mathcal FW_{(s,t)}\mathcal F^{-1}=W_{(-t,s)}. \label{eq:fourier-rotation}
\end{equation}

For an invertible linear map \(\Lambda:\mathbb R^2\to\mathbb R^2\), let \(\langle \Lambda_*H,\varphi\rangle=\langle H,\varphi\circ \Lambda\rangle\) for \(H\in\mathcal S'(\mathbb R^2)\) and \(\varphi\in\mathcal S(\mathbb R^2)\). Pullback by \(\Lambda\) is a Schwartz space isomorphism by \eqref{E107}, and \(\Lambda_*H\) is tempered by the transpose formula \eqref{E74}. For every open \(O\subset\mathbb R^2\), composition with \(\Lambda\) also maps \(C_c^\infty(\Lambda O)\) bijectively onto \(C_c^\infty(O)\), giving
\begin{equation*}
(\Lambda_*H)|_{\Lambda O}=0\quad\Longleftrightarrow\quad H|_O=0,
\qquad \operatorname{supp}(\Lambda_*H)=\Lambda(\operatorname{supp}H).
\tag{\getrefnumber{eq:fourier-rotation}a}\label{E96}
\end{equation*}
To apply this identity, we identify the distribution obtained under the rotation \(R(s,t)=(-t,s)\). For \(u,v\in\mathcal S(\mathbb R)\), apply \eqref{eq:weyl-map} to \(\mathcal F^{-1}u,\mathcal F^{-1}v\). Using unitarity to move \(\mathcal F\) across the inner product and \eqref{eq:fourier-rotation} to identify the conjugated Weyl operator, we obtain
$$
\begin{aligned}
\langle\mathcal FA\mathcal F^{-1}u,v\rangle_2
&=\big\langle U_A(s,t),\langle\mathcal FW_{(s,t)}\mathcal F^{-1}u,v\rangle_2\big\rangle=\big\langle U_A(s,t),\langle W_{(-t,s)}u,v\rangle_2\big\rangle.
\end{aligned}
$$
By the pushforward definition and \eqref{eq:weyl-map}, the last pairing is \(\langle W(R_*U_A)u,\overline v\rangle\). Since \(u,v\) are arbitrary, we have
\begin{equation}
W(R_*U_A)=\mathcal FA\mathcal F^{-1}
\label{E73}
\end{equation}
as maps \(\mathcal S(\mathbb R)\to\mathcal S'(\mathbb R)\). By uniqueness in Lemma \ref{L3}, \eqref{E73} identifies \(R_*U_A\) with \(U_B\). The support transformation \eqref{E96} now carries the half-plane bound for \(U_B\) back to \(U_A\) through the identity
$$
R(\operatorname{supp}U_A)
=\operatorname{supp}(R_*U_A)
\subset\{(s,t):s\ge0\}.
$$
The first coordinate of \(R(s,t)\) is \(-t\), so \(\operatorname{supp}U_A\subset\{t\le0\}\). Intersecting with the bound \(s\ge0\) already obtained for \(A\) proves \eqref{eq:quadrant}.
\end{proof}

\begin{proof}[Proof of Corollary \ref{cor:parabolic}]

Let \(A,B\in\mathcal A_p\) satisfy \(AB=0\). With \(e=(1,-1)\), each \((s,t)\in Q\) satisfies
$$
e\cdot(s,t)=s-t=s+|t|\ge(s^2+t^2)^{1/2}.
$$
Lemma \ref{lem:quadrant} places both \(\operatorname{supp}U_A\) and \(\operatorname{supp}U_B\) in \(Q\), so the displayed bound verifies \eqref{eq:cone} for both distributions with \(e=(1,-1)\), \(\eta=1\), and \(M_{U_A}=M_{U_B}=0\). By \eqref{E78} and \eqref{E116}, \(A,B\) are the bounded extensions of \(W(U_A),W(U_B)\). Corollary \ref{cor:bounded-zero-product} therefore applies to \(AB=0\), giving \(A=0\) or \(B=0\). The inclusions in \eqref{E1} carry this conclusion to \(\mathcal A_{\mathrm n}\) and \(\mathcal A_{\mathrm{alg}}\). The proof is complete.
\end{proof}

\raggedbottom
\section{Toeplitz zero products on the one-dimensional Fock space}\label{sec:concluding}

To apply Corollary \ref{cor:bounded-zero-product}, we use the Toeplitz--Weyl correspondence \cite{BergerCoburn1994} to relate the Fourier support of \(f\) to the spreading support of the Bargmann conjugate of \(T_f\). With inner products on \(L^2(\mathbb C,d\gamma)\) and its Fock subspace linear in the first variable, the Bargmann transform \cite{Bargmann1961} is the map \(\mathcal B:L^2(\mathbb R)\to\mathcal F^2(\mathbb C)\) given by
$$
(\mathcal Bu)(z)=\pi^{-1/4}\int_{\mathbb R}
e^{-r^2/2+\sqrt2\,rz-z^2/2}u(r)\,dr
$$
for \(u\in L^2(\mathbb R)\) and \(z\in\mathbb C\).

In our normalization the change from Fourier to spreading coordinates is the linear map \(L:\mathbb R^2\to\mathbb R^2\) given by \(L(\xi_1,\xi_2)=(\xi_2/2,-\xi_1/2)\). With the pushforward convention used in the proof of Lemma \ref{lem:quadrant}, we have
\begin{equation}
|L\xi|=|\xi|/2,\qquad
\langle L_*H,\varphi\rangle=\langle H,\varphi\circ L\rangle
\label{E49}
\end{equation}
for \(\xi=(\xi_1,\xi_2)\in\mathbb R^2\), \(H\in\mathcal S'(\mathbb R^2)\), and \(\varphi\in\mathcal S(\mathbb R^2)\). In particular \(L_*H\) is tempered.

Zhu's Bargmann transform \(\mathcal B_1:L^2(\mathbb R)\to\mathcal F^2(\mathbb C)\) at \(\alpha=1\) differs from \(\mathcal B\) by the unitary dilation \((\mathcal Uu)(x)=2^{1/4}u(\sqrt2x)\) on \(L^2(\mathbb R)\). Indeed, taking \(\mathcal Uu\) as input in \cite[(6.9)]{Zhu2012} and substituting \(r=\sqrt2x\) gives the kernel of \(\mathcal B\) above, so
\begin{equation*}
\mathcal B=\mathcal B_1\mathcal U.
\tag{\getrefnumber{E49}a}\label{E114}
\end{equation*}
The unitarity of \(\mathcal B\) now follows from that of \(\mathcal B_1\) in \cite[Theorem 6.8]{Zhu2012}.

We now express the Weyl quantization conjugated by \(\mathcal U\) in the spreading coordinates of \eqref{E49}. In Zhu's normalization at \(\alpha=1\), write \(\mathsf D_0=(2i)^{-1}d/dx\) and \(Xw(x)=xw(x)\) on \(\mathcal S(\mathbb R)\).

\begin{lemma}[Spreading distribution in the dilated Weyl coordinates]\label{L7}
For \(h\in\mathcal S'(\mathbb R^2)\), let \(\sigma_h(\zeta,x)=h(x,-\zeta)\) denote the pullback under \((\zeta,x)\mapsto(x,-\zeta)\). We use the distributional Weyl quantization \(\sigma_h(\mathsf D_0,X):\mathcal S(\mathbb R)\to\mathcal S'(\mathbb R)\), whose action on tests is given in \eqref{E60}. For all \(u,v\in\mathcal S(\mathbb R)\), this quantization satisfies
$$
\big\langle\sigma_h(\mathsf D_0,X)\mathcal Uu,\overline{\mathcal Uv}\big\rangle
=\big\langle W\bigl((4\pi^2)^{-1}L_*\widehat h\bigr)u,\overline v\big\rangle.
$$
\end{lemma}

\begin{proof}
Suppose first that \(h\in\mathcal S(\mathbb R^2)\). In Zhu's Fourier representation, conjugation by \(\mathcal U\) turns each exponential into a Weyl operator. To determine its parameter, use \((\mathcal U^*u)(r)=2^{-1/4}u(r/\sqrt2)\) for \(u\in\mathcal S(\mathbb R)\) and \(r\in\mathbb R\). In Zhu's formula \cite[(1.21)]{Zhu2012} at \(\alpha=1\) for \(e^{2i(p\mathsf D_0+qX)}\), the parameters \(p=-\xi_2/2\) and \(q=\xi_1/2\) yield
$$
\bigl(e^{i(\xi_1X-\xi_2\mathsf D_0)}w\bigr)(x)
=e^{-i\xi_1\xi_2/4}e^{i\xi_1x}w(x-\xi_2/2)
\qquad(w\in\mathcal S(\mathbb R)).
$$
Taking \(w=\mathcal Uu\) and evaluating at \(x=r/\sqrt2\) cancels the two dilation factors. Formula \eqref{E24} then identifies the conjugate as
\begin{equation}
\begin{aligned}
(\mathcal U^*e^{i(\xi_1X-\xi_2\mathsf D_0)}\mathcal Uu)(r)
&=e^{-i\xi_1\xi_2/4}e^{i\xi_1r/\sqrt2}u(r-\xi_2/\sqrt2)\\
&=(W_{L\xi}u)(r).
\end{aligned}
\label{E31}
\end{equation} With \(p,q\in\mathbb R\) now denoting the Fourier variables in Zhu's representation, the substitution \(y=-\zeta\) and our convention \eqref{E7} express the \(2\pi\) Fourier transform of \(\sigma_h\) as
$$
\int_{\mathbb R^2}\sigma_h(\zeta,x)e^{-2\pi i(p\zeta+qx)}\,d\zeta\,dx
=\widehat h(2\pi q,-2\pi p).
$$ In Zhu's Fourier representation \cite[(1.26)]{Zhu2012}, set \(\xi=(2\pi q,-2\pi p)\). Then \(2\pi(p\mathsf D_0+qX)=\xi_1X-\xi_2\mathsf D_0\), and \eqref{E31} identifies the conjugated exponential applied to \(u\) with \(W_{L\xi}u\). By \eqref{E24} the vector \(W_{L\xi}u\) depends continuously on \(\xi\) in \(L^2(\mathbb R)\) and has norm \(\|u\|_2\). Since \(\widehat h\in L^1(\mathbb R^2)\), the Bochner integral of \(\widehat h(\xi)W_{L\xi}u\) exists in \(L^2(\mathbb R)\). For \(u,v\in\mathcal S(\mathbb R)\), the corresponding scalar integral is absolutely convergent, since
$$
\int_{\mathbb R^2}|\widehat h(\xi)|\,
|\langle W_{L\xi}u,v\rangle_2|\,d\xi
\le\|\widehat h\|_1\|u\|_2\|v\|_2<\infty.
$$
We may therefore take the inner product with \(v\) under the integral. With \(dp\,dq=(4\pi^2)^{-1}d\xi\), the resulting coefficient is
\begin{equation}
\begin{aligned}
&\quad\langle\mathcal U^*\sigma_h(\mathsf D_0,X)\mathcal Uu,v\rangle_2\\
&=\int_{\mathbb R^2}\widehat h(2\pi q,-2\pi p)
\langle\mathcal U^*e^{2\pi i(p\mathsf D_0+qX)}\mathcal Uu,v\rangle_2\,dp\,dq\\
&=\frac1{4\pi^2}\int_{\mathbb R^2}\widehat h(\xi)\langle W_{L\xi}u,v\rangle_2\,d\xi=\frac1{4\pi^2}\big\langle\widehat h(\xi),\langle W_{L\xi}u,v\rangle_2\big\rangle.
\end{aligned}
\label{E32}
\end{equation} For the last equality we use \eqref{E44} with the test \(\phi(\xi)=\langle W_{L\xi}u,v\rangle_2\). Indeed, the coefficient in \eqref{E88} is Schwartz, and pullback under the invertible map \(L\) preserves this space by \eqref{E107}.

By \eqref{E15}, the right side of \eqref{E32} tests \(h\) against \((4\pi^2)^{-1}\widehat\phi\). To pass from Schwartz to tempered symbols, we determine the test for the operator on the left and show that it is the same function. For \(u,v\in\mathcal S(\mathbb R)\), put \(u_0=\mathcal Uu\) and \(v_0=\mathcal Uv\), both in \(\mathcal S(\mathbb R)\), and define \(\psi:\mathbb R^2\to\mathbb C\) by
\begin{equation*}
\psi(\zeta,x)=\frac1\pi\int_{\mathbb R}
 u_0(x-r/2)\overline{v_0(x+r/2)}e^{2ir\zeta}\,dr.
\tag{\getrefnumber{E32}a}\label{E125}
\end{equation*}
For \(h\in\mathcal S(\mathbb R^2)\), the absolute integral for the coefficient of Zhu's kernel formula \cite[(1.27)]{Zhu2012} is bounded by \(\pi^{-1}\|u_0\|_1\|v_0\|_1\int_{\mathbb R}\sup_x|\sigma_h(\zeta,x)|\,d\zeta<\infty\). We may therefore apply Fubini's theorem and set \(x=(a+b)/2\) and \(r=a-b\). The inverse substitution \(a=x+r/2\), \(b=x-r/2\) has absolute Jacobian one and yields
\begin{equation}
\begin{aligned}
&\quad\langle\sigma_h(\mathsf D_0,X)u_0,\overline{v_0}\rangle\\
&=\frac1\pi\int_{\mathbb R^3}\sigma_h\bigl(\zeta,\tfrac{a+b}{2}\bigr)
u_0(b)\overline{v_0(a)}e^{2i(a-b)\zeta}\,da\,db\,d\zeta\\
&=\frac1\pi\int_{\mathbb R^3}\sigma_h(\zeta,x)
u_0(x-r/2)\overline{v_0(x+r/2)}e^{2ir\zeta}\,dr\,dx\,d\zeta\\
&=\langle\sigma_h,\psi\rangle.
\end{aligned}
\label{E33}
\end{equation} The last equality uses the integral pairing \eqref{E44}. The test \(\psi\) has the Wigner form in \cite[(1.91)--(1.92)]{Folland1989}. Its Schwartz regularity and its continuous dependence on \(u_0,\overline{v_0}\) follow from the same construction. Start with \(u_0\otimes\overline{v_0}\), which depends continuously on the two factors by \eqref{E118}. Pullback under \((x,r)\mapsto(x-r/2,x+r/2)\) and partial Fourier transformation in \(r\) at \(-2\zeta\) give \(\pi\psi\). Both operations preserve the Schwartz topology by \eqref{E107}. Since the change \((\zeta,x)\mapsto(x,-\zeta)\) has absolute Jacobian one, for \(h\in\mathcal S'(\mathbb R^2)\) we may now define \(\langle\sigma_h(\mathsf D_0,X)u_0,\overline{v_0}\rangle\) by
\begin{equation}
\langle\sigma_h,\psi\rangle
=\langle h(x,y),\psi(-y,x)\rangle.
\label{E60}
\end{equation}
For fixed \(h\), the right side of \eqref{E60} is jointly continuous and bilinear in \(u_0,\overline{v_0}\). Applying \eqref{E129} with \(n=1\) gives a continuous map from \(\mathcal S(\mathbb R)\) into its strong dual, agreeing with \eqref{E33} for Schwartz symbols. To check that this extension is the Weyl quantization of a tempered symbol, write \(\mathsf D=(2\pi i)^{-1}d/dx\) as in Folland, so \(\mathsf D_0=\pi\mathsf D\). Apply \cite[Proposition (2.5)]{Folland1989} to the symbol \((\omega,x)\mapsto\sigma_h(\pi\omega,x)\), understood by pullback. For a Schwartz test \(\Theta(\omega,x)\), this pullback pairs \(\sigma_h(\zeta,x)\) with \(\pi^{-1}\Theta(\zeta/\pi,x)\). Taking \(\Theta\) to be Folland's Wigner function of \(u_0,v_0\), set \(p=-r\) in \cite[(1.91)]{Folland1989}. Its integrand at \((\zeta/\pi,x)\) becomes \(u_0(x-r/2)\overline{v_0(x+r/2)}e^{2ir\zeta}\), so \(\pi^{-1}\Theta(\zeta/\pi,x)=\psi(\zeta,x)\) by \eqref{E125}. This identifies the coefficient in Folland's proposition with \eqref{E60}.

It remains to compare the test \(\psi(-y,x)\) for the quantization in \eqref{E60} with the Fourier test from \eqref{E32}. For fixed \(u,v\), these correspond to the two functionals
$$
h\longmapsto\langle h(x,y),\psi(-y,x)\rangle,
\qquad h\longmapsto\frac1{4\pi^2}\langle h,\widehat\phi\rangle.
$$
For Schwartz \(h\), the \(L^2\)-integral in \eqref{E32} represents the same coefficient as \eqref{E33}. Indeed, using \eqref{E44} for the first equality and unitarity of \(\mathcal U\) for the second, we have
$$
\begin{aligned}
\langle\sigma_h(\mathsf D_0,X)u_0,\overline{v_0}\rangle
&=\langle\sigma_h(\mathsf D_0,X)\mathcal Uu,\mathcal Uv\rangle_2\\
&=\langle\mathcal U^*\sigma_h(\mathsf D_0,X)\mathcal Uu,v\rangle_2.
\end{aligned}
$$
By \eqref{E60}, the left side is \(\langle h(x,y),\psi(-y,x)\rangle\), while \eqref{E32} and the Fourier pairing \eqref{E15} identify the right side with \((4\pi^2)^{-1}\langle h,\widehat\phi\rangle\). These pairings agree for every Schwartz \(h\), in particular for \(h\in C_c^\infty(\mathbb R^2)\). By \eqref{E44}, the difference of the two Schwartz tests represents the zero distribution, so \eqref{E37} gives the pointwise identity
$$
\psi(-y,x)=\frac1{4\pi^2}\widehat\phi(x,y)
\qquad ((x,y)\in\mathbb R^2).
$$
Pairing this identity with any \(h\in\mathcal S'(\mathbb R^2)\) and using \eqref{E60} on the left and \eqref{E15} on the right, we obtain
\begin{equation}
\begin{aligned}
\langle\sigma_h(\mathsf D_0,X)u_0,\overline{v_0}\rangle
&=\langle h(x,y),\psi(-y,x)\rangle=\frac1{4\pi^2}\langle\widehat h,\phi\rangle,
\qquad h\in\mathcal S'(\mathbb R^2).
\end{aligned}
\label{E76}
\end{equation}
By \eqref{E49} and the definition of \(\phi\), the last pairing in \eqref{E76} satisfies
\(\langle\widehat h,\phi\rangle=\langle L_*\widehat h,a\mapsto\langle W_a u,v\rangle_2\rangle\).
The definition \eqref{eq:weyl-map} now identifies the Weyl map in the statement.
\end{proof}

For \(f\in L^\infty(\mathbb C)\), write \(A_f=\mathcal B^*T_f\mathcal B\). Zhu's Toeplitz formula expresses \(A_f\) through the Weyl quantization of a Gaussian average of \(f\), to which Lemma \ref{L7} applies. In spreading coordinates the resulting multiplier is \(\Gamma(a)=e^{-|a|^2/2}\), \(a\in\mathbb R^2\). Its derivatives are bounded, so \eqref{E89} and \eqref{E35} define multiplication by \(\Gamma\) on \(\mathcal S'(\mathbb R^2)\).

\begin{lemma}[Spreading distribution of a Fock Toeplitz operator]\label{L1}

For \(f\in L^\infty(\mathbb C)\), we use the Fourier convention \eqref{E15} and define
\begin{equation}
U_f=\frac{1}{4\pi^2}\Gamma\,L_*\widehat f.
\label{E30}
\end{equation}
Then
\begin{equation}
W(U_f)=\mathcal B^*T_f\mathcal B.
\label{eq:fock-spreading}
\end{equation}
The equality holds as maps \(\mathcal S(\mathbb R)\to\mathcal S'(\mathbb R)\), and the right side is the bounded extension of \(W(U_f)\) to \(L^2(\mathbb R)\).

\end{lemma}

\begin{proof}

Write \(G(z)=(2/\pi)e^{-2|z|^2}\) for \(z\in\mathbb C\), and let \(h_f\) be the Gaussian average \(G*f\). Since \(\|G\|_1=1\) and \(f\in L^\infty(\mathbb C)\), the average is defined at every point. The two bounds needed for Weyl quantization and for Zhu's Toeplitz formula are
$$
\begin{aligned}
&\|h_f\|_\infty\le(2/\pi)\|e^{-2|\cdot|^2}\|_1\|f\|_\infty=\|f\|_\infty,\\
&\int_{\mathbb C}|f(z+w)|^2\,d\gamma(z)\le\|f\|_\infty^2\quad(w\in\mathbb C).
\end{aligned}
$$
The first estimate places \(h_f\) in \(\mathcal S'(\mathbb R^2)\) by \eqref{E44}, and the second verifies condition \((I_2)\) in \cite[Section 3.2]{Zhu2012}. The Gaussian average in \cite[(3.8)]{Zhu2012} is \(B_2 f=h_f\). With \(z=x+i\zeta\), its value at \(\overline z=x-i\zeta\) is \(h_f(x,-\zeta)\), so \cite[Theorem 6.17]{Zhu2012} with \(\alpha=1\) and \(g=f\) gives
\begin{equation}
\mathcal B_1^*T_f\mathcal B_1=\sigma_{h_f}(\mathsf D_0,X),
\qquad \sigma_{h_f}(\zeta,x)=h_f(x,-\zeta)
\label{E75}
\end{equation}
on \(\mathcal S(\mathbb R)\), with \((\zeta,x)\in\mathbb R^2\).
By \eqref{E114}, \(A_f=\mathcal U^*(\mathcal B_1^*T_f\mathcal B_1)\mathcal U\). Evaluating \eqref{E75} on \(\mathcal Uu,\mathcal Uv\) and applying Lemma \ref{L7} to \(h_f\in\mathcal S'(\mathbb R^2)\) therefore give
\begin{equation}
\begin{aligned}
\langle A_fu,v\rangle_2
&=\big\langle\sigma_{h_f}(\mathsf D_0,X)\mathcal Uu,\overline{\mathcal Uv}\big\rangle\\
&=\big\langle W\bigl((4\pi^2)^{-1}L_*\widehat h_f\bigr)u,\overline v\big\rangle
\end{aligned}
\label{E77}
\end{equation}
for all \(u,v\in\mathcal S(\mathbb R)\).

To identify the spreading distribution in \eqref{E77} with \(U_f\), we compute \(L_*\widehat h_f\). For \(\tau\in\mathbb R\), the change \(x=\sqrt{\pi/2}\,y\) and \cite[Lemma 6.11]{Zhu2012} with parameter \(\tau/\sqrt{8\pi}\) give
$$
\int_{\mathbb R}e^{-2x^2-i\tau x}\,dx
=\sqrt{\pi/2}\,e^{-\tau^2/8}.
$$
Our Fourier convention \eqref{E7} and the product of the two one-dimensional integrals give \(\widehat G(\xi)=e^{-|\xi|^2/8}\). To identify the average \(h_f\) with the distributional convolution \(G*f\), take \(\varphi\in\mathcal S(\mathbb R^2)\). Since
$$
\iint_{\mathbb R^2\times\mathbb R^2}|G(x-y)f(y)\varphi(x)|\,dy\,dx
\le \|G\|_1\|f\|_\infty\|\varphi\|_1,
$$
Fubini's theorem in \eqref{E44} and the evenness of \(G\) yield \(\langle h_f,\varphi\rangle=\langle f,G*\varphi\rangle\). This is the action of the distributional convolution \(f*G\) on \(\varphi\). In \cite[Example 30.1 and Theorem 30.4, (30.7)]{Treves1967}, take \(S=f\in\mathcal S'\) and \(T=G\), a rapidly decreasing distribution since \(G\) is Schwartz. To express the convolution identity in our normalization, note that \(G*\widehat\varphi=\widehat{\widehat G\,\varphi}\) by \eqref{E7}, Fubini's theorem and the evenness of \(G\). Thus \eqref{E15} gives \(\langle\widehat h_f,\varphi\rangle=\langle f,G*\widehat\varphi\rangle=\langle\widehat f,\widehat G\,\varphi\rangle\). With the Gaussian transform computed above, this reads
$$
\widehat h_f=e^{-|\cdot|^2/8}\widehat f
=(\Gamma\circ L)\widehat f
\quad\text{in }\mathcal S'(\mathbb R^2),
$$
where \((\Gamma\circ L)(\xi)=e^{-|\xi|^2/8}\) because \(|L\xi|=|\xi|/2\) by \eqref{E49}.
For \(\varphi\in\mathcal S(\mathbb R^2)\), the definition of \(L_*\) in \eqref{E49} tests this identity on \(\varphi\circ L\). Since \((\Gamma\circ L)(\varphi\circ L)=(\Gamma\varphi)\circ L\), the multiplication rule \eqref{E35} yields
$$
\langle L_*\widehat h_f,\varphi\rangle
=\langle\widehat f,(\Gamma\varphi)\circ L\rangle
=\langle\Gamma L_*\widehat f,\varphi\rangle.
$$
Since \(\varphi\) is arbitrary, \((4\pi^2)^{-1}L_*\widehat h_f=U_f\) by \eqref{E30}. Substitution in \eqref{E77} proves \eqref{eq:fock-spreading} as maps \(\mathcal S(\mathbb R)\to\mathcal S'(\mathbb R)\), since the coefficients agree for every \(u,v\in\mathcal S(\mathbb R)\). Finally, \(\|A_f\|=\|T_f\|\) by unitarity of \(\mathcal B\), and \(T_f=P(f\,\cdot)\) with \(P\) an orthogonal projection gives \(\|T_f\|\le\|f\|_\infty\). Hence \(A_f\) is the asserted bounded extension to \(L^2(\mathbb R)\).
\end{proof}

\begin{proof}[Proof of Corollary \ref{cor:fock-zero-product}]
We first relate \(\operatorname{supp}U_f\) to \(\operatorname{supp}\widehat f\). This relation will detect vanishing of the symbol and transfer its Fourier support condition to the spreading distribution. For \(H\in\mathcal S'(\mathbb R^2)\) and \(\varphi\in C_c^\infty(\mathbb R^2)\), the nonvanishing of \(\Gamma\) gives \(\varphi/\Gamma\in C_c^\infty(\mathbb R^2)\) with the same support as \(\varphi\). By \eqref{E35},
$$
\langle H,\varphi\rangle
=\left\langle\Gamma H,\frac{\varphi}{\Gamma}\right\rangle.
$$
For every open \(O\subset\mathbb R^2\), multiplication by \(\Gamma\) maps \(C_c^\infty(O)\) onto itself, with inverse \(\varphi\mapsto\varphi/\Gamma\). Thus \(H|_O=0\) if and only if \((\Gamma H)|_O=0\). By Definition \ref{def:distribution-support}, the two distributions have the same support. Taking \(H=L_*\widehat f\) in \eqref{E30} and using \eqref{E96} for the invertible map \(L\), we obtain
\begin{equation}
\operatorname{supp}U_f=L(\operatorname{supp}\widehat f).
\label{eq:fock-support}
\end{equation}
Since \(L\) is invertible, \eqref{eq:fock-support} makes emptiness of \(\operatorname{supp}U_f\) equivalent to emptiness of \(\operatorname{supp}\widehat f\). A tempered distribution with empty support vanishes on \(C_c^\infty\), hence on \(\mathcal S\) by \eqref{E18}. Thus \(U_f=0\) is equivalent to \(\widehat f=0\), or, by Fourier inversion \eqref{E13}, to vanishing of \(f\) as a distribution.
The functional \(\varphi\mapsto\int f\varphi\) is continuous on \(L^1(\mathbb R^2)\). If \(f\) vanishes as a distribution, \eqref{E44} makes this functional zero on \(C_c^\infty(\mathbb R^2)\), which is dense in \(L^1(\mathbb R^2)\). Testing on \(\overline f\,\mathbf1_{\{|x|\le R\}}\in L^1\) now yields \(\int_{|x|\le R}|f(x)|^2\,dx=0\) for every \(R>0\). Conversely, an almost everywhere zero symbol defines the zero distribution by \eqref{E44}. We have proved
\begin{equation}
U_f=0\quad\Longleftrightarrow\quad f=0\ \text{almost everywhere}.
\label{E2}
\end{equation}

We transfer the common cone condition from \(\widehat f,\widehat g\) to \(U_f,U_g\) by setting \(e'=L^{-T}e\), where \(L^{-T}\) denotes the transpose of \(L^{-1}\). For \(\xi\in\operatorname{supp}\widehat f\), we use \(|\xi|=2|L\xi|\) from \eqref{E49} to write the cone bound \eqref{eq:cone} as
$$
e'\cdot L\xi=e\cdot\xi
\ge\eta|\xi|-M_{\widehat f}
=2\eta|L\xi|-M_{\widehat f}.
$$
The same calculation applies to \(\widehat g\), with its own support constant. Formula \eqref{eq:fock-support} places \(U_f,U_g\) in \(\mathcal C_{e',2\eta}\), and \eqref{eq:fock-spreading} identifies \(A_f,A_g\) as bounded extensions of \(W(U_f),W(U_g)\). If \(T_fT_g=0\), unitarity of \(\mathcal B\) gives
$$
A_fA_g=\mathcal B^*T_fT_g\mathcal B=0.
$$
Corollary \ref{cor:bounded-zero-product} implies \(A_f=0\) or \(A_g=0\). By \eqref{eq:fock-spreading} and injectivity in \eqref{E40}, the corresponding spreading distribution vanishes. By \eqref{E2}, it follows that \(f=0\) or \(g=0\) almost everywhere.
\end{proof}

For bounded symbols on the Fock space over \(\mathbb C^n\), Bauer and Le obtain the zero-product conclusion when one symbol is radial \cite[Corollary 3.6]{BauerLe2011}. Bommier-Hato proves a corresponding result on generalized Fock spaces under the same radiality assumption \cite[Theorem B]{BommierHato2020}. In one dimension, Bauer, Choe and Koo treat harmonic symbols whose holomorphic and antiholomorphic parts are of exponential type, assuming that the product vanishes in both orders \cite[Theorem 4.2]{BauerChoeKoo2015}.

In complex dimension \(n\ge2\), Qin constructed two bounded nonradial symbols whose nonzero Toeplitz operators have zero product \cite[Theorem 1.1]{Qin2026}. In one dimension, removing the Fourier support hypothesis from Corollary \ref{cor:fock-zero-product} leads to the bounded-symbol case of Bauer and Le's Question A \cite[Section 6]{BauerLe2011}.

\begin{question}\label{Q2}
If \(f,g\in L^\infty(\mathbb C)\) and \(T_fT_g=0\) on \(\mathcal F^2(\mathbb C)\), must one of \(f,g\) vanish almost everywhere?
\end{question}

\end{document}